# Partial Dynamical Systems and the KMS Condition

Ruy Exel* and Marcelo Laca**

Departamento de Matemática; Universidade Federal de Santa Catarina; 88040-900 Florianópolis SC; Brazil, and
Department of Mathematics; University of Newcastle; NSW 2308 – Australia.

ABSTRACT. Given a countably infinite 0–1 matrix $A$ without identically zero rows, let $\mathcal{O}_A$ be the Cuntz–Krieger algebra recently introduced by the authors and $\mathcal{T}_A$ be the Toeplitz extension of $\mathcal{O}_A$, once the latter is seen as a Cuntz–Pimsner algebra, as recently shown by Szymański. We study the KMS equilibrium states of $C^*$-dynamical systems based on $\mathcal{O}_A$ and $\mathcal{T}_A$, with dynamics satisfying $\sigma_t(s_x) = N_x^{it} s_x$ for the canonical generating partial isometries $s_x$ and arbitrary real numbers $N_x > 1$. The KMS$_\beta$ states on both $\mathcal{O}_A$ and $\mathcal{T}_A$ are completely characterized for certain values of the inverse temperature $\beta$, according to the position of $\beta$ relative to three critical values, defined to be the abscissa of convergence of certain Dirichlet series associated to $A$ and the $N(x)$. Our results for $\mathcal{O}_A$ are derived from those for $\mathcal{T}_A$ by virtue of the former being a covariant quotient of the latter. When the matrix $A$ is finite, these results give theorems of Olesen and Pedersen for $\mathcal{O}_n$ and of Enomoto, Fujii and Watatani for $\mathcal{O}_A$ as particular cases.

CONTENTS

Our motivation for the present work stems from the perception that the Cuntz-Krieger algebras for infinite matrices studied in [**EL**] naturally provide $C^*$-dynamical systems with interesting KMS state structure. The main phenomena in which we are interested are those intrinsically associated to the infinite dimensionality of the matrices but our approach also gives fresh insight into some salient features that have not been observed or emphasized enough even in the finite dimensional case, particularly with respect to the consideration of nonperiodic dynamics and the symmetries of the equilibrium states.

*Date: June 21, 2000.*

* Partially supported by CNPq.

** Supported by the Australian Research Council.

The earliest ancestor of our results is an intriguing theorem of Olesen and Pedersen's [**OP**], which appeared in the late seventies amidst a flurry of examples and counterexamples triggered by the advent of the Cuntz algebras $\mathcal{O}_n$ [**C**], stating that the periodic gauge action on $\mathcal{O}_n$ admits a unique KMS (equilibrium) state, whose inverse temperature is $\beta = \log n$. Shortly afterwards Cuntz and Krieger came up with their $C^*$-algebras $\mathcal{O}_A$ [**CK**] and, in the ensuing flurry, the theorem of Olesen and Pedersen was duly generalized by Enomoto, Fujii, and Watatani [**EFW**], who proved, among other things, that when the matrix $A$ is irreducible (and not a permutation) the gauge action on $\mathcal{O}_A$ admits a unique KMS state, at inverse temperature equal to the logarithm of the spectral radius of $A$. Despite the explicitness of the computations involved those early results have up to until recently been largely regarded as curious counterexamples rather than as sources of interesting new phenomena to be explored. Part of the reason for this derives from [**OP**: Theorem 1], according to which the dynamics involved are "nonphysical" because they have no (weak, approximate) Hamiltonian. Recently, however, the interest in KMS states of Cuntz-Krieger algebras has been renewed, mainly in [**PWY**] where results along the lines of [**EFW**] have been obtained for periodic full dynamics on unital $C^*$-algebras, and where the KMS condition is linked to a variational principle for the entropy.

The main purpose and the methods of the present work are of a different nature: we aim to study the KMS equilibrium states of $C^*$-dynamical systems that are inspired on the periodic gauge action of $\mathbb{R}$ on $\mathcal{O}_A$, but which are more general in three important aspects:

- we allow the matrix $A$ to be countably infinite as in [**EL**];
- we focus on the Toeplitz extension $\mathcal{T}_A$ rather than on $\mathcal{O}_A$ itself; and,
- while still dealing with dynamics having the generating partial isometries as eigenvectors, we allow the possibility of different eigenvalues and thus of nonperiodicity.

In order to deal effectively with the new situation we must first spend some effort in developing the necessary approach and technical tools, and this involves realizing our $C^*$-algebras as crossed products by partial actions of a countably generated free group, and characterizing KMS states of such crossed products in terms of a certain invariance property of probability measures under the partial action, using techniques analogous to those of [**L**] for semigroup actions. There is no significant extra cost in carrying out this first task in the slightly more general context of $C^*$-algebras that are topologically graded over free groups, and we do so in the first few sections.

We then specialize to our main setting, which we would like to describe briefly next. Given a matrix $A$ of zeros and ones over a countable set $\mathcal{G}$, we consider certain one-parameter groups of gauge automorphisms of three closely related $C^*$-algebras, namely $\mathcal{TO}_A$, $\mathcal{T}_A$, and $\mathcal{O}_A$.

$\mathcal{O}_A$ is the generalized Cuntz–Krieger algebra introduced by the authors in [**EL**] for an arbitrary infinite 0–1 matrix $A$. $\mathcal{TO}_A$ was also introduced in [**EL**] as an auxiliary tool to study $\mathcal{O}_A$, and $\mathcal{T}_A$ is the Toeplitz extension of $\mathcal{O}_A$, once the latter is seen as a Cuntz–Pimsner algebra as shown by Szymański [**Sz**].

All of the above three algebras have canonical generating sets consisting of partial isometries, say $\{s_x\}_{x\in\mathcal{G}}$ and, given a choice of positive real numbers $\{N(x)\}_{x\in\mathcal{G}}$, there are one-parameter groups of gauge automorphisms satisfying

$$\sigma_t(s_x) = N(x)^{it} s_x, \quad t \in \mathbb{R}.$$

Clearly these are subgroups of the canonical gauge action of the torus $\mathbb{T}^{\mathcal{G}}$. Since KMS states are known to be $\sigma$-invariant, one expects them to be invariant under the (compact) closure of $\{\sigma_t\}$ inside this torus, and thus to factor through the conditional expectation onto the fixed point algebra. One of the biggest surprises we find here is that the KMS states under analysis are shown (Theorem 8.2) to factor through the conditional expectation onto a much smaller subalgebra which can be identified as the fixed point algebra for a coaction of the infinitely generated free group. We do not explore this coaction here except for the fact that it leads to a highly useful conditional expectation. Nevertheless it is remarkable that the KMS condition seems to impose the preservation of symmetries way beyond what is expected at first.

The small subalgebra mentioned above is actually a commutative algebra and hence the search for KMS states boils down to a study of measures on its spectrum. Specially when dealing with the case of $\mathcal{T}_A$, to which we dedicate the biggest share of our attention, there is a natural dichotomy breaking the spectrum into a "finite" part, denoted $\Omega_f$, and an "infinite" part $\Omega_\infty$ (see 8.6). The measures considered are therefore classified in *finite type* or *infinite type* according to whether they assign full mass to the finite or to the infinite part of the spectrum.

The behavior of KMS states at inverse temperature $\beta$, and hence also of the measures which determine them, strongly depends on the relative position of $\beta$ with respect to three critical inverse temperatures. In order to describe these let us denote by $P_A$ the set of all admissible words (with respect to the given matrix $A$) in the alphabet $\mathcal{G}$, by $P_A^y$ the admissible words ending in $y$, and by $P_A^{xy}$ the admissible words beginning in $x$ and ending in $y$. We then introduce three (families of) Dirichlet series of one variable $\beta$, namely

$$Z(\beta) = \sum_{\mu \in P_A} N(\mu)^{-\beta}, \quad Z_y(\beta) = \sum_{\mu \in P_A^y} N(\mu)^{-\beta}, \quad \text{and} \quad Z_{xy}(\beta) = \sum_{\mu \in P_A^{xy}} N(\mu)^{-\beta},$$

where $N(\mu)$ is defined by $N(\mu_1)\cdots N(\mu_k)$ when $\mu$ is the admissible word $\mu = \mu_1 \cdots \mu_k$.

Every Dirichlet series has an abscissa of convergence which marks the lower end of its interval of convergence. Accordingly we denote by $\beta_c$ the abscissa of convergence of $Z(\beta)$ and this turns out to be the first important critical inverse temperature. We prove that all $\text{KMS}_\beta$ states correspond to finite type measures for $\beta$ above this critical point. In addition we are able to describe these measures in very concrete terms and hence all $\text{KMS}_\beta$ states are concretely exhibited (see 9.7).

In the case of an irreducible matrix $A$ the Dirichlet series $Z_y$, for $y \in \mathcal{G}$, satisfy a "solidarity" property in the sense that, for a given $\beta$, either they all converge or they all diverge. Therefore there is a single abscissa of convergence, denoted $\dot\beta_c$, which does not depend on $y$. Another solidarity property holds among the $Z_{xy}$, in turn defining a third critical value $\ddot\beta_c$. Since each of $Z(\beta)$, $Z_y(\beta)$, and $Z_{xy}(\beta)$ is a subseries of the previous one it is clear that their abscissa of convergence satisfy $\ddot\beta_c \leq \dot\beta_c \leq \beta_c$.

Still speaking of the irreducible case we prove (Theorem 10.6) that all $\text{KMS}_\beta$ states are of infinite type for $\beta$ below $\dot\beta_c$ and that there are no $\text{KMS}_\beta$ states at all for $\beta$ below $\ddot\beta_c$ (Theorem 14.5).
The following diagram illustrates these results:

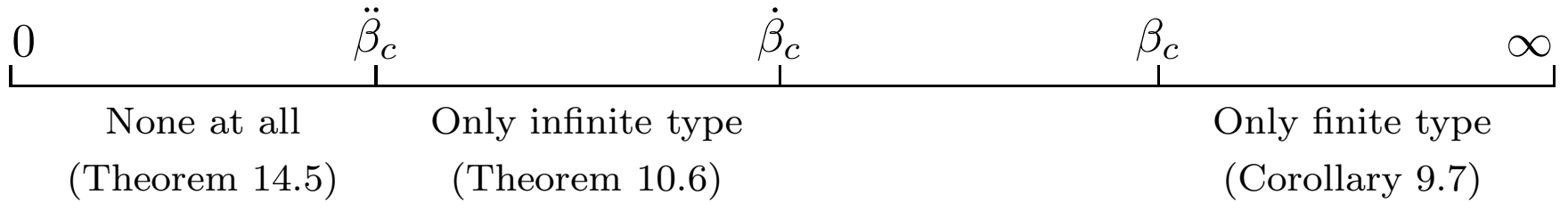

*KMS states and critical inverse temperatures*

The results sketched in this diagram are the strongest results we can offer under the sole assumption that $A$ is irreducible but there are several strengthenings we can provide under extra hypotheses. For instance, we show (Theorem 15.2) that there are no $\beta$-scaling states at all for $\beta < \dot\beta_c$ under the hypothesis that $\inf_{x\in\mathcal{G}} N(x) > 1$.

Unfortunately there is nothing we can say about the interval between $\dot\beta_c$ and $\beta_c$, which remains conspicuously absent from our conclusions. Not being able to deal with it we may at least identify a rather common situation in which it collapses (Proposition 15.4), and this is when there is a finite target set, i.e. a finite set $\{y_1, \ldots, y_n\} \subseteq \mathcal{G}$ such that for every $x \in \mathcal{G}$ one has $A(x, y_i) = 1$ for at least one $i$ (see 10.1.(FTS)). A related problem which we could not resolve is whether or not finite type KMS states can coexist with infinite type ones in the case of an irreducible matrix. We therefore leave these as open problems.

Nevertheless under all of the hypotheses mentioned so far our theory gives a complete description of KMS states for all inverse temperatures, except for the critical inverse temperature $\beta_c$. At $\beta_c$ quite different things can happen. There are examples in which there is a single $\text{KMS}_{\beta_c}$ state (Theorem 18.4) but there are also examples in which infinitely many such states exist (Section 16).

Even though our main focus is on $\mathcal{T}_A$ we can provide some useful information about the KMS states on $\mathcal{O}_A$ as well. Since $\mathcal{O}_A$ is a covariant quotient of $\mathcal{T}_A$, the set of KMS states on $\mathcal{O}_A$ correspond to the set of KMS states on $\mathcal{T}_A$ which factor through $\mathcal{O}_A$. We therefore take up the problem of characterizing the latter set (Theorem 17.2) giving several equivalent ways to describe it.

For a finite matrix $A$ one may obviously say a lot more than in the general case. As it turns out we give a complete characterization of all KMS states on $\mathcal{T}_A$ for a finite irreducible matrix (Theorem 18.4). With respect to the KMS states on $\mathcal{O}_A$ for a finite $A$, we completely characterize its KMS states even if $A$ is not irreducible (Theorem 18.5).

In particular when considering the gauge action, i.e. when $N(x) = e$ for all $x$, our methods can be easily applied to recover the result of Olesen and Pedersen [**OP**] on the uniqueness of the KMS state on $\mathcal{O}_n$ as well as the result obtained by Enomoto, Fujii, and Watatani [**EFW**] for the gauge action on $\mathcal{O}_A$.

In the course of the research reported here we were deeply influenced by Vere–Jones paper [**V**] in which he generalizes the classical Perron–Frobenius theorem to the case of infinite matrices. But, because of the difference between our emphasis on Dirichlet series and Vere-Jones's emphasis on power series, among other reasons, we were often impeded to use his results in a straightforward way.

The present work culminates a project started in December 1996 and continued through several short visits of R.E. to Newcastle and of M.L. to Florianópolis, and we would both like to thank the members of both departments for the hospitality provided to the visitor of turn. We also gratefully acknowledge funding from CNPq and from the Australian Research Council.

## 1. KMS states for graded $C^*$-algebras.

Let $B$ be a $C^*$-algebra and let $G$ be a discrete group. We shall say $B$ is *graded over* $G$, or simply $G$*-graded*, if we are given a linearly independent family $\{B_g\}_{g\in G}$ of closed linear subspaces of $B$ such that, for all $g, h \in G$,

(i) $B_g B_h \subseteq B_{gh}$,
(ii) $B_g^* = B_{g^{-1}}$, and
(iii) $\bigoplus_{g\in G} B_g$ is dense in $B$.

The main examples are given by crossed products algebras, including the case of partial actions which will be discussed in some detail below.

We shall fix throughout a $G$-graded $C^*$-algebra $B$. Moreover we will fix a strongly continuous one-parameter group $\sigma = \{\sigma_t\}_{t\in\mathbb{R}}$ of automorphisms of $B$ such that each $\sigma_t$ restricted to each $B_g$ is a multiple of the identity. One necessarily has that for each $g$ in $G$ there is a positive real number $N(g)$ such that

$$\sigma_t(b) = N(g)^{it} b, \quad t \in \mathbb{R},\, b \in B_g. \tag{1.1}$$

We shall also assume that $N$ is a group homomorphism from $G$ to the multiplicative group $\mathbb{R}_+^*$ of positive reals (this is in fact necessarily the case if $B_g B_h \neq \{0\}$ for all $g$ and $h$).

Obviously $\sigma$ is determined by $N$. However under the present hypothesis it is not clear whether there exists a one-parameter group $\sigma$ satisfying 1.1 for each group homomorphism $N : G \to \mathbb{R}_+^*$. This existence question will be dealt with when we discuss the case of crossed-products by partial actions below.

Nevertheless, even though partial actions are among our main examples, we stress that we are only assuming, for the time being, that $B$ is $G$-graded and that $\sigma$ satisfies 1.1 for some group homomorphism $N : G \to \mathbb{R}_+^*$.

Recall that an element $b \in B$ is said to be *entire analytic* (with respect to $\sigma$) [**BR**: 2.5.20] if the map $t \in \mathbb{R} \mapsto \sigma_t(b)$ extends to an entire analytic function on the complex plane. For simplicity we shall refer to entire analytic elements simply as *analytic* elements.

Also recall from [**BR**: 5.3.1] that a state $\psi$ on $B$ is said to be a *$\sigma$-KMS state at value* $\beta \in (0, \infty)$ — the *inverse temperature* in Mathematical Physics terminology — or simply a *KMS$_\beta$ state* if for any pair of elements $a$ and $b$ in a given norm-dense $\sigma$-invariant *-subalgebra of analytic elements of $B$ one has

$$\psi(a\sigma_{i\beta}(b)) = \psi(ba). \tag{1.2}$$

By the last sentence of [**BR**: 5.3.7] one has that, for a KMS state, the identity above in fact holds for every $a$ in $B$ and every analytic element $b$ in $B$.

In the special case $\beta = \infty$, KMS$_\beta$ states are called *ground states* and are defined in a slightly different fashion but by [**BR**: 5.3.19] they are the states on $B$ such that for every pair of analytic elements $a$ and $b$ in $B$,

$$\sup_{\mathrm{Im} z \geq 0} |\psi(a\sigma_z(b))| < \infty. \tag{1.3}$$

In order to verify a state to be a ground state the most economical way is to verify 1.3 only for $a$ and $b$ in a norm-dense subspace of analytic elements of $B$ (*cf.* [**Pe**: 8.12.3] and [**L**: Remark 11]).

Given $g \in G$ it is clear from 1.1 that each $b$ in $B_g$ is analytic and moreover

$$\sigma_z(b) = N(g)^{iz} b, \tag{1.4}$$

for $z \in \mathbb{C}$. It follows by linearity that the algebraic direct sum $\bigoplus_{g\in G} B_g$ consists of analytic elements. In addition the latter is clearly a norm-dense $\sigma$-invariant *-subalgebra of $B$. We will therefore use this algebra whenever we need to verify the KMS condition, both for finite and infinite values of $\beta$.

The following is a characterization of $\sigma$-KMS states on $B$.

**1.5. Proposition.** *Suppose that $B$ is $G$-graded and that $\sigma$ satisfies 1.1 for a group homomorphism $N : G \to \mathbb{R}_+^*$. Let $\psi$ be a state on $B$ and let $\beta \in (0,\infty)$. Then*

(i) *$\psi$ is a KMS$_\beta$ state if and only if $\psi(ab) = N(g)^\beta \psi(ba)$ whenever $a \in B$, $g \in G$, and $b \in B_g$.*

(ii) *$\psi$ is a ground state if and only if $\psi(BB_g) = \{0\}$ whenever $N(g) < 1$.*

*Proof.* It will be convenient to keep in mind that plugging $z = i\beta$ in 1.4 gives

$$\sigma_{i\beta}(b) = N(g)^{-\beta} b.$$

Suppose that $\psi$ is a KMS$_\beta$ state and let $a \in B$ and $b \in B_g$. We then have

$$\psi(ba) = \psi(a\sigma_{i\beta}(b)) = N(g)^{-\beta}\psi(ab),$$

proving the forward implication in (i).

Conversely, suppose that $\psi$ satisfies the equality in (i). Given $a \in B$ and $b \in B_g$ we then have $\psi(ba) = N(g)^{-\beta}\psi(ab) = \psi(a\sigma_{i\beta}(b))$, proving that 1.2 holds for our choice of $a$ and $b$. By linearity we see that the same is true for any $b \in \bigoplus_{g\in G} B_g$. Since the latter is a norm-dense $\sigma$-invariant *-subalgebra of analytic elements we see that $\psi$ is a KMS$_\beta$ state.

Let us now deal with (ii). Given $a \in B$ and $b \in B_g$ we have

$$|\psi(a\sigma_z(b))| = |N(g)^{iz}\psi(ab)| = N(g)^{-\mathrm{Im}z}|\psi(ab)|. \qquad (\dagger)$$

Observe that this is bounded on the upper half plane as a function of $z$ if and only if either $N(g) \geq 1$ or $\psi(ab) = 0$. Thus if $\psi$ is a ground state and $N(g) < 1$ we must have $\psi(ab) = 0$, proving the forward implication in (ii).

Conversely, if $\psi(BB_g) = 0$ whenever $N(g) < 1$ then ($\dagger$) is always bounded on Im$z \geq 0$. It follows by linearity that 1.3 holds for any $b \in \bigoplus_{g\in G} B_g$. Since the latter is a norm-dense set of analytic elements we see that $\psi$ is a ground state. ☐

It should be noted that 1.5.i implies that, for $\beta < \infty$, a KMS$_\beta$ state restricted to $B_e$ must be a trace. In contrast, ground states need not restrict to traces on $B_e$. In fact, if $\sigma$ is the trivial action of $\mathbb{R}$ on $B$, corresponding to $N(g) \equiv 1$, then any state on $B$ is a ground state and one can clearly manufacture examples in which $\psi|_{B_e}$ is not a trace.

Recall from [**E2**: 3.4] that $B$ is said to be *topologically $G$-graded* if there exists a positive contractive conditional expectation

$$E : B \to B_e$$

which vanishes on every $B_g$ for $g \neq e$. From [**E2**: 3.5] it follows that $\bigoplus_{g\in G} B_g$ is a *topological direct sum* in the sense that the canonical projections onto the $B_g$'s extend to bounded linear maps

$$E_g : B \to B_g.$$

If $B$ is topologically $G$-graded and we are given a state $\phi$ on $B_e$ the composition $\psi := \phi \circ E$ is a state on $B$. Our next result is intended to discuss the conditions under which $\psi$ is a $\sigma$-KMS state on $B$.

**1.6. Proposition.** *Assume that $B$ is topologically $G$-graded with conditional expectation $E$ and that $\sigma$ satisfies 1.1 for a group homomorphism $N : G \to \mathbb{R}_+^*$. Let $\phi$ be a state on $B_e$ and set $\psi = \phi \circ E$. Also let $\beta \in (0,\infty)$. Then*

(i) *$\psi$ is a KMS$_\beta$ state if and only if $\phi(ab) = N(g)^\beta \phi(ba)$ for every $g \in G$, $a \in B_{g^{-1}}$, and $b \in B_g$.*

(ii) *$\psi$ is a ground state if and only if $\phi(B_{g^{-1}}B_g) = \{0\}$ whenever $N(g) < 1$.*

*Proof.* The forward implication in (i) follows immediately from 1.5.i. Conversely let $a \in B$ and $b \in B_g$. By considering first the case in which $a \in \bigoplus_{g\in G} B_g$ it is easy to see that

$$E(ab) = E_{g^{-1}}(a)b, \quad \text{and} \quad E(ba) = bE_{g^{-1}}(a).$$

Therefore

$$\psi(ab) = \phi(E(ab)) = \phi(E_{g^{-1}}(a)b) = N(g)^\beta \phi(bE_{g^{-1}}(a)) = N(g)^\beta \phi(E(ba)) = N(g)^\beta \psi(ba).$$

That $\psi$ is a KMS$_\beta$ state now follows from 1.5.i.

The forward implication in (ii) follows immediately from 1.5.ii. Conversely let $a \in B$ and $b \in B_g$ with $N(g) < 1$. We then have

$$\psi(ab) = \phi(E(ab)) = \phi(E_{g^{-1}}(a)b) = 0.$$

Thus $\psi$ is a ground state by 1.5.ii. □

## 2. Graded algebras given by partial actions.

One of the main sources of examples of topologically $G$-graded algebras is the theory of partial actions. Accordingly this section is devoted to reviewing the results of this theory which are relevant to us. The reader should consult [**E1**], [**M**], [**E3**], and [**E4**] for more details.

Recall that a partial action of a discrete group $G$ on a $C^*$-algebra $A$ is a pair

$$\Theta = (\{D_g\}_{g\in G}, \{\theta_g\}_{g\in G})$$

such that, for each $g$ in $G$, $D_g$ is a closed two sided ideal of $A$,

$$\theta_g : D_{g^{-1}} \to D_g$$

is a *-isomorphism, and for all $g$ and $h$ in $G$ one has

(i) $D_e = A$ and $\theta_e$ is the identity automorphism of $A$,
(ii) $\theta_g(D_{g^{-1}} \cap D_h) = D_g \cap D_{gh}$, and
(iii) $\theta_g(\theta_h(a)) = \theta_{gh}(a)$ for all $a \in D_{h^{-1}} \cap D_{h^{-1}g^{-1}}$.

Throughout this section we shall fix a partial action as above. Our goal is to construct a $G$-graded algebra from this data.

Let $\mathcal{L}$ denote the set of all functions $f : G \to A$ such that $f(g) \in D_g$ for all $g \in G$ and moreover $\sum_{g\in G} \|f(g)\| < \infty$. Clearly $\mathcal{L}$ is a Banach space under the norm

$$\|f\| = \sum_{g\in G} \|f(g)\|, \quad f \in \mathcal{L}.$$

It is often convenient to denote by $\sum_{g\in G} a_g\delta_g$ the function $f$ such that $f(g) = a_g$. Employing this notation we define a Banach *-algebra structure on $\mathcal{L}$ by means of the operations

$$(a\delta_g) * (b\delta_h) = \theta_g\left(\theta_g^{-1}(a)b\right)\delta_{gh}, \quad \text{and} \tag{2.1}$$

$$(a\delta_g)^* = \theta_g^{-1}(a^*)\delta_{g^{-1}}$$

for $a \in D_g$ and $b \in D_h$. See the references given above for the proof that $\mathcal{L}$ is indeed a Banach *-algebra under these operations.

**2.2. Definition.** The *crossed-product of the $C^*$-algebra $A$ by the group $G$ under the partial action $\Theta$* is the enveloping $C^*$-algebra of the Banach *-algebra $\mathcal{L}$ described above. We denote it as $A\rtimes_\Theta G$, or simply $A\rtimes G$ if $\Theta$ is understood.

**2.3. Proposition.** *The collection of subspaces $\{B_g\}_{g\in G}$ of $A\rtimes G$ given by $B_g = D_g\delta_g$ makes $A\rtimes G$ into a topologically $G$-graded algebra.*

*Proof.* The construction of $A\rtimes G$ is precisely that of the cross-sectional $C^*$-algebra of the Fell bundle formed by the $B_g$'s under the operations defined by 2.1. The result is then a consequence of [**E2**: 3.2 and 2.9]. □

Observe that $A$ is canonically isomorphic to $B_e$ via the map $a \mapsto a\delta_e$. We will therefore identify $A$ and $B_e$ and hence think of the conditional expectation as the unique bounded map $E : A\rtimes G \to A$ such that

$$E\left(\sum_{g\in G} a_g\delta_g\right) = a_e.$$

Let us now deal with the question of defining the dynamics on $A\rtimes G$ in terms of a given group homomorphism $N$.

**2.4. Proposition.** *Let $N : G \to \mathbb{R}_+^*$ be a group homomorphism. Then there exists a strongly continuous one-parameter group $\sigma$ of *-automorphisms of $A\rtimes G$ such that*

$$\sigma_t(b) = N(g)^{it} b$$

*for all $t \in \mathbb{R}$, $g \in G$, and $b \in B_g$.*

*Proof.* For each $t$ in $\mathbb{R}$ consider the linear operator $\sigma_t$ on $\mathcal{L}$ given by

$$\sigma_t\left(\sum_{g\in G} a_g\delta_g\right) = \sum_{g\in G} N(g)^{it} a_g\delta_g.$$

It is easy to see that each $\sigma_t$ is an isometric *-isomorphism of $\mathcal{L}$. It is also clear that $\sigma_t\sigma_s = \sigma_{t+s}$ for all $t$ and $s$ in $\mathbb{R}$, so that $\sigma$ gives a group homomorphism

$$\sigma : \mathbb{R} \to \mathrm{Aut}(\mathcal{L}),$$

which is obviously strongly continuous. The result now follows by extending each $\sigma_t$ to a *-isomorphism of the enveloping $C^*$-algebra $A\rtimes G$. □

This puts us in the context of section 1 and so we may apply 1.6 to characterize the KMS states on $A\rtimes G$ that factor through the conditional expectation $E$ described above. The following result facilitates checking conditions (i) and (ii) of 1.6:

**2.5. Proposition.** *Let $\Theta$ be a partial action of the discrete group $G$ on a $C^*$-algebra $A$ and consider the standard grading $\{B_g\}_g = \{D_g\delta_g\}_g$ of $A\rtimes G$. Let $N : G \to \mathbb{R}_+^*$ be a group homomorphism.*

(i) *If $\beta \in (0, \infty)$ and $g \in G$ then*

$$\phi(ab) = N(g)^\beta \phi(ba), \quad \forall a \in B_{g^{-1}},\ b \in B_g \qquad (\dagger)$$

*if and only if $\phi$ is a trace and*

$$\phi(\theta_g(a)) = N(g)^{-\beta} \phi(a), \quad \forall a \in D_{g^{-1}}. \qquad (\ddagger)$$

(ii) $B_{g^{-1}}B_g = D_{g^{-1}}$.

*Proof.* We begin with the forward implication in (i). That $\phi$ is a trace follows from (†) applied to the case $g = e$. Given $a \in D_{g^{-1}}$ choose an approximate identity $\{u_i\}_i$ for $D_g$ and observe that

$$\phi(\theta_g(a)) = \lim_i \phi(u_i\theta_g(a)) = \lim_i \phi((u_i\delta_g)\,(a\delta_{g^{-1}})) =$$

$$= N(g)^{-\beta} \lim_i \phi((a\delta_{g^{-1}})\,(u_i\delta_g)) = N(g)^{-\beta} \lim_i \phi(a\theta_{g^{-1}}(u_i)) = N(g)^{-\beta}\phi(a).$$

Conversely, suppose that $\phi$ is a trace satisfying (‡). Given $a \in D_{g^{-1}}$ and $b \in D_g$ we have that $a\delta_{g^{-1}} \in B_{g^{-1}}$ and $b\delta_g \in B_g$ and

$$\phi\big((a\delta_{g^{-1}})(b\delta_g)\big) = \phi\big(\theta_{g^{-1}}(\theta_g(a)b)\big) = N(g)^{\beta}\phi\big(\theta_g(a)b\big) =$$

$$= N(g)^{\beta}\phi\big(b\,\theta_g(a)\big) = N(g)^{\beta}\phi\big((b\delta_g)(a\delta_{g^{-1}})\big).$$

Since $a\delta_{g^{-1}}$ and $b\delta_g$ are generic elements of $B_{g^{-1}}$ and $B_g$, respectively, we have proven (†). We leave the proof of (ii) to the reader. ☐

Observe that 2.5.i.(‡) says that, on the suitable domain, $\phi$ is scaled when composed with $\theta_g$. If one considers a global action (i.e. a partial action for which each $D_g = A$, as in the classical situation) then this scaling property cannot hold in non-trivial situations because the norm of $\phi$ is necessarily invariant. Nevertheless if one deals with *partial actions* this obstruction disappears allowing for many interesting examples as we shall see below. See also [**L**] for examples arising as semigroup crossed products.

## 3. Algebras graded over the free group.

We will be mostly interested in the case where the group $G$ is the free group $\mathbb{F}$ on a (possibly infinite) set $\mathcal{G}$ of generators. When speaking of $\mathbb{F}$ we will often employ the following standard notations:

- If $g \in \mathbb{F}$ we will denote by $|g|$ the length of $g$, that is, the number of generators appearing in the reduced decomposition of $g$.
- $\mathbb{F}_+$ will refer to the positive cone of $\mathbb{F}$, that is, the subsemigroup of $\mathbb{F}$ generated by $\mathcal{G}$, including the unit group element.
- We will usually denote the elements of $\mathbb{F}$ by $g, h, k$; the elements of $\mathcal{G}$ by $x, y, z$; and the elements of $\mathbb{F}_+$ by $\mu, \nu$.

Gradings over $\mathbb{F}$ occasionally satisfy two additional properties which we would now like to recall from [**E2**].

**3.1. Definition.** A grading $\{B_g\}_{g\in\mathbb{F}}$ of a $C^*$-algebra $B$ is said to be:

(i) *semi-saturated* if $B_{gh} = B_g B_h$ (closed linear span) whenever $g, h \in \mathbb{F}$ are such that $|gh| = |g| + |h|$.

(ii) *orthogonal* if $B_x^* B_y = \{0\}$ whenever $x, y \in \mathcal{G}$ and $x \neq y$.

The following Lemma is the main result of this section and is the key to our characterization of KMS states on Cuntz–Krieger algebras.

**3.2. Lemma.** *Let $B$ be a $C^*$-algebra which is $\mathbb{F}$-graded by means of a semi-saturated orthogonal grading $\{B_g\}_{g\in\mathbb{F}}$. Also let $\sigma$ be a strongly continuous one-parameter group of automorphisms of $B$ satisfying 1.1 for some group homomorphism $N : \mathbb{F} \to \mathbb{R}_+^*$. Suppose that $N(x) > 1$ for all $x \in \mathcal{G}$ and let $\psi$ be a KMS$_\beta$ state on $B$ for $\beta$ in the interval $(0, \infty]$. Then $\psi(B_g) = \{0\}$ for all $g \in \mathbb{F}$ with $g \neq e$.*

*Proof.* Initially observe that the hypothesis on $N$ implies that $N(\mu) > 1$ for all $\mu \in \mathbb{F}_+ \setminus \{e\}$.

Since there exists a slight asymmetry between the cases of finite and infinite $\beta$ let us first assume that $\beta \in (0, \infty)$. In this case we claim that $\psi(B_\mu) = \{0\}$ for every $\mu \in \mathbb{F}_+ \setminus \{e\}$. To see this note that $B_\mu = B_e B_\mu$ (closed linear span) by [**BMS**: 1.7], so it suffices to show that $\psi(ab) = 0$ whenever $a \in B_e$ and $b \in B_\mu$. We then have, using 1.5.i twice, that

$$\psi(ab) = N(\mu)^\beta \psi(ba) = N(\mu)^\beta N(e)^\beta \psi(ab) = N(\mu)^\beta \psi(ab).$$

But, since $N(\mu) \neq 1$ and $\beta \neq 0$, we have that $\psi(ab) = 0$ as claimed. Clearly we also have that

$$\psi(B_{\mu^{-1}}) = \psi(B_\mu^*) = \overline{\psi(B_\mu)} = \{0\},$$

so that we have proven the statement for all $g \in \mathbb{F}_+ \cup \mathbb{F}_+^{-1} \setminus \{e\}$.

Since our grading is semi-saturated and orthogonal we have that $B_g = \{0\}$ for every $g \in \mathbb{F}$ which is not of the form $g = \mu\nu^{-1}$, where $\mu, \nu \in \mathbb{F}_+$ and $|g| = |\mu| + |\nu|$ as a moment's reflection will easily show. We therefore assume that $g$ is of the above form.

We will proceed by induction on $m = \min\{|\mu|, |\nu|\}$. If $m = 0$ then either $g = \mu$ or $g = \nu^{-1}$ and the conclusion follows as above. So assume that $m \geq 1$ and write $\mu = x\mu'$ and $\nu = y\nu'$, where $x, y \in \mathcal{G}$ and $\mu', \nu' \in \mathbb{F}_+$. Since the grading is semi-saturated we have that $B_g = B_\mu B_\nu^* = B_x B_{\mu'} B_{\nu'}^* B_y^*$ and hence we will be done once we prove that $\psi(b_x b_{\mu'} b_{\nu'}^* b_y^*) = 0$ whenever $b_i \in B_i$ for $i = x, y, \mu', \nu'$. We have by 1.5.i that

$$\psi(b_x b_{\mu'} b_{\nu'}^* b_y^*) = N(y^{-1}) \psi(b_y^* b_x b_{\mu'} b_{\nu'}^*). \tag{‡}$$

If on the one hand $x \neq y$ then $b_y^* b_x = 0$ by orthogonality and (‡) vanishes. If on the other hand $x = y$ then

$$b_y^* b_x b_{\mu'} b_{\nu'}^* \in B_y^* B_x B_{\mu'} B_{\nu'}^* \subseteq B_{\mu'} B_{\nu'}^*.$$

By the induction hypothesis $\psi$ vanishes on $B_{\mu'} B_{\nu'}^*$ and so again (‡) vanishes. This concludes the proof of the case $\beta < \infty$.

Assume now that $\beta = \infty$ and hence that $\psi$ is a ground state. As before we only need to prove that $\psi(B_g) = 0$ for a nontrivial $g$ of the form $\mu\nu^{-1}$ with $\mu, \nu \in \mathbb{F}_+$ and $|g| = |\mu| + |\nu|$.

If $\nu \neq e$ then $N(\nu) < 1$ and $B_g = B_\mu B_{\nu^{-1}}$ by semi-saturatedness. So we have that $\psi(B_g) = \psi(B_\mu B_{\nu^{-1}}) = \{0\}$ by 1.5.ii.

If $\nu = e$ then $\mu \neq e$ and $\psi(B_g) = \overline{\psi(B_\mu^*)} = \overline{\psi(B_{\mu^{-1}})} = \{0\}$ as seen above. □

Observe that in the above proof, when considering finite values of $\beta$, we only needed that $N(\mu) \neq 1$ for $\mu \in \mathbb{F}_+ \setminus \{e\}$ and hence the result does hold under this weakened hypothesis. However we have not seen how to reach the same conclusion for ground states.

In the topologically graded case we may use 1.6 and 3.2 to get the following very precise characterization of KMS states:

**3.3. Theorem.** *Let $B$ be a $C^*$-algebra which is topologically $\mathbb{F}$-graded by means of a semi-saturated orthogonal grading $\{B_g\}_{g \in \mathbb{F}}$. Also let $\sigma$ be a strongly continuous one-parameter group of automorphisms of $B$ satisfying 1.1 for a given group homomorphism $N : \mathbb{F} \to \mathbb{R}_+^*$. Suppose that $N(x) > 1$ for all $x \in \mathcal{G}$ and let $\beta \in (0, \infty]$. Then the formula*

$$\phi \mapsto \phi \circ E$$

*defines an affine homeomorphism from the set $\mathcal{S}$ defined below onto the set of KMS$_\beta$ states on $B$.*

(i) *If $\beta \in (0, \infty)$ then $\mathcal{S}$ is the set of traces $\phi$ on $B_e$ such that*

$$\phi(ab) = N(x)^\beta \phi(ba), \quad \forall x \in \mathcal{G},\ a \in B_{x^{-1}},\ b \in B_x.$$

(ii) *If $\beta = \infty$ then $\mathcal{S}$ is the set of states $\phi$ of $B_e$ such that*

$$\phi(B_x B_{x^{-1}}) = \{0\}, \quad \forall x \in \mathcal{G}.$$

*Proof.* Under the hypothesis that $\phi$ lies in the set $\mathcal{S}$ described in (i) we claim (see 1.6.i) that

$$\phi(ab) = N(g)^\beta \phi(ba), \quad \forall g \in \mathbb{F},\ a \in B_{g^{-1}},\ b \in B_g.$$

Clearly this holds for $g = e$ because $\phi$ is a trace. Next consider the case $g = x^{-1} \in \mathcal{G}^{-1}$. So let $a \in B_{g^{-1}}$ and $b \in B_g$. Then $b \in B_{x^{-1}}$ and $a \in B_x$ so that $\phi(ba) = N(x)^\beta \phi(ab)$ from where we obtain $\phi(ab) = N(x^{-1})^\beta \phi(ba) = N(g)^\beta \phi(ba)$.

Proceeding by induction on $|g|$ let $|g| > 1$ and write $g = hx$ with $x \in \mathcal{G} \cup \mathcal{G}^{-1}$ and $|g| = |h| + 1$. Given $a \in B_{g^{-1}}$ and $b \in B_g$ suppose that $b = b_h b_x$ where $b_h \in B_h$ and $b_x \in B_x$. Then $ab_h \in B_{x^{-1}}$ and

$$\phi(ab) = \phi(ab_h b_x) = N(x)^\beta \phi(b_x a b_h) = \ldots$$

Moreover $b_x a \in B_{h^{-1}}$ and hence, by the induction hypothesis, the above equals

$$\ldots = N(x)^\beta N(h)^\beta \phi(b_h b_x a) = N(g)^\beta \phi(ba).$$

Since the grading is semi-saturated the linear combinations of the $b$'s considered is dense in $B_g$ and hence the claim is proven. By 1.6.i $\phi$ is a KMS$_\beta$ state.

On the other hand, under the hypothesis that $\phi$ satisfies the property described in (ii) we claim (see 1.6.ii) that

$$\phi(B_{g^{-1}} B_g) = \{0\}, \quad \forall g \in G \text{ such that } N(g) < 1.$$

We have already mentioned that, as a consequence of semi-saturatedness, $B_g = \{0\}$ unless $g = \mu\nu^{-1}$, where $\mu, \nu \in \mathbb{F}_+$ and $|g| = |\mu| + |\nu|$. We therefore suppose that $g$ has this form. Since $N(g) = N(\mu)N(\nu)^{-1} < 1$ and $N(\mu) \geq 1$ we must have $\nu \neq e$. So write $\nu = x\nu'$ with $x \in \mathcal{G}$ and $\nu' \in \mathbb{F}_+$. By semi-saturatedness we have

$$B_{g^{-1}} B_g = B_\nu B_{\mu^{-1}} B_\mu B_{\nu^{-1}} = B_x B_{\nu'} B_{\mu^{-1}} B_\mu B_{\nu'^{-1}} B_{x^{-1}} \subseteq B_x B_e B_{x^{-1}} \subseteq B_x B_{x^{-1}},$$

whence $\phi(B_{g^{-1}} B_g) = \{0\}$. By 1.6.ii $\phi$ is a ground state.

Since $\phi = (\phi \circ E)\,|_{B_e}$ the correspondence is 1–1. In order to show surjectiveness let $\psi$ be a KMS$_\beta$ state on $B$. Then, letting $\phi = \psi|_{B_e}$, we have that $\psi = \phi \circ E$, which is easily proven by checking on $\bigoplus_{g \in G} B_g$ with the help of 3.2. Finally we have that $\phi$ is in $\mathcal{S}$ by the forward implications in 1.6.

Clearly $\phi \mapsto \phi \circ E$ is an affine map. Moreover this is obviously a continuous map (with respect to the weak topologies as usual). Since $\mathcal{S}$ is compact we actually have a homeomorphism. ◻

## 4. Partial actions of the free group.

In this section we will show how to put together the results of the previous sections in order to determine the KMS states on the crossed product $C^*$-algebra resulting from a partial action of the free group with suitable properties. We therefore fix throughout this section a $C^*$-algebra $A$ and a partial action

$$\Theta = (\{D_g\}_{g \in \mathbb{F}}, \{\theta_g\}_{g \in \mathbb{F}})$$

of the free group $\mathbb{F}$ on a possibly infinite set of generators $\mathcal{G}$. The crossed-product algebra $A \rtimes \mathbb{F}$ is therefore topologically graded via the subspaces $B_g = D_g \delta_g$ by 2.3. We begin by determining conditions for this grading to be semi-saturated and orthogonal.

**4.1. Proposition.** *The above grading of $A \rtimes \mathbb{F}$ is:*

(i) *semi-saturated if and only if $D_{gh} \subseteq D_g$ whenever $|gh| = |g| + |h|$.*

(ii) *orthogonal if and only if $D_x \cap D_y = \{0\}$ for all $x, y \in \mathcal{G}$ with $x \neq y$.*

*Assume, in addition, that $A$ is abelian with spectrum a locally compact space $X$ and that $\Theta$ is obtained by means of a partial action*

$$\alpha = \big(\{U_g\}_{g \in \mathbb{F}}, \{\alpha_g\}_{g \in \mathbb{F}}\big)$$

*of $\mathbb{F}$ on $X$ (cf.* [**EL**: Section 2], [**E4**]*). Then the grading of $A \rtimes \mathbb{F}$ is:*

*(a) semi-saturated if and only if $U_{gh} \subseteq U_g$ whenever $|gh| = |g| + |h|$.*

*(b) orthogonal if and only if $U_x \cap U_y = \emptyset$ for all $x, y \in \mathcal{G}$ with $x \neq y$.*

*Proof.* Given $g, h \in \mathbb{F}$ we have

$$B_g B_h = (D_g \delta_g)(D_h \delta_h) = \theta_g\big(\theta_g^{-1}(D_g) D_h\big)\delta_{gh} = \theta_g\big(D_{g^{-1}} D_h\big)\delta_{gh},$$

so that our grading is semi-saturated if and only if $\theta_g\big(D_{g^{-1}} \cap D_h\big) = D_{gh}$ when $|gh| = |g| + |h|$. However it is an axiom for partial actions that $\theta_g(D_{g^{-1}} \cap D_h) = D_g \cap D_{gh}$. So semi-saturatedness is equivalent to $D_g \cap D_{gh} = D_{gh}$ which is the same as saying that $D_{gh} \subseteq D_g$.

If we plug $g = x^{-1}$ and $h = y$ in the equation displayed above, where $x, y \in \mathcal{G}$ and $x \neq y$, then we see that our grading is orthogonal if and only if $D_x \cap D_y = \{0\}$.

As for the part in which $A$ is assumed abelian note that $D_g = C_0(U_g)$ and also that $D_g \cap D_h = C_0(U_g \cap U_h)$ for all $g$ and $h$. Therefore

$$D_{gh} \subseteq D_g \Longleftrightarrow U_{gh} \subseteq U_g$$

and

$$D_x \cap D_y = \{0\} \Longleftrightarrow U_x \cap U_y = \{0\}. \qquad \square$$

It is perhaps worth mentioning that the condition $|gh| = |g| + |h|$ in the free group means that $g \leq gh$ in the sense that $g$ is an initial segment in the reduced decomposition of $gh$. Therefore condition 4.1.i above can be rephrased as saying that $D_g$ is decreasing as a function of $g$ in the sense that $g \leq k \Rightarrow D_k \subseteq D_g$.

**4.2. Definition.** We shall say that $\Theta$ is a *semi-saturated partial action* if condition 4.1.i above holds. We shall say that $\Theta$ is an *orthogonal partial action* when 4.1.ii is satisfied.

Having understood the relationship between graded algebras and partial actions regarding semi-saturatedness and orthogonality we may now present our main abstract result.

**4.3. Theorem.** *Let $\Theta$ be a semi-saturated orthogonal partial action of the free group $\mathbb{F}$ on a $C^*$-algebra $A$. On $A \rtimes \mathbb{F}$ consider the standard grading $\{B_g\}_{g\in\mathbb{F}}$ and conditional expectation $E : A \rtimes \mathbb{F} \to A$. Also let $\{N(x)\}_{x\in\mathcal{G}}$ be a collection of real numbers in the interval $(1, \infty)$. Then there exists a unique strongly continuous one-parameter group $\sigma$ of automorphisms of $A \rtimes \mathbb{F}$ such that*

$$\sigma_t(b) = N(x)^{it} b, \quad \forall x \in \mathcal{G},\ b \in B_x.$$

*The $\sigma$-KMS states at inverse temperature $\beta$ on $A \rtimes \mathbb{F}$ are precisely those of the form $\psi = \phi \circ E$ where $\phi$ is a state on $A$ satisfying:*

(i) *if $\beta < \infty$ : $\phi$ is a trace and $\phi(\theta_x(a)) = N(x)^{-\beta}\phi(a)$ for all $x \in \mathcal{G}$ and all $a \in D_{x^{-1}}$.*

(ii) *if $\beta = \infty$ : $\phi(D_x) = \{0\}$ for all $x \in \mathcal{G}$.*

*Proof.* Extend $N$ to a group homomorphism $N : \mathbb{F} \to \mathbb{R}_+^*$ and hence we may apply 2.4 to conclude that $\sigma$ exists as required. Given that the action is semi-saturated, so is the grading by 4.1.i. With this it is easy to see that $\sigma$ is uniquely determined by the identity displayed in the statement.

We may now apply 3.3 and thus all we need to do is show that the set $\mathcal{S}$ described in 3.3.i–ii can be alternatively characterized by conditions (i–ii) here, but this is precisely the purpose of 2.5. $\square$

## 5. Toeplitz–Cuntz–Krieger algebras for infinite matrices.

Throughout this and the remaining sections of this work we shall fix a countable (meaning finite or countably infinite) set $\mathcal{G}$ and a matrix $A = \{A(x,y)\}_{x,y\in\mathcal{G}}$ with entries in the set $\{0,1\}$, having no identically zero rows. It is in terms of $A$ that we will introduce the algebras which will be the object of our main applications.

Recall from [**EL**] that the Toeplitz–Cuntz–Krieger algebra[1] $\widetilde{\mathcal{TO}}_A$ is the universal unital $C^*$-algebra generated by a family of partial isometries $\{s_x\}_{x\in\mathcal{G}}$ subject to the requirement that their initial projections $q_x = s_x^* s_x$ and final projections $p_y = s_y s_y^*$ satisfy the following conditions for all $x$ and $y$ in $\mathcal{G}$:

CK$_1$) $q_x q_y = q_y q_x$,

CK$_2$) $s_x^* s_y = 0$, if $x \neq y$,

CK$_3$) $q_x s_y = A(x,y) s_y$.

According to [**EL**: 4.6] $\widetilde{\mathcal{TO}}_A$ is the crossed-product $C^*$-algebra for a partial group action which we would now like to briefly describe. See [**EL**] for full details.

**5.1. Definition.** Let[2] $\widetilde{\Omega}_{\mathcal{TO}_A}$ be the closed subset of the compact topological space $2^{\mathbb{F}}$ given by

$$\widetilde{\Omega}_{\mathcal{TO}_A} = \left\{ \begin{array}{ll} \xi \in 2^{\mathbb{F}}: & e \in \xi,\ \xi \text{ is convex,} \\ & \text{if } g \in \xi \text{ there is at most one } x \in \mathcal{G} \text{ such that } gx \in \xi, \text{ and} \\ & \text{if } g \in \xi,\ y \in \mathcal{G}, \text{ and } gy \in \xi \text{ then } gx^{-1} \in \xi \Leftrightarrow A(x,y) = 1 \end{array} \right\}.$$

For each $g \in \mathbb{F}$ let $\Delta_g$ be the clopen subset of $\widetilde{\Omega}_{\mathcal{TO}_A}$ given by

$$\Delta_g = \{\xi \in \widetilde{\Omega}_{\mathcal{TO}_A} : g \in \xi\}$$

and consider the homeomorphism

$$\alpha_g : \xi \in \Delta_{g^{-1}} \longmapsto g\xi \in \Delta_g.$$

Then

$$\alpha := \big(\{\Delta_g\}_{g\in\mathbb{F}}, \{\alpha_g\}_{g\in\mathbb{F}}\big)$$

is a partial group action of $\mathbb{F}$ on $\widetilde{\Omega}_{\mathcal{TO}_A}$ in the sense of [**EL**: Section 2] (see also [**M**], [**E3**], and [**E4**]). This induces a partial action of $\mathbb{F}$ on $C(\widetilde{\Omega}_{\mathcal{TO}_A})$, namely

$$\Theta = \left(\{D_g\}_{g\in G}, \{\theta_g\}_{g\in G}\right),$$

given by $D_g = C_0(\Delta_g)$ and

$$\theta_g : f \in D_{g^{-1}} \mapsto f \circ \alpha_{g^{-1}} \in D_g.$$

The already mentioned Theorem 4.6 of [**EL**] asserts that $\widetilde{\mathcal{TO}}_A$ is isomorphic to $C(\widetilde{\Omega}_{\mathcal{TO}_A}) \rtimes \mathbb{F}$ under an isomorphism that maps each $s_x$ to $1_{\Delta_x}\delta_x$, where $1_{\Delta_x}$ is the characteristic function of $\Delta_x$.

In order to define the next two algebras which are relevant to our study it is convenient to introduce the following notation: given finite subsets $X, Y \subseteq \mathcal{G}$ we let

$$A(X,Y,z) = \prod_{x\in X} A(x,z) \prod_{y\in Y} (1 - A(y,z)), \quad z \in \mathcal{G},$$

and

$$q(X,Y) = \prod_{x\in X} q_x \prod_{y\in Y} (1 - q_y).$$

[1] This algebra was denoted $\mathcal{TO}_A$ in [**EL**] but will be denoted $\widetilde{\mathcal{TO}}_A$ here for reasons which will become apparent soon.

[2] This space was denoted $\Omega_A^\tau$ in [**EL**] but will be denoted $\widetilde{\Omega}_{\mathcal{TO}_A}$ here.

Observe that the 0–1 vector $\big(A(X,Y,z)\big)_{z\in\mathcal{G}}$ is simply the coordinatewise product of the row-vectors of $A$ indexed by $X$, and the *boolean negation* of the row-vectors indexed by $Y$.

Recall from [**EL**] that the (unital) Cuntz–Krieger algebra $\widetilde{\mathcal{O}}_A$ is the quotient of $\widetilde{\mathcal{TO}}_A$ obtained by imposing the following extra relation in addition to $\mathrm{CK}_{1-3}$:

$\mathrm{CK}_4$) $q(X,Y) = \sum_{z\in\mathcal{G}} A(X,Y,z)p_z$ whenever $X,Y \subseteq \mathcal{G}$ are finite sets such that $A(X,Y,z)$ is finitely supported as a function of $z$.

As explained in the first section of [**EL**], condition $\mathrm{CK}_4$ is formally derived from multiplying together sufficiently many occurrences of the the Cuntz–Krieger relations [**CK**] so that the infinite sums involved become finite.

Theorem 7.10 of [**EL**] asserts that $\widetilde{\mathcal{O}}_A \simeq C(\widetilde{\Omega}_{\mathcal{O}_A})\rtimes\mathbb{F}$, where $\widetilde{\Omega}_{\mathcal{O}_A}$ is the $\alpha$-invariant subset of $\widetilde{\Omega}_{\mathcal{TO}_A}$ obtained by taking the closure of the set of unbounded elements (*cf.* Definition 5.5 in [**EL**]). As before there is an isomorphism which maps each $s_x$ to $1_{\Delta_x\cap\widetilde{\Omega}_{\mathcal{O}_A}}\delta_x$.

We shall be concerned here with yet another $C^*$-algebra, denoted $\widetilde{\mathcal{T}}_A$ (see also [**FLR**], [**Sz**]), which sits in between $\widetilde{\mathcal{TO}}_A$ and $\widetilde{\mathcal{O}}_A$ in the sense that the quotient map $\widetilde{\mathcal{TO}}_A \to \widetilde{\mathcal{O}}_A$ alluded to above factors through $\widetilde{\mathcal{T}}_A$.

**5.2. Definition.** Given a countable set $\mathcal{G}$ and a 0–1 matrix $A = \{A(x,y)\}_{x,y\in\mathcal{G}}$ with no identically zero rows we denote by $\widetilde{\mathcal{T}}_A$ the universal unital $C^*$-algebra generated by a family of partial isometries $\{s_x\}_{x\in\mathcal{G}}$ subject to conditions $\mathrm{CK}_{1-3}$ and

$\mathrm{CK}_4^0$) $q(X,Y) = 0$ whenever $X$ and $Y$ are finite subsets of $\mathcal{G}$ such that $A(X,Y,z)$ is identically zero as a function of $z$.

In [**Sz**: Theorem 5] Szymański has realized $\widetilde{\mathcal{O}}_A$ as the Cuntz–Pimsner algebra [**Pi**] of a bimodule, and has shown that $\widetilde{\mathcal{T}}_A$ is the corresponding Toeplitz extension [**Sz**: Theorem 10] (the case in which $A$ is the edge matrix of a graph had been dealt with in [**FLR**]). Because of this and other reasons related to the interesting features of KMS states on $\widetilde{\mathcal{T}}_A$, we will gradually concentrate our attention on $\widetilde{\mathcal{T}}_A$.

Observe that $\mathrm{CK}_4^0$, seen as a set of relations, is a subset of $\mathrm{CK}_4$ since the identity in the latter reduces to $q(X,Y) = 0$ when $A(X,Y,\,\cdot\,) \equiv 0$. Given that $\mathrm{CK}_4^0$ is less restrictive than $\mathrm{CK}_4$ we therefore have that the algebras $\widetilde{\mathcal{TO}}_A$, $\widetilde{\mathcal{T}}_A$, and $\widetilde{\mathcal{O}}_A$ are "decreasing" in the sense that each is a quotient of the preceding one.

We would now like to describe $\widetilde{\mathcal{T}}_A$ as the crossed-product algebra for a partial group action in order to be able to study its KMS states using the tools developed in the previous sections. In preparation for this we need to recall some terminology from [**EL**]. Given $\xi \in \widetilde{\Omega}_{\mathcal{TO}_A}$ and $g \in \xi$ recall from [**EL**: 5.5 and 5.6] that

- the *root of $g$ relative to $\xi$* is the subset of $\mathcal{G}$ defined by $R_\xi(g) = \{x \in \mathcal{G} : gx^{-1} \in \xi\}$.
- the *stem* of an element $\xi \in \widetilde{\Omega}_{\mathcal{TO}_A}$ is the unique maximal (finite or infinite) word in the alphabet $\mathcal{G}$ such that all of its finite initial subwords (interpreted as elements of $\mathbb{F}_+$) belong to $\xi$.
- $\xi$ is said to be *bounded* if its stem is finite. Otherwise $\xi$ is said to be *unbounded*.

We shall also make use of the topological space $\Sigma_A$ obtained by taking the closure, within the compact space $2^{\mathcal{G}}$, of the set of columns of $A$. Observe that a column of $A$ is a 0–1 vector and hence it may indeed be seen as an element of $2^{\mathcal{G}}$. Moreover, any subset of $\mathcal{G}$ such as $R_\xi(g)$ can, and often will, also be interpreted as belonging to $2^{\mathcal{G}}$ in the usual way.

**5.3. Theorem.** *Let $\widetilde{\Omega}_{\mathcal{T}_A}$ be the set of all $\xi \in \widetilde{\Omega}_{\mathcal{TO}_A}$ such that either $\xi$ is unbounded, or it is bounded and $R_\xi(\omega) \in \Sigma_A$, where $\omega$ is the stem of $\xi$. Then $\widetilde{\Omega}_{\mathcal{T}_A}$ is a compact subspace of $\widetilde{\Omega}_{\mathcal{TO}_A}$. Moreover, letting for each $g \in \mathbb{F}$,*

$$\Delta_g^\tau := \Delta_g \cap \widetilde{\Omega}_{\mathcal{T}_A} = \{\xi \in \widetilde{\Omega}_{\mathcal{T}_A} : g \in \xi\}$$

*and*

$$\alpha_g : \xi \in \Delta^\tau_{g^{-1}} \mapsto g\xi \in \Delta^\tau_g$$

*we have that* $\big(\{\Delta^\tau_g\}_{g\in\mathbb{F}}, \{\alpha_g\}_{g\in\mathbb{F}}\big)$ *is a partial action of* $\mathbb{F}$ *on* $\widetilde\Omega_{\mathcal{T}_A}$ *such that the crossed-product* $C(\widetilde\Omega_{\mathcal{T}_A})\rtimes\mathbb{F}$ *is isomorphic to* $\widetilde{\mathcal{T}}_A$ *under an isomorphism which maps each* $1_{\Delta^\tau_x}\delta_x$ *to* $s_x$*.*

*Proof.* Given finite subsets $X$ and $Y$ of $\mathcal{G}$ let $f_{X,Y} : 2^{\mathbb{F}} \to \{0,1\}$ be defined by

$$f_{X,Y}(\xi) = \prod_{x\in X}\big[x^{-1} \in \xi\big] \prod_{y\in Y}\big[y^{-1} \notin \xi\big],$$

where the brackets correspond to the obvious boolean valued function. Following [**EL**: Section 7] and [**ELQ** : 4.4] it suffices to show that $\widetilde\Omega_{\mathcal{T}_A}$ consist of the set of all $\xi \in \widetilde\Omega_{\mathcal{TO}_A}$ such that $f_{X,Y}(g^{-1}\xi) = 0$ for all $g \in \xi$ and all pairs $X, Y$ of finite subsets of $\mathcal{G}$ such that $A(X, Y, \;\cdot\;) \equiv 0$.

We begin by proving the inclusion "$\subseteq$". So let $\xi \in \widetilde\Omega_{\mathcal{T}_A}$ and suppose by contradiction that $f_{X,Y}(g^{-1}\xi) = 1$ for some $g \in \xi$, where $X$ and $Y$ are as above. We have

$$1 = f_{X,Y}(g^{-1}\xi) = \prod_{x\in X}\big[gx^{-1} \in \xi\big] \prod_{y\in Y}\big[gy^{-1} \notin \xi\big],$$

and so for all $x \in X$ and $y \in Y$ one has that $gx^{-1} \in \xi$ and $gy^{-1} \notin \xi$. These translate to $x \in R_\xi(g)$ and $y \notin R_\xi(g)$. Consider the neighborhood of $R_\xi(g)$ within $\Sigma_A$ given by

$$V(X, Y) = \{c \in 2^{\mathcal{G}} : x \in c,\; y \notin c, \quad \forall x \in X,\; \forall y \in Y\}.$$

We claim that $V(X, Y)$ contains at least one column of $A$. The argument here breaks into two cases: if, on the one hand, $g$ is the (finite) stem of $\xi$ then the claim follows from the hypothesis that $\xi \in \widetilde\Omega_{\mathcal{T}_A}$ and hence that $R_\xi(g)$ is in the closure $\Sigma_A$ of the set of columns of $A$. If on the other hand $g$ is not the stem of $\xi$ then there exists $j \in \mathcal{G}$ such that $gj \in \xi$. It follows from the fact that $\xi \in \widetilde\Omega_{\mathcal{TO}_A}$ that $gx^{-1} \in \xi \Leftrightarrow A(x, j) = 1$ and hence that $R_\xi(g)$ coincides with the $j^{th}$ column of $A$, which therefore lies in $V(X, Y)$.

The claim is therefore proven and we may then pick $j$ such that the $j^{th}$ column of $A$ belongs to $V(X, Y)$. It follows that $A(x, j) = 1$ for $x \in X$ and $A(y, j) = 0$ for $y \in Y$ so that $A(X, Y, j) = 1$. This contradicts the fact that $A(X, Y, \;\cdot\;) \equiv 0$.

In order to prove the reverse inclusion let $\xi \in \widetilde\Omega_{\mathcal{TO}_A}$ be such that $f_{X,Y}(g^{-1}\xi) = 0$ whenever $g \in \xi$ and $A(X, Y, \;\cdot\;) \equiv 0$. We want to prove that $\xi \in \widetilde\Omega_{\mathcal{T}_A}$. If $\xi$ is unbounded then this follows by definition. So we assume that $\xi$ is bounded and we must show that $R_\xi(\omega) \in \Sigma_A$, where $\omega$ is the stem of $\xi$. Suppose by contradiction that this is not so and hence there exists a "basic" neighborhood of $R_\xi(\omega)$ of the form $V(X, Y)$ containing no column of $A$.

Since $R_\xi(\omega)$ is obviously in $V(X, Y)$, one has that $\omega x^{-1} \in \xi$ and $\omega y^{-1} \notin \xi$ for all $x$ in $X$ and $y$ in $Y$, which says that

$$f_{X,Y}(\omega^{-1}\xi) = \prod_{x\in X}\big[\omega x^{-1} \in \xi\big] \prod_{y\in Y}\big[\omega y^{-1} \notin \xi\big] = 1. \tag{†}$$

Consider the equation

$$\prod_{x\in X} A(x, j) \prod_{y\in Y}(1 - A(y, j)) = 1$$

in the unknown $j$. The solutions consist, of course, of those $j$'s such that for all $x \in X$ and $y \in Y$ one has that $A(x, j) = 1$ and $A(y, j) = 0$. Thus $j$ is a solution if and only if the $j^{th}$ column of $A$ belongs to $V(X, Y)$. By assumption there is no such column and hence neither are there solutions. In other words $A(X, Y, \;\cdot\;) \equiv 0$. By hypothesis we therefore have that $f_{X,Y}(w^{-1}\xi) = 0$ which contradicts (†). □

There is a subspace of $\widetilde\Omega_{\mathcal{T}_A}$ which will be relevant to us later and this is perhaps the right place to introduce it.

**5.4. Proposition.** *Let* $\widetilde\Omega_e$ *be the subset of* $\widetilde\Omega_{\mathcal{T}_A}$ *consisting of all* $\xi \in \widetilde\Omega_{\mathcal{T}_A}$ *whose stem is equal to* $e$*. Then* $\widetilde\Omega_e$ *is a retract of* $\widetilde\Omega_{\mathcal{T}_A}$ *in the sense that there exists a continuous function* $r : \widetilde\Omega_{\mathcal{T}_A} \to \widetilde\Omega_e$ *such that* $r(\xi) = \xi$ *for all* $\xi \in \widetilde\Omega_e$*. Moreover such a function can be chosen so that* $R_{r(\xi)}(e) = R_\xi(e)$*.*

*Proof.* Given any $\xi$ in $\widetilde{\Omega}_{\mathcal{T}_A}$ let $\eta \in \widetilde{\Omega}_{\mathcal{T}\mathcal{O}_A}$ have trivial stem and be such that $R_\eta(e) = R_\xi(e)$. Such an element exists by [**EL**: 5.14]. We claim that $\eta \in \widetilde{\Omega}_{\mathcal{T}_A}$. To see this we have to consider two cases: on the one hand if $\xi$ has trivial stem then $\eta = \xi$ by the uniqueness part of [**EL**: 5.14] and hence obviously $\eta \in \widetilde{\Omega}_{\mathcal{T}_A}$. If, on the other hand, the stem of $\xi$ is not trivial then there exists some $j \in \mathcal{G} \cap \xi$. It follows from

$$x^{-1} \in \xi \Leftrightarrow A(x, j) = 1$$

that $R_\xi(e)$ coincides with the $j^{th}$ column of $A$ and hence $R_\xi(e) \in \Sigma_A$. Thus $R_\eta(e) \in \Sigma_A$ implying that $\eta \in \widetilde{\Omega}_{\mathcal{T}_A}$.

Define $r(\xi) = \eta$ thus obtaining a function from $\widetilde{\Omega}_{\mathcal{T}_A}$ to $\widetilde{\Omega}_e$ which clearly restricts to the identity on $\widetilde{\Omega}_e$. It now remains to prove that $r$ is continuous. Given that $\widetilde{\Omega}_e$ has a product topology it is enough to verify that the map

$$\xi \in \widetilde{\Omega}_{\mathcal{T}_A} \quad \mapsto \quad \big[g \in r(\xi)\big] \in \{0, 1\}$$

is continuous for all $g \in \mathbb{F}$. Write $g = xg'$ with $x \in \mathcal{G} \cap \mathcal{G}^{-1}$ and $|g| = 1 + |g'|$. Suppose first that $x \in \mathcal{G}$. Since $r(\xi)$ is convex by being in $\widetilde{\Omega}_{\mathcal{T}\mathcal{O}_A}$ and since $x \notin r(\xi)$ because the stem of $r(\xi)$ is trivial we must have $\big[g \in r(\xi)\big] = 0$ for all $\xi$. Suppose now that $x \in \mathcal{G}^{-1}$. Using convexity as well as [**EL**: 5.11] we may prove that $\big[g \in r(\xi)\big] = \big[g \in \xi\big]$ and we again have continuity of our map. ∎

## 6. Partial representations.

We will now consider the map

$$S : \mathbb{F} \to \widetilde{\mathcal{T}}_A$$

defined (*cf.* [**EL**: Section 3]) as follows: if $x$ is in $\mathcal{G}$ put $S(x) = s_x$ and $S(x^{-1}) = s_x^*$. For a general $g \in \mathbb{F}$ write $g = x_1 x_2 \ldots x_n$ in reduced form, that is, each $x_k \in \mathcal{G} \cup \mathcal{G}^{-1}$ and $x_{k+1} \neq x_k^{-1}$, and set

$$S(g) = S(x_1) \cdots S(x_n).$$

The key feature of $S$ (*cf.* [**EL**: 3.2]) is that it is a *partial representation* of $\mathbb{F}$ in the sense that

- $S(e) = I$,
- $S(g^{-1}) = S(g)^*$, and
- $S(g)S(h)S(h^{-1}) = S(gh)S(h^{-1})$,

for all $g, h \in \mathbb{F}$ (see also [**E4**]). We will use, for each $g \in \mathbb{F}$, the notation

- $p_g = S(g)S(g)^*$, and
- $q_g = S(g)^*S(g)$.

One may prove [**E4**: 2.4.iii] that the $p_g$ and $q_g$ form a commutative set. Also each $S(g)$ is a partial isometry with initial and final projections $q_g$ and $p_g$, respectively.

It is easy to see (*cf.* also [**EL**: 3.2]) from the definition of $S$ that it is a *semi-saturated* partial representation in the sense that

- $S(g)S(h) = S(gh)$ whenever $|gh| = |g| + |h|$.

Moreover it is clearly also an *orthogonal* partial representation in the sense that

- $S(x)^*S(y) = 0$ whenever $x, y \in \mathcal{G}$ are such that $x \neq y$.

Partial isometries in general tend to behave very badly, rarely satisfying any algebraic properties at all. For instance the square of a partial isometry may not be a partial isometry. However, as we shall see, the fact that our partial isometries $s_x$ are assembled into a partial representation will make it much easier for us to deal with them. Actually this was the main technical tool in bringing the Cuntz–Krieger algebras for infinite matrices to life [**EL**].

Let us now relate $S$ to the crossed-product structure of $\widetilde{\mathcal{T}}_A$. Recall from 5.3 that $\big(\{\Delta_g^\tau\}_{g\in\mathbb{F}}, \{\alpha_g\}_{g\in\mathbb{F}}\big)$ is a partial action of $\mathbb{F}$ on $\widetilde{\Omega}_{\mathcal{T}_A}$ whose crossed-product is isomorphic to $\widetilde{\mathcal{T}}_A$ in such a way that each $s_x$ corresponds to $1_{\Delta_x^\tau}\delta_x$. At the algebra level let us agree to denote by

$$\theta_g : C_0(\Delta_{g^{-1}}^\tau) \to C_0(\Delta_g^\tau)$$

the *-isomorphism given by $\theta_g(f) = f \circ \alpha_{g^{-1}}$ for all $f \in C_0(\Delta_{g^{-1}}^\tau)$.

**6.1. Proposition.** *Let $g \in \mathbb{F}$. Then*

(i) $S(g) = 1_{\Delta_g^\tau}\delta_g$,

(ii) $p_g = 1_{\Delta_g^\tau}$,

(iii) $q_g = 1_{\Delta_{g^{-1}}^\tau}$,

(iv) $S(g)aS(g)^* = \theta_g(a)$ *for all* $a \in C_0(\Delta_{g^{-1}}^\tau)$,

(v) $S(g)aS(g)^* = \theta_g(q_g a)$ *for all* $a \in C(\widetilde{\Omega}_{\mathcal{T}_A})$, *and*

(vi) $\theta_g(q_g) = p_g$.

*Proof.* In order to prove (i) we will use induction on $|g|$. If $|g| = 0$ the result is obvious. If $g = x \in \mathcal{G}$ we have $S(x) = s_x = 1_{\Delta_x^\tau}\delta_x$, as already mentioned. It is an easy exercise to show that this implies the result also for $g \in \mathcal{G}^{-1}$. So suppose that $|g| > 1$ and write $g = rs$ with $|g| = |r| + |s|$ and $|r|, |s| < |g|$. Using that $S$ is semi-saturated and the induction hypothesis we have

$$S(g) = S(r)S(s) = (1_{\Delta_r^\tau}\delta_r)(1_{\Delta_s^\tau}\delta_s) = \theta_r\big(\theta_r^{-1}(1_{\Delta_r^\tau})1_{\Delta_s^\tau}\big)\delta_{rs} =$$

$$= \theta_r\left(1_{\Delta_{r^{-1}}^\tau}1_{\Delta_s^\tau}\right)\delta_{rs} = \theta_r\left(1_{\Delta_{r^{-1}}^\tau\cap\Delta_s^\tau}\right)\delta_{rs} = 1_{\alpha_r\left(\Delta_{r^{-1}}^\tau\cap\Delta_s^\tau\right)}\delta_{rs}.$$

On the other hand observe that

$$\alpha_r\left(\Delta_{r^{-1}}^\tau \cap \Delta_s^\tau\right) = \Delta_r^\tau \cap \Delta_{rs}^\tau = \Delta_{rs}^\tau,$$

where the first equality follows from the fact that $\alpha$ is a partial action and the second by semi-saturatedness (see 4.1.a). Therefore

$$S(g) = 1_{\Delta_{rs}^\tau}\delta_{rs} = 1_{\Delta_g^\tau}\delta_g,$$

concluding the proof of (i). In order to prove (ii) we have

$$p_g = S(g)S(g)^* = (1_{\Delta_g^\tau}\delta_g)(1_{\Delta_{g^{-1}}^\tau}\delta_{g^{-1}}) = \theta_g\big(\theta_{g^{-1}}(1_{\Delta_g^\tau})1_{\Delta_{g^{-1}}^\tau}\big)\delta_e = \theta_g\big(1_{\Delta_{g^{-1}}^\tau}\big)\delta_e = 1_{\Delta_g^\tau}\delta_e = 1_{\Delta_g^\tau},$$

where we identify, as usual, $a$ and $a\delta_e$ for all $a$ in $C(\widetilde{\Omega}_{\mathcal{T}_A})$. Point (iii) follows from (ii) simply by replacing $g$ by $g^{-1}$. Regarding (iv) let $a \in C_0(\Delta_{g^{-1}}^\tau)$. We have

$$S(g)a = (1_{\Delta_g^\tau}\delta_g)(a\delta_e) = \theta_g\big(\theta_{g^{-1}}(1_{\Delta_g^\tau})a\big)\delta_g = \theta_g\big(1_{\Delta_{g^{-1}}^\tau}a\big)\delta_g = \theta_g(a)\delta_g.$$

Therefore

$$S(g)aS(g)^* = (\theta_g(a)\delta_g)(1_{\Delta_{g^{-1}}^\tau}\delta_{g^{-1}}) = \theta_g(a1_{\Delta_{g^{-1}}^\tau})\delta_e = \theta_g(a)\delta_e = \theta_g(a).$$

As for (v) we have for all $a \in C(\widetilde{\Omega}_{\mathcal{T}_A})$ that

$$S(g)aS(g)^* = S(g)S(g)^*S(g)aS(g)^* = S(g)q_g aS(g)^* = \theta_g(q_g a),$$

by (iv) because $q_g a = 1_{\Delta_{g^{-1}}^\tau}a \in C_0(\Delta_{g^{-1}}^\tau)$. Finally (vi) follows by plugging $a = 1$ in (v). ∎

We now wish to name a subalgebra of $\widetilde{\mathcal{T}}_A$ which will play an important role alongside the partial representation $S$ above.

**6.2. Definition.** We will let $\widetilde{\mathcal{Q}}$ be the subalgebra of $\widetilde{\mathcal{T}}_A$ generated by the set $\{q_x : x \in \mathcal{G}\} \cup \{1\}$.

Note that $\widetilde{\mathcal{T}}_A$ is in fact a subalgebra of $C(\widetilde{\Omega}_{\mathcal{T}_A})$. We will now discuss certain properties of $S$ in relation to $\widetilde{\mathcal{Q}}$.

**6.3. Proposition.** *For $\mu, \nu$ in $\mathbb{F}_+$ one has*

(i) *If $|\mu| = |\nu|$ but $\mu \neq \nu$ then $S(\mu)^* S(\nu) = 0$.*

(ii) *If $|\mu| \geq 1$ and $z$ is the last generator in the reduced decomposition of $\mu$ then $q_\mu = \varepsilon q_z$, where $\varepsilon$ is either 1 or 0 according to whether $\mu$ is admissible (i.e. $A(\mu_i, \mu_{i+1}) = 1$ for all $i = 1, \ldots, |\mu| - 1$) or not.*

(iii) *If $|\mu|$ and $z$ are as in (ii) then $S(\mu)^* \widetilde{\mathcal{Q}} S(\mu) \subseteq \mathbb{C} q_z \subseteq \widetilde{\mathcal{Q}}$.*

(iv) *If $|\mu| \leq |\nu|$ then $\big(S(\mu)\widetilde{\mathcal{Q}}S(\mu)^*\big)\big(S(\nu)\widetilde{\mathcal{Q}}S(\nu)^*\big) \subseteq S(\nu)\widetilde{\mathcal{Q}}S(\nu)^*$.*

*Proof.* Statements (i–ii) follow from claims 2 and 1, respectively, in the proof of [**EL**: Proposition 3.2]. To prove (iii) let $x \in \mathcal{G}$ and observe that

$$S(\mu)^* q_x S(\mu) = S(x\mu)^* S(x\mu) = q_{x\mu} = \varepsilon q_z,$$

by (ii). Now let $x_1, \ldots, x_n \in \mathcal{G}$ and observe that, since $S(\mu)$ is a partial isometry, we have that $S(\mu)^* = S(\mu)^* S(\mu) S(\mu)^*$ and hence

$$S(\mu)^* q_{x_1} \ldots q_{x_n} S(\mu) = S(\mu)^* q_{x_1} S(\mu) S(\mu)^* q_{x_2} \ldots q_{x_{n-1}} S(\mu) S(\mu)^* q_{x_n} S(\mu).$$

Since each $S(\mu)^* q_{x_i} S(\mu)$ belongs to $\mathbb{C} q_z$, the same holds for $S(\mu)^* q_{x_1} \ldots q_{x_n} S(\mu)$. Since $\widetilde{\mathcal{Q}}$ is linearly spanned by the set of products $q_{x_1} \ldots q_{x_n}$, we have proven (iii).

With respect to (iv) write $\nu = \nu'\nu''$, with $|\nu'| = |\mu|$ and notice that $S(\mu)^* S(\nu) = S(\mu)^* S(\nu') S(\nu'')$ vanishes by (i) if $\mu \neq \nu'$. So, assuming that $\mu = \nu'$, we have

$$\big(S(\mu)\widetilde{\mathcal{Q}}S(\mu)^*\big)\big(S(\nu)\widetilde{\mathcal{Q}}S(\nu)^*\big) = S(\mu)\widetilde{\mathcal{Q}}S(\mu)^* S(\mu) S(\nu'')\widetilde{\mathcal{Q}}S(\nu)^* = S(\mu)\widetilde{\mathcal{Q}}S(\nu'')\widetilde{\mathcal{Q}}S(\nu)^* =$$

$$= S(\mu)S(\nu'')S(\nu'')^*\widetilde{\mathcal{Q}}S(\nu'')\widetilde{\mathcal{Q}}S(\nu)^* = S(\nu)S(\nu'')^*\widetilde{\mathcal{Q}}S(\nu'')\widetilde{\mathcal{Q}}S(\nu)^* \subseteq S(\nu)\widetilde{\mathcal{Q}}S(\nu)^*,$$

where we have used (iii) in the last step. ☐

The next result gives a total set for $C(\widetilde{\Omega}_{\mathcal{T}_A})$ in terms of $\widetilde{\mathcal{Q}}$ and $S(\mathbb{F}_+)$.

**6.4. Proposition.** *$C(\widetilde{\Omega}_{\mathcal{T}_A})$ coincides with the closed linear span of the set*

$$\{S(\mu) a S(\mu)^* \colon \mu \in \mathbb{F}_+,\ a \in \widetilde{\mathcal{Q}}\}.$$

*Proof.* By 6.1.ii we have that $p_g = 1_{\Delta_g^{\mathcal{T}}}$. The set of all $p_g$'s therefore separates points of $\widetilde{\Omega}_{\mathcal{T}_A}$ and hence generates $C(\widetilde{\Omega}_{\mathcal{T}_A})$ as a $C^*$-algebra. We claim that every nonzero $p_g$ belongs to the set in the statement. To see this note that $S(g) = 0 = p_g$ unless $g = \mu\nu^{-1}$, where $\mu, \nu \in \mathbb{F}_+$ are admissible words such that $|g| = |\mu| + |\nu|$, because $S$ is semi-saturated and orthogonal (see [**EL**: 3.1]). Let us therefore suppose that $g$ is of this form. If $|\nu| = 0$ then the claim is obvious. Otherwise let $z$ be the last generator in the reduced decomposition of $\nu$ and observe that

$$p_g = S(g)S(g)^* = S(\mu)S(\nu)^*S(\nu)S(\mu)^* = S(\mu)q_\nu S(\mu)^* = S(\mu)q_z(\mu)^*,$$

by 6.3.ii. It is now enough to show that the set in the statement is closed under multiplication, but this follows immediately from 6.3.iv. ☐

## 7. Unital and non-unital algebras.

In this section we propose to extend part of the discussion about units found in [**EL**: Section 8] to $\widetilde{\mathcal{TO}}_A$ and $\widetilde{\mathcal{T}}_A$. Recall that $\widetilde{\mathcal{TO}}_A$, $\widetilde{\mathcal{T}}_A$, and $\widetilde{\mathcal{O}}_A$ were defined via universal properties in the category of *unital* $C^*$-algebras and hence they are obviously unital. However all of them have possibly non-unital counterparts which are also of interest. If $\widetilde{B}$ denotes any one of these algebras and $\{s_x\}_{x\in\mathcal{G}}$ is the canonical set of generating partial isometries we shall also consider the (non-necessarily unital) sub-$C^*$-algebra $B$ of $\widetilde{B}$ generated by the set $\{s_x : x \in \mathcal{G}\}$. $B$ will be denoted, respectively, by $\mathcal{TO}_A$, $\mathcal{T}_A$, and $\mathcal{O}_A$, filling the third column of the following table:

| Relations | algebras $\widetilde{B}$ | algebras $B$ | spaces $\widetilde{\Omega}$ | spaces $\Omega$ |
|---|---|---|---|---|
| $\text{CK}_{1-3}$ | $\widetilde{\mathcal{TO}}_A$ | $\mathcal{TO}_A$ | $\widetilde{\Omega}_{\mathcal{TO}_A}$ | $\Omega_{\mathcal{TO}_A}$ |
| $\text{CK}_{1-3} + \text{CK}_4^0$ | $\widetilde{\mathcal{T}}_A$ | $\mathcal{T}_A$ | $\widetilde{\Omega}_{\mathcal{T}_A}$ | $\Omega_{\mathcal{T}_A}$ |
| $\text{CK}_{1-4}$ | $\widetilde{\mathcal{O}}_A$ | $\mathcal{O}_A$ | $\widetilde{\Omega}_{\mathcal{O}_A}$ | $\Omega_{\mathcal{O}_A}$ |

*Table 7.1*

It may or may not happen that $B = \widetilde{B}$ but it is nevertheless true that $B \cup \{1\}$ generates $\widetilde{B}$. Therefore $\widetilde{B}$ can be seen as the unitization of $B$. It also follows that the codimension of $B$ in $\widetilde{B}$ is at most one.

Let $\widetilde{\Omega}$ be either one of $\widetilde{\Omega}_{\mathcal{TO}_A}$, $\widetilde{\Omega}_{\mathcal{T}_A}$, or $\widetilde{\Omega}_{\mathcal{O}_A}$ according to the fourth column of table 7.1 so that $\widetilde{B} \simeq C(\widetilde{\Omega}) \rtimes \mathbb{F}$ as seen above. Let $\epsilon = \{e\}$ be seen as a subset of $\mathbb{F}$ and hence as an element of $2^{\mathbb{F}}$, which one may easily show lies in $\widetilde{\Omega}_{\mathcal{TO}_A}$. It may or may not happen that $\epsilon \in \widetilde{\Omega}$ but we shall nevertheless let $\Omega = \widetilde{\Omega} \setminus \{\epsilon\}$, leading to the space indicated in the last column of table 7.1. $\Omega$ is then a locally compact topological space which is clearly invariant under the partial action $\alpha$ of $\mathbb{F}$ on $\widetilde{\Omega}_{\mathcal{TO}_A}$.

**7.2. Proposition.** *Choose a row in table 7.1 and let $\widetilde{B}$, $B$, $\widetilde{\Omega}$ and $\Omega$ be accordingly chosen. Then $\widetilde{B} = C(\widetilde{\Omega}) \rtimes \mathbb{F}$ $B = C_0(\Omega) \rtimes \mathbb{F}$.*

*Proof.* The three cases corresponding to the column labeled $\widetilde{B}$ have already been dealt with, whereas the case of $\mathcal{O}_A$ is treated in [**EL**: 8.4.iii]. We therefore focus on the remaining cases, i.e. $\mathcal{TO}_A$ and $\mathcal{T}_A$.

Clearly $C_0(\Omega)$ is an invariant ideal in $C(\widetilde{\Omega})$ so that we may use [**ELQ**: 3.1] to conclude that $C_0(\Omega) \rtimes \mathbb{F}$ is an ideal in $C(\widetilde{\Omega}) \rtimes \mathbb{F}$. Given $x \in \mathcal{G}$ let $1_x := 1_{\Delta_x \cap \widetilde{\Omega}} \in C(\widetilde{\Omega})$ and observe that $1_x(\epsilon) = [x \in \epsilon] = 0$. Therefore $1_x \in C_0(\Omega)$ and hence $s_x = 1_x \delta_x \in C_0(\Omega) \rtimes \mathbb{F}$. So

$$B \subseteq C_0(\Omega) \rtimes \mathbb{F} \subseteq C(\widetilde{\Omega}) \rtimes \mathbb{F} = \widetilde{B}.$$

Suppose by contradiction that $B$ is a proper subset of $C_0(\Omega) \rtimes \mathbb{F}$. As observed above the codimension of $B$ in $\widetilde{B}$ is at most one, hence $C_0(\Omega) \rtimes \mathbb{F} = C(\widetilde{\Omega}) \rtimes \mathbb{F}$. Applying [**ELQ**: 3.1] we conclude that $\Omega = \widetilde{\Omega}$ and hence that $\epsilon \notin \widetilde{\Omega}$. In the case of $\mathcal{TO}_A$ this is already a contradiction because $\epsilon$ does belong to $\widetilde{\Omega}_{\mathcal{TO}_A}$ as already mentioned. So it remains to consider $\mathcal{T}_A$. The characterization of $\widetilde{\Omega}_{\mathcal{T}_A}$ given in 5.3 says that $R_\epsilon(e)$, namely the zero vector, is not in the closure of the set of columns of $A$. This implies that there exists a finite set $Y \subseteq \mathcal{G}$ such that the basic neighborhood of the zero vector in $2^{\mathcal{G}}$ given by

$$V(\emptyset, Y) = \{c \in 2^{\mathcal{G}} : \ y \notin c, \ \forall y \in Y\}$$

contains no column of $A$. It follows that $A(\emptyset, Y, \ \cdot\ ) \equiv 0$ and hence we have by $\text{CK}_4^0$ that

$$0 = q(\emptyset, Y) = \textstyle\prod_{y \in Y}(1 - q_y).$$

Upon expanding the right hand side above we discover that 1 is in the algebra generated by the $q_y$'s and hence also that $1 \in B$. This implies that $B = \widetilde{B}$ leading to a contradiction. ∎

Table 7.1 therefore displays six algebras (for each matrix $A$) which admit a crossed-product structure and hence we may use the results of Section 2 to study their KMS states, or rather at least those which factor through the conditional expectation. However we wish to be able to apply the much stronger Theorem 4.3 which requires the corresponding partial actions to be semi-saturated and orthogonal.

**7.3. Proposition.** *The partial action of $\mathbb{F}$ on each one of the six spaces appearing in the last two columns of table 7.1 is semi-saturated and orthogonal.*

*Proof.* We start by verifying conditions 4.1.a-b for the case of $\widetilde{\Omega}_{\mathcal{TO}_A}$. Beginning with 4.1.a let $g,h\in\mathbb{F}$ be such that $|gh|=|g|+|h|$. In this case the shortest path from $e$ to $gh$ in the Cayley graph of $\mathbb{F}$ must pass through $g$. If $\xi\in\Delta_{gh}$ then both $e$ and $gh$ lie in $\xi$ and hence so does $g$ by convexity [**EL**: 4.4]. Therefore $g\in\xi$ and $\xi\in\Delta_g$. This proves that $\Delta_{gh}\subseteq\Delta_g$.

Let us now check 4.1.b. Suppose that $x,y\in\mathcal{G}$, $x\neq y$, and $\xi\in\Delta_x\cap\Delta_y$. Then $\{e,x,y\}\subseteq\xi$ which contradicts the penultimate property defining $\widetilde{\Omega}_{\mathcal{TO}_A}$ in 5.1. Therefore $\Delta_x\cap\Delta_y=\emptyset$.

Considering the other five partial actions under analysis observe that they are all obtained by restricting the one for $\widetilde{\mathcal{TO}}_A$ to invariant subsets. Properties 4.1.a-b immediately follow and so the proof is concluded. □

## 8. Scaling states and the partition function $Z(\beta)$.

We are already working under the choice of a fixed 0–1 matrix $A$ and now we are about to make other important standing hypotheses. For ease of reference we record them here.

**8.1. Standing Hypothesis.** From now on and throughout the rest of this work we will let

(i) $\mathcal{G}$ be a countable set (meaning finite or countably infinite),

(ii) $A=\{A(x,y)\}_{x,y\in\mathcal{G}}$ be a 0–1 matrix having no identically zero rows,

(iii) $\{N(x)\}_{x\in\mathcal{G}}$ be a collection of real numbers in the interval $(1,\infty)$,

(iv) $\sigma$ be the unique one-parameter group of automorphisms of each one of the algebras in table 7.1 (by abuse of notation) satisfying $\sigma_t(s_x)=N(x)^{it}s_x$, and

(v) all references to KMS states will be with respect to the one-parameter group $\sigma$ above.

The existence of $\sigma$, in any one of its versions, may be deduced either from 4.3 or from the universal properties of our algebras.

We begin with an important consequence of 4.3 and 7.3:

**8.2. Corollary.** *Under 8.1 let $B$ be any one of the $C^*$-algebras:*

$$\widetilde{\mathcal{TO}}_A,\ \mathcal{TO}_A,\ \widetilde{\mathcal{T}}_A,\ \mathcal{T}_A,\ \widetilde{\mathcal{O}}_A,\ \text{and}\ \mathcal{O}_A$$

*and let $\Omega$ be the respective space chosen from*

$$\widetilde{\Omega}_{\mathcal{TO}_A},\ \Omega_{\mathcal{TO}_A},\ \widetilde{\Omega}_{\mathcal{T}_A},\ \Omega_{\mathcal{T}_A},\ \widetilde{\Omega}_{\mathcal{O}_A},\ \text{and}\ \Omega_{\mathcal{O}_A},$$

*so that $B\simeq C_0(\Omega)\rtimes\mathbb{F}$ as seen above. Given $\beta\in(0,\infty]$ the correspondence $\phi\mapsto\phi\circ E$ is an affine homeomorphism between the set of states $\phi$ on $C_0(\Omega)$ satisfying*

(i) *if $\beta<\infty$ : $\phi(\theta_x(a))=N(x)^{-\beta}\phi(a)$ for all $x\in\mathcal{G}$ and for all $a\in C_0(\Delta_{x^{-1}}\cap\Omega)$,*

(ii) *if $\beta=\infty$ : $\phi(C_0(\Delta_x\cap\Omega))=\{0\}$ for all $x\in\mathcal{G}$,*

*and the KMS$_\beta$ states on $B$. If $\lambda$ is the probability measure on $\Omega$ corresponding via the Riesz Representation Theorem to $\phi$ then (i–ii) are respectively equivalent to:*

*(i′) if $\beta<\infty$ : $\lambda(\alpha_x(S))=N(x)^{-\beta}\lambda(S)$ for all Borel subsets $S\subseteq\Delta_{x^{-1}}\cap\Omega$.*

*(ii′) if $\beta=\infty$ : $\lambda(\Delta_x\cap\Omega)=\{0\}$ for all $x\in\mathcal{G}$.*

The states and measures of 8.2 will evidently become the main players in this theory and hence they deserve a name:

**8.3. Definition.** Let $\Omega$ be as in 8.2, let $\phi$ be a state on $C_0(\Omega)$ and let $\lambda$ be the probability measure on $\Omega$ associated to $\phi$ via the Riesz Representation Theorem. Then $\phi$ will be called a *$\beta$-scaling state*, and $\lambda$ will be called a *$\beta$-scaling measure*, where $\beta \in (0, \infty]$, if the conditions of 8.2 are satisfied.

It is now perhaps the right time for us to make a choice among the six algebras of Theorem 8.2. From now on we shall concentrate our study on $\mathcal{T}_A$ for several reasons, namely:

- because the study of KMS states is most interesting in this case,
- because $\mathcal{T}_A$ is the Toeplitz extension of $\mathcal{O}_A$, viewing the latter as a Cuntz–Pimsner algebra [**Sz**], and
- because every KMS state on $\mathcal{O}_A$ gives a KMS state on $\mathcal{T}_A$, by composing with the canonical quotient map, and hence we include $\mathcal{O}_A$ in the process.

The KMS states of $\mathcal{TO}_A$ are likely to be interesting as well, since they include everything else, again by considering the quotient maps. Moreover they could be studied with much the same tools we shall use here. But, alas, we wont be looking at them in this work.

Regarding $\widetilde{\mathcal{T}}_A$ recall from the beginning of Section 7 that $\Omega_{\mathcal{T}_A} = \widetilde{\Omega}_{\mathcal{T}_A} \setminus \{\epsilon\}$. By definition of $\widetilde{\Omega}_{\mathcal{T}_A}$ (see 5.3) one has that $\epsilon \in \widetilde{\Omega}_{\mathcal{T}_A}$ if and only $R_\epsilon(e)$ (which is clearly the empty set, or the zero vector in $2^{\mathbb{F}}$) belongs to the closure of the set of columns of $A$, namely $\Sigma_A$. So when the zero vector is not in $\Sigma_A$ one has that $\Omega_{\mathcal{T}_A} = \widetilde{\Omega}_{\mathcal{T}_A}$ and hence also $\mathcal{T}_A = \widetilde{\mathcal{T}}_A$.

Nevertheless it is quite possible that the zero vector lie in $\Sigma_A$, and hence to study $\beta$-scaling measures on $\Omega_{\mathcal{T}_A}$ is, strictly speaking, not the same as to do so for $\widetilde{\Omega}_{\mathcal{T}_A}$. The difference however is not very deep in the sense that a probability measure on $\widetilde{\Omega}_{\mathcal{T}_A}$ is given, in an essentially unique way, by a convex combination of a probability measure on $\Omega_{\mathcal{T}_A}$ and the Dirac measure $\delta_\epsilon$, i.e. the measure on $\widetilde{\Omega}_{\mathcal{T}_A}$ assigning mass one to the point $\epsilon$.

**8.4. Proposition.** *Suppose that $\epsilon \in \widetilde{\Omega}_{\mathcal{T}_A}$. Then:*

(i) *$\delta_\epsilon$ is $\beta$-scaling for all $\beta$ in $(0, \infty]$,*

(ii) *for $\beta$ in $(0, \infty]$ the $\beta$-scaling measures on $\widetilde{\Omega}_{\mathcal{T}_A}$ consist precisely of $\delta_\epsilon$ and the convex combinations of a $\beta$-scaling measure on $\Omega_{\mathcal{T}_A}$ and $\delta_\epsilon$.*

*Proof.* Part (i) follows easily from the fact that $\delta_\epsilon(\Delta_g^\tau) = 0$ for all $g \in \mathbb{F} \setminus \{e\}$. Part (ii) is then evident. □

Therefore, once we classify all $\beta$-scaling measures on $\Omega_{\mathcal{T}_A}$, we will be able to transfer that knowledge to $\widetilde{\Omega}_{\mathcal{T}_A}$. By 8.2 we will therefore have classified the KMS states on $\widetilde{\mathcal{T}}_A$. This said we shall now restrict our attention to studying the case of $\mathcal{T}_A$.

A property of $\beta$-scaling states which will be often used is described in our next:

**8.5. Proposition.** *Let $\beta \in (0, \infty)$ and let $\phi$ be a $\beta$-scaling state on $C_0(\Omega_{\mathcal{T}_A})$. Then for every $\mu \in \mathbb{F}_+$ one has that*

$$\phi(p_\mu) = N(\mu)^{-\beta} \phi(q_\mu).$$

*Proof.* By 6.1.vi we have $\phi(p_\mu) = \phi(\theta_\mu(q_\mu)) = N(\mu)^{-\beta}\phi(q_\mu)$. □

We shall now collect some notations to be used sooner or later in this and the following sections.

**8.6. Definition.** We will denote by:

(i) $\Omega_\mu$ the subset of $\Omega_{\mathcal{T}_A}$ formed by the $\xi \in \Omega_{\mathcal{T}_A}$ whose stem coincides with $\mu$ for $\mu \in \mathbb{F}_+$,

(ii) $\Omega_e^\mu$ the intersection $\Omega_e \cap \Delta_{\mu^{-1}}^\tau$ for $\mu \in \mathbb{F}_+$,

(iii) $\Omega_f$ the set of bounded elements of $\Omega_{\mathcal{T}_A}$, i.e. elements with finite stem,

(iv) $\Omega_\infty$ the set formed by the unbounded elements of $\Omega_{\mathcal{T}_A}$,

(v) $P_A$ the subset of $\mathbb{F}_+$ formed by all admissible words, and

(vi) $P_A^n$ the subset of $P_A$ formed by the admissible words of length $n$.

Regarding 8.6.ii note that, again by definition of $\widetilde{\Omega}_{\mathcal{TO}_A}$, when $\mu$ is a non-trivial admissible word one has that $\Delta^\tau_{\mu^{-1}} = \Delta^\tau_{x^{-1}}$, where $x$ is the last generator in the reduced decomposition of $\mu$. So $\Omega_e^\mu = \Omega_e \cap \Delta^\tau_{x^{-1}}$.

Recall from 5.4 that $\widetilde{\Omega}_e$ is the set of all $\xi \in \widetilde{\Omega}_{\mathcal{T}_A}$ whose stem is equal to $e$. Therefore, since $\Omega_{\mathcal{T}_A} = \widetilde{\Omega}_{\mathcal{T}_A} \setminus \{\epsilon\}$, we have that $\Omega_e = \widetilde{\Omega}_e \setminus \{\epsilon\}$. Also observe that by definition of $\widetilde{\Omega}_{\mathcal{TO}_A}$, and since $\Omega_\mu \subseteq \Omega_{\mathcal{T}_A} \subseteq \widetilde{\Omega}_{\mathcal{TO}_A}$, if $\mu$ is not admissible then $\Omega_\mu = \emptyset$ (see [**EL**: 5.4]).

The following are a few easy consequences of the definition:

**8.7. Proposition.** *Indicating by $\dot\cup$ the disjoint union of sets we have:*

(i) $\Omega_{\mathcal{T}_A} = \Omega_\infty \,\dot\cup\, \Omega_f$

(ii) $\Omega_f = \dot\bigcup_{\mu\in\mathbb{F}_+} \Omega_\mu = \dot\bigcup_{\mu\in P_A} \Omega_\mu$

(iii) $\Omega_{\mathcal{T}_A} = \Omega_e \,\dot\cup \left(\dot\bigcup_{x\in\mathcal{G}} \Delta^\tau_x\right)$

(iv) *If $\mu \in \mathbb{F}_+$ is admissible then $\Omega_\mu = \alpha_\mu(\Omega_e^\mu)$.*

**8.8. Definition.** Let $\phi$ be a state on $C_0(\Omega_{\mathcal{T}_A})$ and let $\lambda$ be the probability measure on $\Omega_{\mathcal{T}_A}$ associated to $\phi$ via the Riesz Representation Theorem. Then both $\phi$ and $\lambda$ will be said to be of

(i) *finite type* if $\lambda(\Omega_f) = 1$,

(ii) *infinite type* if $\lambda(\Omega_\infty) = 1$.

Observe that both $\Omega_f$ and $\Omega_\infty$ are invariant under $\alpha$. Therefore, given any $\beta$-scaling measure $\lambda$, where $\beta \in (0,\infty]$, the restriction of $\lambda$ to either one of $\Omega_f$ and $\Omega_\infty$ satisfies (i′) or (ii′). This yields:

**8.9. Proposition.** *Every $\beta$-scaling measure $\lambda$ on $\Omega_{\mathcal{T}_A}$ which is not of finite nor of infinite type can be written in a unique way as a convex combination of a finite type $\beta$-scaling measure $\lambda_f$ and an infinite type $\beta$-scaling measure $\lambda_\infty$.*

It is easy to characterize the infinite type $\beta$-scaling measures:

**8.10. Proposition.** *Let $\beta \in (0,\infty)$ and let $\lambda$ be a $\beta$-scaling measure on $\Omega_{\mathcal{T}_A}$. Then the following are equivalent:*

(i) $\lambda(\Omega_e) = 0$.

(ii) *$\lambda$ is of infinite type.*

*Proof.* (i)⇒(ii): By 8.7.iv one has that

$$\lambda(\Omega_\mu) = \lambda(\alpha_\mu(\Omega_e^\mu)) = N(\mu)^{-\beta}\lambda(\Omega_e \cap \Delta^\tau_{x^{-1}}) = 0$$

for all admissible words $\mu \in \mathbb{F}_+$ ending in $x$. So $\lambda(\Omega_f) = 0$ by 8.7.ii and the assumption that $\mathcal{G}$ is countable.

That (ii)⇒(i) follows from: $\lambda(\Omega_e) \le \lambda(\Omega_f) = 1 - \lambda(\Omega_\infty) = 0$. □

The appropriate form of the above result for the case $\beta = \infty$ is given by:

**8.11. Proposition.** *A probability measure $\lambda$ on $\Omega_{\mathcal{T}_A}$ is an $\infty$-scaling measure if and only if $\lambda(\Omega_e) = 1$. In particular every $\infty$-scaling measure is of finite type.*

*Proof.* Follows immediately from 8.7.iii. □

Given a $\beta$-scaling measure, regardless of it being of finite or infinite type, it is possible to compute the measure of $\Omega_\infty$ as follows:

**8.12. Lemma.** *Let $\beta \in (0,\infty)$ and let $\phi$ be a $\beta$-scaling state on $C_0(\Omega_{\mathcal{T}_A})$ with associated measure $\lambda$. Then*

$$\lambda(\Omega_\infty) = \lim_{n\to\infty} \sum_{\mu\in P_A^n} N(\mu)^{-\beta}\phi(q_\mu).$$

*Proof.* For each integer $n$ consider the subset $S_n$ of $\Omega_{\mathcal{T}_A}$ formed by all those $\xi \in \Omega_{\mathcal{T}_A}$ whose stem has length bigger than or equal to $n$. Observing that

$$S_n = \dot{\bigcup}_{\mu \in P_A^n} \Delta_\mu^\tau,$$

that $\Delta_\mu^\tau = \emptyset$ unless $\mu$ is in $P_A$, and using 8.5, we have

$$\lambda(S_n) = \sum_{\mu \in P_A^n} \lambda(\Delta_\mu^\tau) = \sum_{\mu \in P_A^n} \phi(p_\mu) = \sum_{\mu \in P_A^n} N(\mu)^{-\beta}\phi(q_\mu).$$

Clearly $\Omega_\infty = \bigcap_{n\in\mathbb{N}} S_n$ and the $S_n$ are decreasing. Therefore

$$\lambda(\Omega_\infty) = \lim_{n\to\infty} \lambda(S_n) = \lim_{n\to\infty} \sum_{\mu \in P_A^n} N(\mu)^{-\beta}\phi(q_\mu). \qquad \square$$

Perhaps the most important consequence to be drawn from 8.12 is the following:

**8.13. Proposition.** *Let $\beta \in (0,\infty)$ and suppose that*

$$\lim_{n\to\infty} \sum_{\mu \in P_A^n} N(\mu)^{-\beta} = 0.$$

*Then every $\beta$-scaling state on $C_0(\Omega_{\mathcal{T}_A})$ is of finite type.*

*Proof.* Let $\phi$ be a $\beta$-scaling state with associated measure $\lambda$. From 8.12 we have

$$\lambda(\Omega_\infty) = \lim_{n\to\infty} \sum_{\mu \in P_A^n} N(\mu)^{-\beta}\phi(q_\mu) \le \lim_{n\to\infty} \sum_{\mu \in P_A^n} N(\mu)^{-\beta} = 0,$$

and hence $\lambda(\Omega_f) = 1 - \lambda(\Omega_\infty) = 1$. $\square$

**8.14. Definition.** The *partition function* for the dynamical system $(\mathcal{T}_A, \sigma, \mathbb{R})$ is the function $Z(\beta)$ given by the Dirichlet series

$$Z(\beta) = \sum_{\mu \in P_A} N(\mu)^{-\beta}.$$

Since $P_A = \dot{\bigcup}_{n\in\mathbb{N}} P_A^n$, it is clear that $Z(\beta) = \sum_{n\in\mathbb{N}} \sum_{\mu\in P_A^n} N(\mu)^{-\beta}$. Therefore the convergence of the series for $Z(\beta)$ implies the hypothesis of 8.13. We therefore get the following special case of 8.13:

**8.15. Corollary.** *Suppose that $\beta \in (0,\infty)$ is such that the series for $Z(\beta)$ converges. Then every $\beta$-scaling state on $C_0(\Omega_{\mathcal{T}_A})$ is of finite type.*

Observe that, by 8.11, the $\infty$-scaling states are also all of finite type. This may be seen as a generalization of the above result if one adopts the convention that $N(x)^{-\infty} = 0$.

The convergence of the series for $Z(\beta)$ is not an extremely rare phenomena. For example:

**8.16. Proposition.** *If $\beta \in (0,\infty)$ and $\sum_{x\in\mathcal{G}} N(x)^{-\beta} < 1$ then*

$$Z(\beta) \le \frac{1}{1 - \sum_{x\in\mathcal{G}} N(x)^{-\beta}}.$$

*Proof.* Observe that $P_A \subseteq \mathbb{F}_+ = \dot{\bigcup_{n\in\mathbb{N}}} \mathbb{F}_+^n$, where $\mathbb{F}_+^n$ denotes the subset of $\mathbb{F}_+$ consisting of words of length $n$. Therefore

$$Z(\beta) \ \leq \ \sum_{\mu\in\mathbb{F}_+} N(\mu)^{-\beta} = \sum_{n=0}^{\infty}\sum_{\mu\in\mathbb{F}_+^n} N(\mu)^{-\beta} \ = \ \sum_{n=0}^{\infty}\left(\sum_{x\in\mathcal{G}} N(x)^{-\beta}\right)^n \ = \frac{1}{1-\sum_{x\in\mathcal{G}} N(x)^{-\beta}}. \qquad \square$$

For every Dirichlet series there exists a critical value $\bar\beta$ such that the series converges for $\beta > \bar\beta$ and diverges for $\beta < \bar\beta$. The behavior for $\beta = \bar\beta$ depending of further analysis of the series under consideration. This critical value is often referred to as the *abscissa of convergence*.

**8.17. Definition.** The abscissa of convergence of $Z(\beta)$ will be called the *critical inverse temperature* and will be denoted $\beta_c$. The set

$$I_c = \{\beta \in (0,\infty) : Z(\beta) < \infty\} \ \cup \ \{\infty\}$$

will be called the *interval of super-critical inverse temperatures*.

The possibilities for $I_c$ are therefore $(\beta_c,\infty]$ or $[\beta_c,\infty]$ when $\beta_c<\infty$. If $\beta_c=\infty$ then we must necessarily have $I_c=\{\infty\}$.

We therefore obtain:

**8.18. Corollary.** *For $\beta\in I_c$ every $\beta$-scaling state on $C_0(\Omega_{\mathcal{T}_A})$ is of finite type.*

## 9. Existence of finite type scaling states.

So far we have studied scaling states, and therefore KMS states on $\mathcal{T}_A$, under the assumption that they exist. In this section we shall obtain our first nontrivial existence results. Our main tool will be a parametrization of finite type scaling measures by means of their restriction to $\Omega_e$.

**9.1. Proposition.** *Let $\beta\in(0,\infty)$ and let $\lambda$ be a $\beta$-scaling measure on $\Omega_{\mathcal{T}_A}$. Then*

$$\lambda(\Omega_f) = \sum_{\mu\in P_A} N(\mu)^{-\beta}\lambda(\Omega_e^\mu).$$

*Proof.* Given $\mu\in P_A$ we have by 8.7.iv that $\lambda(\Omega_\mu) = \lambda(\alpha_\mu(\Omega_e^\mu)) = N(\mu)^{-\beta}\lambda(\Omega_e^\mu)$. The conclusion then follows from 8.7.ii. $\square$

The right hand side expression in 9.1 will be of crucial importance both for measures defined in $\Omega_e$ (observe that $\Omega_e^\mu\subseteq\Omega_e$ for all $\mu$) and for measures on $\Omega_{\mathcal{T}_A}$. This motivates the following:

**9.2. Definition.** For a measure[3] $\gamma$ defined on some measure space containing $\Omega_e$ we let

$$Z(\beta,\gamma) = \sum_{\mu\in P_A} N(\mu)^{-\beta}\gamma(\Omega_e^\mu), \quad \beta\in(0,\infty).$$

Recall from our discussion immediately after 8.6 that, for a non-trivial admissible word $\mu$, one has $\Omega_e^\mu = \Omega_e\cap\Delta_{x^{-1}}^\tau$, where $x$ is the last generator in the reduced decomposition of $\mu$. So $\Omega_e^\mu$ depends only on $x$. This said, given a measure $\gamma$ on some measure space containing $\Omega_e$, observe that

$$Z(\beta,\gamma) = \gamma(\Omega_e) + \sum_{x\in\mathcal{G}}\left(\sum_{\mu\in P_A^x} N(\mu)^{-\beta}\right)\gamma(\Omega_e^x), \quad \beta\in(0,\infty),$$

where, for each $x\in\mathcal{G}$, $P_A^x$ is the set of all admissible words ending in $x$. This motivates the introduction of our second (family of) partition function:

[3] We assume all measures are positive regular Borel measures.

**9.3. Definition.** Let $x \in \mathcal{G}$. The *fixed-target partition function* relative to the generator $x$ for the dynamical system $(\mathcal{T}_A, \sigma, \mathbb{R})$ is the function $Z_x(\beta)$ given by the Dirichlet series

$$Z_x(\beta) = \sum_{\mu \in P_A^x} N(\mu)^{-\beta}, \quad \beta \in (0, \infty).$$

For future reference we record the following:

**9.4. Proposition.** *For every measure $\gamma$ defined on some measure space containing $\Omega_e$ we have:*

$$Z(\beta, \gamma) = \gamma(\Omega_e) + \sum_{x \in \mathcal{G}} Z_x(\beta)\gamma(\Omega_e^x), \quad \beta \in (0, \infty).$$

Regarding 9.1, the observation that the series there converges and that the summands only depend on the restriction of $\lambda$ to $\Omega_e$ lead us to our next step.

**9.5. Proposition.** *Let $\beta \in (0, \infty)$ and let $\gamma$ be a measure on $\Omega_e$ such that $Z(\beta, \gamma) = 1$. Let $\lambda$ be the measure on $\Omega_{\mathcal{T}_A}$ given for every measurable subset $S \subseteq \Omega_{\mathcal{T}_A}$ by*

$$\lambda(S) = \sum_{\mu \in P_A} N(\mu)^{-\beta} \gamma\big(\alpha_{\mu^{-1}}(S \cap \Omega_\mu)\big).$$

*Then $\lambda$ is a finite type $\beta$-scaling probability measure on $\Omega_{\mathcal{T}_A}$. The correspondence $\gamma \mapsto \lambda$ gives a one-to-one affine map from the set of all measures on $\Omega_e$ such that $Z(\beta, \gamma) = 1$ onto the set of finite type $\beta$-scaling measures on $\Omega_{\mathcal{T}_A}$.*

*Proof.* Given a measurable $S \subseteq \Omega_{\mathcal{T}_A}$ observe that $S \cap \Omega_\mu \subseteq \Omega_\mu \subseteq \Delta_\mu^\tau$, which is the domain of $\alpha_{\mu^{-1}}$. Also $\alpha_{\mu^{-1}}(S \cap \Omega_\mu) \subseteq \alpha_{\mu^{-1}}(\Omega_\mu) = \Omega_e^\mu$ by 8.7.iv. Since $\Omega_e^\mu \subseteq \Omega_e$ and $\gamma$ is defined on $\Omega_e$ we see that each summand in the definition of $\lambda$ above is indeed well defined. Moreover

$$\gamma\big(\alpha_{\mu^{-1}}(S \cap \Omega_\mu)\big) \leq \gamma\big(\alpha_{\mu^{-1}}(\Omega_\mu)\big) = \gamma(\Omega_e^\mu),$$

which implies that the series defining $\lambda(S)$ is dominated by $Z(\beta, \gamma)$ and hence converges. For $S = \Omega_{\mathcal{T}_A}$ one has that

$$\lambda(S) = \sum_{\mu \in P_A} N(\mu)^{-\beta} \gamma\big(\alpha_{\mu^{-1}}(\Omega_\mu)\big) = \sum_{\mu \in P_A} N(\mu)^{-\beta} \gamma\big(\Omega_e^\mu\big) = Z(\beta, \gamma) = 1,$$

and hence $\lambda$ is indeed a probability measure. It is clearly of finite type. In order to show that $\lambda$ is $\beta$-scaling we must show that $\lambda$ satisfies

$$\lambda(\alpha_x(S)) = N(x)^{-\beta} \lambda(S)$$

for all $x \in \mathcal{G}$ and all Borel subsets $S \subseteq \Delta_{x^{-1}}^\tau$. Observing that both $\Omega_f$ and $\Omega_\infty$ are invariant under $\alpha$, and that $\lambda$ vanishes on $\Omega_\infty$, we may suppose that $S \subseteq \Omega_f$. By 8.7.ii we may in fact assume that $S \subseteq \Omega_\mu$ for some $\mu \in P_A$.

Discarding the trivial case "$S = \emptyset$" we have that

$$\emptyset \neq S \subseteq \Delta_{x^{-1}}^\tau \cap \Omega_\mu$$

and hence $x\mu$ is admissible. Moreover $\alpha_x(S) \subseteq \Omega_{x\mu}$ and

$$\lambda(\alpha_x(S)) = N(x\mu)^{-\beta} \gamma\big(\alpha_{(x\mu)^{-1}}(\alpha_x(S))\big) = N(x)^{-\beta} N(\mu)^{-\beta} \gamma\big(\alpha_{\mu^{-1}}(S)\big) =$$

$$= N(x)^{-\beta} \lambda(S).$$

It is clear that the restriction of $\lambda$ to $\Omega_e$ coincides with $\gamma$ and hence our correspondence is injective. On the other hand given any finite type $\beta$-scaling measure $\lambda$ on $\Omega_{\mathcal{T}_A}$ it is an easy exercise to show that the restriction of $\lambda$ to $\Omega_e$, say $\gamma$, is a measure that satisfies $Z(\beta, \gamma) = 1$ and is mapped to $\lambda$ under our correspondence. This proves surjectivity. Finally it is clear that we have an affine map. ☐

Suppose we are given a nonzero measure $\gamma$ on $\Omega_e$ such that $Z(\beta,\gamma) < \infty$. Note that such a measure must necessarily be finite because $\gamma(\Omega_e) \le Z(\beta,\gamma)$. By normalizing it we obtain a measure

$$\gamma' = \frac{1}{Z(\beta,\gamma)}\,\gamma$$

that satisfies the hypothesis of 9.5 and hence gives rise to a $\beta$-scaling measure. Of course many different $\gamma$'s are mapped to the same $\lambda$ but this happens if and only if the $\gamma$'s involved are multiples of each other.

Even if $Z(\beta,\gamma) = \infty$, actually even if $\gamma$ is an infinite measure, one could attempt to define a $\beta$-scaling *infinite* measure on $\Omega_{\mathcal{T}_A}$ using the method of 9.5. This combined with a likely generalization of 8.2 for infinite measures would perhaps lead to interesting KMS *weights* on $\mathcal{T}_A$. However we will not pursue these ideas in the present work.

**9.6. Definition.** Given $\beta \in (0,\infty)$ and a nonzero measure $\gamma$ on $\Omega_e$ such that $Z(\beta,\gamma) < \infty$ we will denote by $T_\beta(\gamma)$ the finite type $\beta$-scaling (probability) measure $\lambda$ obtained by applying the construction of 9.5 to $\gamma/Z(\beta,\gamma)$. If $\beta = \infty$ and $\gamma$ is any nonzero finite measure on $\Omega_e$ we will let $T_\beta(\gamma)$ be the measure on $\Omega_{\mathcal{T}_A}$ given simply by

$$T_\beta(\gamma)(S) = \frac{\gamma(S\cap\Omega_e)}{\gamma(\Omega_e)}$$

for all Borel subsets $S \subseteq \Omega_{\mathcal{T}_A}$.

Recall that $I_c$ is the interval of convergence of the Dirichlet series $\sum_{\mu\in P_A} N(\mu)^{-\beta}$. Given $\beta \in I_c$ note that the convergence of that series implies the convergence of

$$\sum_{\mu\in P_A} N(\mu)^{-\beta}\gamma(\Omega_e^\mu),$$

which defines $Z(\beta,\gamma)$, irrespective of which finite measure $\gamma$ we have in mind. This says that $Z(\beta,\gamma) < \infty$ and hence that $T_\beta(\gamma)$ is defined for every nonzero finite measure $\gamma$ on $\Omega_e$.

Combining what we have just found with 8.18 we obtain a complete characterization of $\beta$-scaling measures on $\Omega_{\mathcal{T}_A}$ (and hence also of KMS$_\beta$ states on $\mathcal{T}_A$ by 8.2) in the interval of super-critical inverse temperatures:

**9.7. Theorem.** *Under 8.1 let $\beta \in I_c$. Then the correspondence $\gamma \mapsto T_\beta(\gamma)$ establishes a surjective map from the set of nonzero finite measures $\gamma$ on $\Omega_e$ to the set of $\beta$-scaling measures on $\Omega_{\mathcal{T}_A}$, all of which are of finite type. This correspondence is not injective but $T_\beta(\gamma_1) = T_\beta(\gamma_2)$ if and only if $\gamma_1$ is a multiple of $\gamma_2$.*

## 10. Irreducible matrices and the fixed-target partition function $Z_y(\beta)$.

From now on we shall occasionally make a few other hypotheses, in addition to 8.1, which should perhaps be listed here for ease of reference:

**10.1. Occasional Hypotheses.**

(IRR) $A$ is *irreducible*, i.e. for every $x$ and $y$ in $\mathcal{G}$ there exists an admissible word $\mu$ with $\mu_1 = x$ and $\mu_{|\mu|} = y$.

(COL) $A$ has no identically zero *columns,*

(FTS) There exists a *finite target set*, i.e. a finite set $\{y_1,\ldots,y_n\} \subseteq \mathcal{G}$ such that for every $x \in \mathcal{G}$ one has $A(x,y_i) = 1$ for at least one $i$.

(INF) The $N(x)$'s are *bounded away from* 1 in the sense that $\inf_{x\in\mathcal{G}} N(x) > 1$.

Except for the implication "(IRR) $\Rightarrow$ (COL)", which is easy to verify, there are no other logical relations between the conditions of 10.1.

It should be stressed that we are assuming 8.1 throughout, often without notice, but we will be explicit when assuming any one of the "occasional hypotheses" above.

In some cases there is a close relationship between the convergence of the series for the various $Z_x(\beta)$ which we would like to present now.

**10.2. Proposition.** *Let $x, y \in \mathcal{G}$. Suppose that there exists an admissible word $\nu \in P_A$ beginning in $x$ and ending in $y$. Then for every $\beta \in (0, \infty)$ one has that*

$$Z_x(\beta) \leq N(x^{-1}\nu)^{\beta} Z_y(\beta).$$

*Proof.* Considering the (obviously injective) map $\mu \in P_A^x \mapsto \mu x^{-1}\nu \in P_A^y$ we have

$$Z_y(\beta) = \sum_{\mu \in P_A^y} N(\mu)^{-\beta} \geq \sum_{\mu \in P_A^x} N(\mu x^{-1}\nu)^{-\beta} =$$

$$= N(x^{-1}\nu)^{-\beta} \sum_{\mu \in P_A^x} N(\mu)^{-\beta} = N(x^{-1}\nu)^{-\beta} Z_x(\beta). \qquad \square$$

In particular, under the conditions above, if the series for $Z_y(\beta)$ converges then so does the series for $Z_x(\beta)$. The case in which this relationship is richest is when $A$ is irreducible (see 10.1.(IRR)), in which case we get the following "solidarity" result for our Dirichlet series:

**10.3. Proposition.** *Let $A$ be an irreducible matrix. Then for every $\beta \in (0, \infty)$ one has that either*

- *$Z_x(\beta) < \infty$ for all $x \in \mathcal{G}$, or*
- *$Z_x(\beta) = \infty$ for all $x \in \mathcal{G}$.*

*Proof.* Follows immediately from 10.2. $\square$

We will therefore assume throughout this section that $A$ is an irreducible matrix. It follows that the set of $\beta$'s for which $Z_x(\beta) < \infty$ does not depend on $x$, motivating the following:

**10.4. Definition.** If $A$ is irreducible the abscissa of convergence for each and every one of the Dirichlet series $Z_x(\beta)$ will be called the *fixed-target critical inverse temperature* and will be denoted $\dot{\beta}_c$. The set of $\beta$'s where each and every one of these series converge, including $\beta = \infty$, will be called the *interval of fixed-target super-critical inverse temperatures* and will be denoted $\dot{I}_c$.

As before $\dot{I}_c$ can be either one of $(\dot{\beta}_c, \infty]$ or $[\dot{\beta}_c, \infty]$ when $\dot{\beta}_c < \infty$, and $\dot{I}_c = \{\infty\}$ when $\dot{\beta}_c = \infty$.

Since $Z_x(\beta)$ is defined as a subseries of $Z(\beta)$ it is obvious that the convergence of the latter implies the convergence of the former. This gives:

**10.5. Proposition.** *One has that $\dot{\beta}_c \leq \beta_c$ and $\dot{I}_c \supseteq I_c$.*

Recall that 8.18 says that for $\beta \in I_c$ every $\beta$-scaling measure is of finite type. Our next result goes in the opposite direction stating that there are no finite type $\beta$-scaling measure for $\beta \notin \dot{I}_c$.

**10.6. Theorem.** *Under 8.1 assume that $A$ is irreducible and let $\beta \notin \dot{I}_c$. Then every $\beta$-scaling measure on $\Omega_{\mathcal{T}_A}$ is of infinite type. Consequently there are no finite type $\beta$-scaling measures.*

*Proof.* Let $\lambda$ be a $\beta$-scaling measure on $\Omega_{\mathcal{T}_A}$. Then by 9.1 and 9.4 we have

$$\lambda(\Omega_f) = \lambda(\Omega_e) + \sum_{x \in \mathcal{G}} Z_x(\beta) \lambda(\Omega_e^x).$$

Given that $\beta \notin \dot{I}_c$ we have that $Z_x(\beta) = \infty$ for all $x$ implying that $\lambda(\Omega_e^x) = 0$ and hence that $\lambda(\Omega_f) = \lambda(\Omega_e)$. As we are working under the assumption that $\mathcal{G}$ is countable we also have that

$$0 = \lambda\left(\bigcup_{x\in\mathcal{G}} \Omega_e^x\right) = \lambda\left(\bigcup_{x\in\mathcal{G}} \Omega_e \cap \Delta^\tau_{x^{-1}}\right) = \lambda\left(\Omega_e \cap \bigcup_{x\in\mathcal{G}} \Delta^\tau_{x^{-1}}\right).$$

We claim that $\Omega_e \subseteq \bigcup_{x\in\mathcal{G}} \Delta^\tau_{x^{-1}}$. To prove it assume by contradiction that $\xi \in \Omega_e$ but $\xi \notin \Delta^\tau_{x^{-1}}$ for all $x$ in $\mathcal{G}$. Then $x^{-1} \notin \xi$ for all $x$, which gives $R_\xi(e) = \emptyset$. Since $\xi \in \Omega_e$ we have that the stem of $\xi$ is $e$. Using [**EL**: 5.12] we conclude that $\xi = \epsilon$. But this is a contradiction since $\Omega_e \subseteq \Omega_{\mathcal{T}_A} = \widetilde{\Omega}_{\mathcal{T}_A} \setminus \{\epsilon\}$. This proves our claim and hence that $0 = \lambda(\Omega_e) = \lambda(\Omega_f)$. ☐

The following diagram subsumes the information about $\beta$-scaling states on $C_0(\Omega_{\mathcal{T}_A})$, and hence also about KMS$_\beta$ states on $\mathcal{T}_A$, that we have gathered so far in the case of an irreducible matrix $A$.

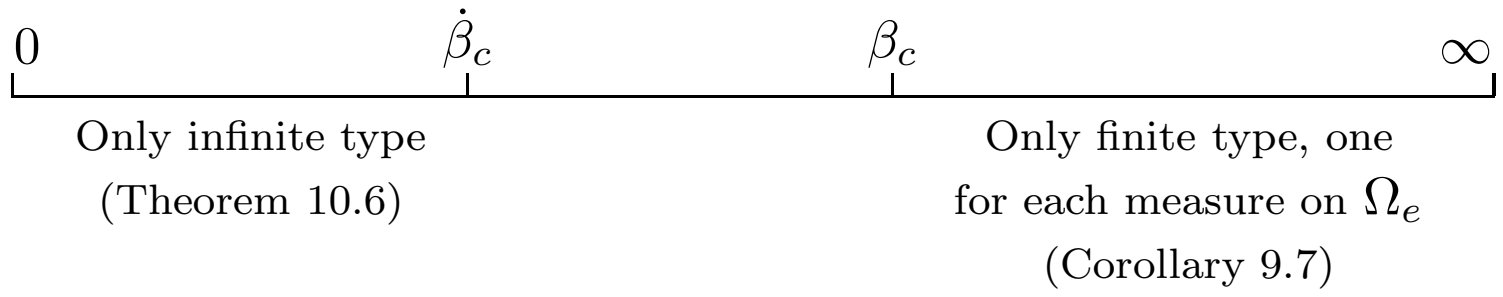

*Diagram 10.7*

With this we essentially exhaust the conclusions that can be drawn from the techniques developed so far. In order to proceed further we need a characterization of $\beta$-scaling measures which, unlike 9.5, includes both finite and infinite type measures.

## 11. The structure of $\widetilde{\mathcal{T}}_A$.

We retain, as always, the hypotheses listed in 8.1. In this section it will be convenient to deal with unital algebras and hence we will mainly consider $\widetilde{\mathcal{T}}_A$ as opposed to $\mathcal{T}_A$.

Our major desire is to describe, for each inverse temperature $\beta$, the KMS$_\beta$ states of $\widetilde{\mathcal{T}}_A$, which we will do by characterizing the simplex formed by all $\beta$-scaling probability measures on $\widetilde{\Omega}_{\mathcal{T}_A}$. As an intermediate goal we will show that these measures are parametrized by certain states on the algebra $\widetilde{\mathcal{Q}}$ defined in 6.2. In preparation for this we will now dive into the study of this and other subalgebras of $\widetilde{\mathcal{T}}_A$.

Recall that the space $\Sigma_A$, introduced shortly before 5.3, is the closure of the set $\{c_x : x \in \mathcal{G}\}$, formed by the columns $c_x$ of $A$, within the topological Cantor space $2^{\mathcal{G}}$,

**11.1. Proposition.** *Consider the map $R : \widetilde{\Omega}_{\mathcal{T}_A} \to \Sigma_A$ given by $R(\xi) = R_\xi(e)$ and let $r : \widetilde{\Omega}_{\mathcal{T}_A} \to \widetilde{\Omega}_e$ be given by 5.4. Then:*

(i) *There exists a homeomorphism $h : \Sigma_A \to \widetilde{\Omega}_e$ such that the diagram*

$$\begin{array}{ccc} \widetilde{\Omega}_{\mathcal{T}_A} & \stackrel{r}{\longrightarrow} & \widetilde{\Omega}_e \\ {\scriptstyle R}\downarrow & \nearrow_h & \\ \Sigma_A & & \end{array}$$

*commutes.*

(ii) *Let $\hat{R} : C(\Sigma_A) \to C(\widetilde{\Omega}_{\mathcal{T}_A})$ and $\hat{r} : C(\widetilde{\Omega}_e) \to C(\widetilde{\Omega}_{\mathcal{T}_A})$ be obtained by transposing $R$ and $r$, respectively. Then the range of both $\hat{R}$ and $\hat{r}$ coincide with $\widetilde{\mathcal{Q}}$.*

(iii) *Both $\hat{R}$ and $\hat{r}$ are isomorphisms onto $\widetilde{\mathcal{Q}}$ and hence*

$$\widetilde{\mathcal{Q}} \simeq C(\Sigma_A) \simeq C(\widetilde{\Omega}_e).$$

(iv) *For every $a \in \widetilde{\mathcal{Q}}$ and every $\xi \in \widetilde{\Omega}_{\mathcal{T}_A}$ one has that $a(\xi) = a(r(\xi))$.*

*Proof.* We claim that for $\xi, \eta \in \widetilde{\Omega}_{\mathcal{T}_A}$ one has that

$$R(\xi) = R(\eta) \ \Leftrightarrow\ r(\xi) = r(\eta).$$

On the one hand by definition of $r$ one has that $R(\xi) = R(\eta) \ \Rightarrow\ r(\xi) = r(\eta)$. On the other hand, observing that $R(r(\xi)) = R(\xi)$, we see that $r(\xi) = r(\eta) \ \Rightarrow\ R(\xi) = R(\eta)$. Since both $R$ and $r$ are clearly surjective, a bijection $h$ exists such that $h \circ R = r$. By compactness of $\widetilde{\Omega}_{\mathcal{T}_A}$ both $R$ and $r$ are quotient maps and hence $h$ is a homeomorphism.

By 6.1.iii we have $q_x = 1_{\Delta^{\tau}_{x^{-1}}}$ so that

$$q_x(\xi) = \big[x^{-1} \in \xi\big], \quad \forall \xi \in \widetilde{\Omega}_{\mathcal{T}_A}.$$

For $a \in \widetilde{\mathcal{Q}}$ it follows that the value of $a(\xi)$ depends only on $\{x \in \mathcal{G} : x^{-1} \in \xi\} = R_\xi(e) = R(\xi)$ in the sense that for $\xi$ and $\eta$ in $\widetilde{\Omega}_{\mathcal{T}_A}$

$$R(\xi) = R(\eta) \quad\Longrightarrow\quad \big(\forall a \in \widetilde{\mathcal{Q}} \quad a(\xi) = a(\eta)\big).$$

This immediately implies (iv) in view of the fact that $R(r(\xi)) = R(\xi)$.

The converse of the above implication also holds, as it can be proved by considering $a = q_x$. Therefore the equivalence relation defined on $\widetilde{\Omega}_{\mathcal{T}_A}$ by $R$ (i.e. having the same image under $R$) coincides with the equivalence relation defined by $\widetilde{\mathcal{Q}}$ (i.e. having the same image under every $a \in \widetilde{\mathcal{Q}}$). These in turn also coincide with the relation defined by $r$, whence (ii).

Since both $R$ and $r$ are surjective we have that both $\hat{R}$ and $\hat{r}$ are injective therefore proving (iii). ☐

We will later need a technical result about approximating positive elements of $\widetilde{\mathcal{Q}}$ which we would now like to present. Let $\xi \in \widetilde{\Omega}_{\mathcal{T}_A}$ correspond to $c \in \Sigma_A$ under $R$ (that is $R(\xi) = c$) and observe that for all $x \in \mathcal{G}$ $q_x(\xi) = \big[x^{-1} \in \xi\big] = \big[x \in c\big]$. Identifying $\widetilde{\mathcal{Q}}$ with $C(\Sigma_A)$ via $\hat{R}$ we may therefore think of $q_x$ as the function

$$q_x(c) = \big[x \in c\big]. \tag{11.2}$$

Given finite subsets $X$ and $Y$ of $\mathcal{G}$ it follows that $q(X, Y)$ is in turn identified with the characteristic function of the set

$$V(X, Y) = \{c \in \Sigma_A : x \in c,\ y \notin c,\ \forall x \in X,\ \forall y \in Y\}.$$

Observe also that these sets form a basis for the product topology on $2^{\mathcal{G}}$, consisting of clopen sets.

**11.3. Lemma.** *For each $a \geq 0$ in $\widetilde{\mathcal{Q}}$ and any $\varepsilon > 0$ there are finite subsets $X_1, \ldots, X_n$ and $Y_1, \ldots, Y_n$ of $\mathcal{G}$ and positive real numbers $\lambda_1, \ldots, \lambda_n$ such that the element*

$$b = \sum_{i=1}^{n} \lambda_i q(X_i, Y_i)$$

*satisfies*

(i) $0 \leq b \leq a$, *and*
(ii) $\|a - b\| < \varepsilon$.

*Proof.* Let $K = \{c \in \Sigma_A : a(c) \geq 2\|a\|/3\}$ and $U = \{c \in \Sigma_A : a(c) > \|a\|/3\}$ so that $K$ is compact, $U$ is an open set, and $K \subseteq U \subseteq \Sigma_A$. Choose a finite covering of $K$ consisting of sets $V(X_i, Y_i) \subseteq U$. It is not hard to show that such a covering can be found so that the $V(X_i, Y_i)$ are pairwise disjoint. Let

$$b_1 = \frac{\|a\|}{3} \sum_{i=1}^{n} q(X_i, Y_i).$$

It is now easy to show that $0 \leq b_1 \leq a$ and that $\|a - b_1\| \leq 2\|a\|/3$. We may then repeat this procedure starting with $a - b_1$ and, after $n$ steps, we will have obtained a sequence $b_1, \ldots, b_n$ of elements of $\widetilde{\mathcal{Q}}$, each of which is a scalar multiple of a sum of $q(X, Y)$'s, and such that $0 \leq b_n \leq a - b_1 - \cdots - b_{n-1}$ and $\|a - b_1 - \cdots - b_{n-1} - b_n\| \leq (2/3)^n \|a\|$. After a finite number of steps the element $b = b_1 + \cdots + b_n$ will satisfy the required properties. ☐

We would now like to study other subalgebras of $\widetilde{\mathcal{T}}_A$.

**11.4. Proposition.** *For each $\mu \in \mathbb{F}_+$ let $\mathcal{I}^\mu = S(\mu)\widetilde{\mathcal{Q}}S(\mu)^*$. Then*

(i) *If $\mu, \nu \in \mathbb{F}_+$ are such that $|\mu| = |\nu|$ but $\mu \neq \nu$ then $\mathcal{I}^\mu\mathcal{I}^\nu = \{0\}$.*

(ii) *If $\mu, \nu \in \mathbb{F}_+$ are such that $|\mu| \leq |\nu|$ then $\mathcal{I}^\mu\mathcal{I}^\nu \subseteq \mathcal{I}^\nu$.*

(iii) *Each $\mathcal{I}^\mu$ is a closed unital *-subalgebra of $\widetilde{\mathcal{T}}_A$ and $S(\mu)S(\mu)^*$ is its unit.*

*Proof.* Part (i) follows at once from 6.3.i, while (ii) is just a restatement of 6.3.iv. That each $\mathcal{I}^\mu$ is a *-algebra follows from (ii) and it is obvious that $S(\mu)S(\mu)^*$ serves as a unit for it. It therefore remains to show that $\mathcal{I}^\mu$ is closed. So let a sequence $\{S(\mu)a_nS(\mu)^*\}_n$, with $a_n \in \widetilde{\mathcal{Q}}$, converge to some $b$ in $\widetilde{\mathcal{T}}_A$. Then

$$b = \lim_n S(\mu)a_nS(\mu)^* = \lim_n S(\mu)S(\mu)^*S(\mu)a_nS(\mu)^*S(\mu)S(\mu)^* =$$

$$= S(\mu)S(\mu)^*bS(\mu)S(\mu)^*,$$

and hence it suffices to show that $S(\mu)^*bS(\mu)$ is in $\widetilde{\mathcal{Q}}$. But

$$S(\mu)^*bS(\mu) = \lim_n S(\mu)^*S(\mu)a_nS(\mu)^*S(\mu),$$

which belongs to $\widetilde{\mathcal{Q}}$ by 6.3.ii. □

From now on $\mathbb{F}_+^n$ will denote the subset of $\mathbb{F}_+$ consisting of elements $\mu$ with $|\mu| = n$. Thus $\mathbb{F}_+^0 = \{e\}$ and $\mathbb{F}_+^1 = \mathcal{G}$.

**11.5. Proposition.** *For each integer $n \geq 0$ let $\mathcal{I}_n$ be the closure of $\bigoplus_{\mu\in\mathbb{F}_+^n} \mathcal{I}^\mu$ within $\widetilde{\mathcal{T}}_A$. Then $\mathcal{I}_n$ is a closed *-subalgebra of $\widetilde{\mathcal{T}}_A$ which is *-isomorphic to the $c_0$ direct sum of the $\mathcal{I}^\mu$, that is, the $C^*$-algebra consisting of families $(a_\mu)_{\mu\in\mathbb{F}_+^n}$ such that $\lim_\mu \|a_\mu\| = 0$. In particular the net of idempotents $\{\sum_{\mu\in J} p_\mu\}_J$, where $J$ ranges in the collection of finite subset of $\mathbb{F}_+^n$, forms an approximate unit for $\mathcal{I}_n$.*

*Proof.* The statement follows easily from the fact that the $\mathcal{I}^\mu$ considered form a collection of pairwise orthogonal $C^*$-algebras by 11.4.i. □

**11.6. Proposition.**

(i) *For every $n$ and $m$ one has that $\mathcal{I}_n\mathcal{I}_m \subseteq \mathcal{I}_{\max\{n,m\}}$.*

(ii) *For $n \geq 1$ and $x \in \mathcal{G}$ one has that $s_x^*\mathcal{I}_ns_x \subseteq \mathcal{I}_{n-1}$.*

(iii) *For $n \geq 0$ and $x \in \mathcal{G}$ one has that $s_x\mathcal{I}_ns_x^* \subseteq \mathcal{I}_{n+1}$.*

*Proof.* The first statement follows from 11.4.ii. As for (ii) let $a \in \widetilde{\mathcal{Q}}$ and $\mu \in \mathbb{F}_+^n$. We then need to show that $s_x^*S(\mu)aS(\mu)^*s_x$ lies in $\mathcal{I}_{n-1}$. Let $y$ be the first generator in the reduced decomposition of $\mu$ so that $\mu = y\mu'$, where $\mu' \in \mathbb{F}_+$.

Observe that unless $x = y$ we have that $s_x^*S(\mu) = 0$. So assume that $x = y$.

The case $n = 1$ is somewhat special so let us treat it first. We then have that $\mu = x$ and thus

$$s_x^*S(\mu)aS(\mu)^*s_x = q_xaq_x \in \widetilde{\mathcal{Q}} = \mathcal{I}_0.$$

If $n \geq 2$ then $|\mu'| \geq 1$ and hence

$$s_x^*S(\mu) = q_xS(\mu') = \varepsilon S(\mu'),$$

where $\varepsilon \in \{0, 1\}$ by $\mathrm{CK}_3$. Therefore

$$s_x^*S(\mu)aS(\mu)^*s_x = \varepsilon S(\mu')aS(\mu')^* \in \mathcal{I}_{n-1}.$$

The third assertion is obvious. □

**11.7. Proposition.** *For each integer $n \geq 0$ let $A_n$ be the closure of $\mathcal{I}_0 + \cdots + \mathcal{I}_n$ within $\widetilde{\mathcal{T}}_A$. Then*

(i) *$A_n$ is a $C^*$-algebra,*

(ii) *$\mathcal{I}_n$ is an ideal in $A_n$,*

(iii) *$A_{n+1} = A_n + \mathcal{I}_{n+1}$,*

(iv) *$C(\widetilde{\Omega}_{\mathcal{T}_A})$ is the closure of $\cup_n A_n$.*

(v) *For $n \geq 1$ and $x \in \mathcal{G}$ one has that $s_x^* A_n s_x \subseteq A_{n-1}$.*

(vi) *For $n \geq 0$ and $x \in \mathcal{G}$ one has that $s_x A_n s_x^* \subseteq A_{n+1}$.*

*Proof.* Clearly (i) and (ii) follow from 11.6.i. Using [**Pe**: 1.5.8] one gets (iii). As for (iv), it follows from 6.4. Finally (v) follows from a combination of 11.6.ii and 6.3.iii while (vi) is a direct consequence of 11.6.iii. ☐

One more technical result is in order:

**11.8. Proposition.** *For each $n \in \mathbb{N}$ one has that $A_n \cap \mathcal{I}_{n+1} = \{0\}$.*

*Proof.* Assume first that $n = 0$, observing that $A_0 = \widetilde{\mathcal{Q}}$. Let $a \in \widetilde{\mathcal{Q}} \cap \mathcal{I}_1$. Using 11.5 write $a = \sum_{z\in\mathcal{G}} a_z$ where $a_z \in \mathcal{I}^z$ and $\lim_z \|a_z\| = 0$. We claim that each $a_z$ is a scalar multiple of $p_z$. In fact observe that, since $a \in \widetilde{\mathcal{Q}}$, we have by 6.3.iii that $s_z^* a s_z = \lambda_z q_z$ for some $\lambda_z \in \mathbb{C}$. Therefore

$$a_z = p_z a p_z = s_z(s_z^* a s_z)s_z^* = \lambda_z s_z q_z s_z^* = \lambda_z p_z.$$

We then have that $a = \sum_{z\in\mathcal{G}} \lambda_z p_z$ with $\lim_z \lambda_z = 0$.

Suppose by way of contradiction that $a \neq 0$. Then there exists at least one $\lambda_{z_0} p_{z_0}$ which is nonzero. Given that $\lim_z \lambda_z = 0$ we see that $\lambda_{z_0}$ is an isolated point in the set of all $\lambda_z$'s. One may therefore take a continuous function $f : \mathbb{C} \to \mathbb{C}$ such that $f(\lambda_{z_0}) = 1$ and $f(\lambda_z) = 0$ whenever $\lambda_z \neq \lambda_{z_0}$. It follows that $f(a) = \sum_{z\in Z} p_z$, where $Z$ is the (necessarily finite) set $Z = \{z \in \mathcal{G} : \lambda_z = \lambda_{z_0}\}$. It is also clear that $0 \neq f(a) \in \widetilde{\mathcal{Q}} \cap \mathcal{I}_1$.

Given $\xi \in \widetilde{\Omega}_{\mathcal{T}_A}$ we have, using 11.1.iv, that

$$f(a)\big|_\xi = f(a)\big|_{r(\xi)} = \sum_{z\in Z} p_z(r(\xi)) = \sum_{z\in Z} \big[z \in r(\xi)\big] = 0,$$

because the stem of $r(\xi)$ is trivial. It follows that $f(a) = 0$, a contradiction.

Now assume that $n \geq 1$ and let $a \in A_n \cap \mathcal{I}_{n+1}$. Given $\nu \in \mathbb{F}_+^n$ we have that $S(\nu)^* a S(\nu) \in A_0 \cap \mathcal{I}_1$ by 11.6.ii and 11.7.v and hence $S(\nu)^* a S(\nu) = 0$. With more reason $S(\mu)^* a S(\mu) = 0$ for $\mu \in \mathbb{F}_+^{n+1}$. Therefore

$$p_\mu a = p_\mu a p_\mu = S(\mu)S(\mu)^* a S(\mu)S(\mu)^* = 0$$

for every $\mu \in \mathbb{F}_+^{n+1}$. The conclusion then follows from the last sentence in 11.5. ☐

## 12. Invariant and subinvariant states on $\widetilde{\mathcal{Q}}$.

Our goal in this section is to show that $\beta$-scaling states on $C(\widetilde{\Omega}_{\mathcal{T}_A})$, and hence also $\mathrm{KMS}_\beta$ states on $\widetilde{\mathcal{T}}_A$, are in 1–1 correspondence with certain states on $\widetilde{\mathcal{Q}}$. These states are best motivated by the following:

**12.1. Proposition.** *Let $\beta \in (0,\infty)$ and let $\phi$ be a $\beta$-scaling state on $C(\widetilde{\Omega}_{\mathcal{T}_A})$. Denote by $\rho$ the restriction of $\phi$ to $\widetilde{\mathcal{Q}}$. Then, for every pair of finite subsets $X$ and $Y$ of $\mathcal{G}$, we have*

$$\sum_{z\in\mathcal{G}} A(X,Y,z)N(z)^{-\beta}\rho(q_z) \le \rho(q(X,Y)).$$

*Proof.* Recall from $\mathrm{CK}_3$ that $q_x s_z = A(x,z)s_z$ and hence also $(1-q_y)s_z = (1-A(y,z))s_z$ so that

$$q(X,Y)s_z = \left(\prod_{x\in X} q_x \prod_{y\in Y}(1-q_y)\right) s_z = \left(\prod_{x\in X} A(x,z) \prod_{y\in Y}(1-A(y,z))\right) s_z = A(X,Y,z)s_z$$

for all $z$ in $\mathcal{G}$. By multiplying this on the right hand side by $s_z^*$ we have that $q(X,Y)p_z = A(X,Y,z)p_z$ and hence that $A(X,Y,z)p_z \le q(X,Y)$ in the usual order of projections. Since the $p_z$ are pairwise orthogonal by $\mathrm{CK}_2$ we conclude that any finite sum $\sum_{z\in Z} A(X,Y,z)p_z$ ($Z$ a finite set) gives a projection dominated by $q(X,Y)$. It follows that

$$\phi\left(\sum_{z\in Z} A(X,Y,z)p_z\right) \le \phi(q(X,Y)).$$

On the other hand recall from 8.5 that $\phi(p_z) = N(z)^{-\beta}\phi(q_z)$. Therefore

$$\sum_{z\in Z} A(X,Y,z)N(z)^{-\beta}\rho(q_z) \le \rho(q(X,Y)).$$

Since $Z$ is arbitrary, the proof is concluded. ☐

It should be noted that the result above covers the case $X = Y = \emptyset$, in which case it says that

$$\sum_{z\in\mathcal{G}} N(z)^{-\beta}\rho(q_z) \le 1. \tag{12.2}$$

The states $\rho$ appearing above will acquire a crucial importance from this point on and hence we make the following:

**12.3. Definition.** Let $\beta \in (0,\infty)$. A state $\rho$ on $\widetilde{\mathcal{Q}}$ is said to be

(i) *$\beta$-subinvariant* when the inequality in 12.1 holds for all finite subsets $X, Y \subseteq \mathcal{G}$.

(ii) *$\beta$-invariant* when the inequality in 12.1 becomes an equality for all finite subsets $X, Y \subseteq \mathcal{G}$.

A probability measure on $\Sigma_A$ is said to be $\beta$-subinvariant (resp. $\beta$-invariant) if integration against it leads to a $\beta$-subinvariant (resp. $\beta$-invariant) state on $C(\Sigma_A) = \widetilde{\mathcal{Q}}$. Every state or measure will be considered $\infty$-subinvariant by default.

Our last result therefore says that the correspondence $\phi \mapsto \phi|_{\widetilde{\mathcal{Q}}}$ maps the set of $\beta$-scaling states on $C(\widetilde{\Omega}_{\mathcal{T}_A})$ to the set of $\beta$-subinvariant states on $\Sigma_A$. We will now seek to prove that this is in fact a bijective correspondence, thus obtaining a new characterization of KMS states which is significantly better than the one obtained in 8.2 in the sense that $\Sigma_A$ is a much more tractable space than $\widetilde{\Omega}_{\mathcal{T}_A}$.

We begin by proving that $\phi \mapsto \phi|_{\widetilde{\mathcal{Q}}}$ defines an injective map.

**12.4. Proposition.** *Let $\beta \in (0,\infty]$ and let $\phi$ and $\phi'$ be $\beta$-scaling states on $C(\widetilde{\Omega}_{\mathcal{T}_A})$ such that $\phi|_{\widetilde{\mathcal{Q}}} = \phi'|_{\widetilde{\mathcal{Q}}}$. Then $\phi = \phi'$.*

*Proof.* We first claim that $\phi$ and $\phi'$ coincide on elements of the form $S(\mu)aS(\mu)^*$, where $\mu \in \mathbb{F}_+$ and $a \in \widetilde{\mathcal{Q}}$. Using 6.1.v we have that

$$\phi\big(S(\mu)aS(\mu)^*\big) = \phi\big(\theta_\mu(q_\mu a)\big) = N(\mu)^{-\beta}\phi(q_\mu a).$$

Since $q_\mu a$ is in $\widetilde{\mathcal{Q}}$ by 6.3.ii the claim is proven. By 6.4 it follows that $\phi$ and $\phi'$ coincide on $C(\widetilde{\Omega}_{\mathcal{T}_A})$. □

In order to prove that the correspondence $\phi \mapsto \phi|_{\widetilde{\mathcal{Q}}}$ is surjective we need the following general result about states.

**12.5. Proposition.** *Let $B$ be a unital $C^*$-algebra containing a closed two sided ideal $I$ and a sub-$C^*$-algebra $A$ such that $1 \in A$ and $B = A + I$. Also let $\phi$ be a state on $A$ and $\psi$ be a positive linear functional on $I$. Denote by $\widetilde{\psi}$ the canonical extension of $\psi$ to a positive functional on $B$ (that is, $\widetilde{\psi}(b) = \lim_i \psi(bu_i)$, where $\{u_i\}_i$ is an approximate unit for $I$). Suppose that*

(i) *$\phi \geq \widetilde{\psi}$ on $A$, and*
(ii) *$\phi = \psi$ on $A \cap I$.*

*Then there exists a state $\rho$ on $B$ such that $\rho|_A = \phi$ and $\rho|_I = \psi$.*

*Proof.* Given $b$ in $B$ write $b = a + x$, where $a \in A$ and $x \in I$, and put $\rho(b) = \phi(a) + \psi(x)$. It follows from (ii) that $\rho$ is a well defined linear functional on $B$. In order to show that $\rho$ is positive let $b = a + x \in B$ and observe that

$$\rho(b^*b) = \rho(a^*a + a^*x + x^*a + x^*x) = \phi(a^*a) + \psi(a^*x + x^*a + x^*x) \geq$$

$$\geq \widetilde{\psi}(a^*a) + \psi(a^*x + x^*a + x^*x) = \widetilde{\psi}(b^*b) \geq 0.$$

Since $1 \in A$ we have $\|\rho\| = \rho(1) = \phi(1) = 1$ and hence $\rho$ is indeed a state. □

With the following result we complete the announced parametrization of $\beta$-scaling states on $C(\widetilde{\Omega}_{\mathcal{T}_A})$ by means of $\beta$-subinvariant states on $\widetilde{\mathcal{Q}}$.

**12.6. Proposition.** *Let $\beta \in (0,\infty]$ and let $\rho$ be a $\beta$-subinvariant state on $\widetilde{\mathcal{Q}}$. Then there exists a (necessarily unique) $\beta$-scaling state $\phi$ on $C(\widetilde{\Omega}_{\mathcal{T}_A})$ such that $\phi|_{\widetilde{\mathcal{Q}}} = \rho$.*

*Proof.* We begin with the case $\beta < \infty$. For each $n \in \mathbb{N}$ we will construct a state $\rho_n$ on the algebra $A_n$ (see 11.7) such that for every $n \geq 1$,

(i) $\rho_0 = \rho$,
(ii) $\rho_n|_{A_{n-1}} = \rho_{n-1}$,
(iii) $\rho_n(s_x a s_x^*) = N(x)^{-\beta}\rho_{n-1}(a s_x^* s_x)$ for $a \in A_{n-1}$ and $x \in \mathcal{G}$,
(iv) $N(x)^{-\beta}\rho_{n-1}\big(s_x^* a s_x\big) = \rho_n(a s_x s_x^*)$ for $a \in A_n$ and $x \in \mathcal{G}$.

We shall proceed by induction and hence let us suppose we are given $m \geq 0$ and $\{\rho_n\}_{0\leq n\leq m}$ satisfying (i–iv) for all $n = 1,\ldots,m$. Define a linear functional $\chi_{m+1}$ on $\mathcal{I}_{m+1}$ by

$$\chi_{m+1}(a) = \sum_{x\in\mathcal{G}} N(x)^{-\beta}\rho_m(s_x^* a s_x), \quad a \in \mathcal{I}_{m+1}.$$

In order to verify that this is well defined observe that, by 11.6.ii, for any $a \in \mathcal{I}_{m+1}$ we have $s_x^* a s_x \in s_x^*\mathcal{I}_{m+1}s_x \subseteq \mathcal{I}_m \subseteq A_m$ so that $\rho_m(s_x^* a s_x)$ is defined. To see that the sum converges it is enough to consider a positive $a$, in which case we have

$$\sum_{x\in\mathcal{G}} N(x)^{-\beta}\rho_m(s_x^* a s_x) \leq \sum_{x\in\mathcal{G}} N(x)^{-\beta}\rho_m(\|a\|q_x) = \|a\|\sum_{x\in\mathcal{G}} N(x)^{-\beta}\rho(q_x) \leq \|a\|,$$

where the last step follows from 12.2. It is then clear that $\chi_{m+1}$ is a well defined positive linear functional on $\mathcal{I}_{m+1}$. By 11.8 we have that $A_m \cap \mathcal{I}_{m+1} = \{0\}$ and hence the expression

$$\rho_{m+1}(a+b) = \rho_m(a) + \chi_{m+1}(b), \quad a \in A_m,\ b \in \mathcal{I}_{m+1}$$

gives a well defined linear functional $\rho_{m+1}$ on $A_m + \mathcal{I}_{m+1} = A_{m+1}$.

We will now prove that (ii–iv) hold for $n = m+1$. By definition $\rho_{m+1}|_{A_m} = \rho_m$, taking care of (ii). In order to check (iii), that is

$$\rho_{m+1}(s_x a s_x^*) = N(x)^{-\beta} \rho_m(a s_x^* s_x), \quad a \in A_m,\ x \in \mathcal{G}, \tag{†}$$

let us first suppose that $m = 0$. Then $a \in A_m = \widetilde{\mathcal{Q}}$ so that $s_x a s_x^* \in \mathcal{I}_{m+1}$ and

$$\rho_{m+1}(s_x a s_x^*) = \chi_{m+1}(s_x a s_x^*) = \sum_{y \in \mathcal{G}} N(y)^{-\beta} \rho_m(s_y^* s_x a s_x^* s_y) =$$

$$= N(x)^{-\beta} \rho_m(s_x^* s_x a s_x^* s_x) = N(x)^{-\beta} \rho_m(a s_x^* s_x).$$

Let us suppose now that $m \geq 1$ and, given that $a \in A_m = A_{m-1} + \mathcal{I}_m$, it is enough to verify (†) separately for $a \in A_{m-1}$ and for $a \in \mathcal{I}_m$.

If $a \in \mathcal{I}_m$ then $s_x a s_x^* \in \mathcal{I}_{m+1}$ and the exact same calculation used to deal with the case $m = 0$ just above gives the conclusion.

If $a \in A_{m-1}$ then $s_x a s_x^* \in A_m$ and, by induction,

$$\rho_{m+1}(s_x a s_x^*) = \rho_m(s_x a s_x^*) = N(x)^{-\beta} \rho_{m-1}(a s_x^* s_x) = N(x)^{-\beta} \rho_m(a s_x^* s_x).$$

This concludes the proof of (†). To prove (iv), that is

$$N(x)^{-\beta} \rho_m\big(s_x^* a s_x\big) = \rho_{m+1}(a s_x s_x^*), \quad a \in A_{m+1},\ x \in \mathcal{G}, \tag{‡}$$

let us again first suppose that $m = 0$. Given that $s_x s_x^* \in \mathcal{I}_1$, which is an ideal in $A_1$ by 11.7.ii, we have that $a s_x s_x^* \in \mathcal{I}_1$ and then

$$\rho_1(a s_x s_x^*) = \chi_{m+1}(a s_x s_x^*) = \sum_{y \in \mathcal{G}} N(y)^{-\beta} \rho(s_y^* a s_x s_x^* s_y) = N(x)^{-\beta} \rho(s_x^* a s_x s_x^* s_x) = N(x)^{-\beta} \rho(s_x^* a s_x),$$

proving (‡) for $m = 0$. Assume now that $m \geq 1$. Given $a \in A_{m+1}$ we have by 11.7.v that $s_x^* a s_x \in A_m$. Therefore, plugging $a := s_x^* a s_x$ into (†) gives

$$\rho_{m+1}(s_x s_x^* a s_x s_x^*) = N(x)^{-\beta} \rho_m(s_x^* a s_x s_x^* s_x)$$

which implies that

$$\rho_{m+1}(a s_x s_x^*) = N(x)^{-\beta} \rho_m(s_x^* a s_x),$$

concluding the proof of (‡) in the general case.

It remains to prove that $\rho_{m+1}$ is a state and we shall derive this from 12.5 applied to the pair $(\rho_m, \chi_{m+1})$. Clearly 12.5.ii holds by 11.8. With respect to checking 12.5.i let us use the approximate unit for $\mathcal{I}_{m+1}$ provided by 11.5. We then have for any $a \geq 0$ in $A_m$ that

$$\widetilde{\chi}_{m+1}(a) = \chi_{m+1}\left(\lim_J \sum_{\mu \in J} p_\mu a\right) = \sum_{\mu \in \mathbb{F}_+^{m+1}} \chi_{m+1}(p_\mu a) =$$

$$= \sum_{\mu\in\mathbb{F}_+^{m+1}} \sum_{x\in\mathcal{G}} N(x)^{-\beta}\rho_m(s_x^* p_\mu a s_x) \leq \sum_{x\in\mathcal{G}} N(x)^{-\beta}\rho_m(s_x^* a s_x),$$

where the last inequality follows from the fact that the $p_\mu$ are pairwise orthogonal (see 6.3.i).

Our present goal is to prove that $\widetilde{\chi}_{m+1}(a) \leq \rho_m(a)$ for all $a \in A_m$. In order to accomplish this let us first suppose that $m = 0$ and that $a = q(X, Y)$, where $X$ and $Y$ are finite subsets of $\mathcal{G}$. By $\mathrm{CK}_3$ we have for any $z \in \mathcal{G}$ that

$$s_z^* a s_z = s_z^* q(X, Y) s_z = A(X, Y, z) q_z$$

(see also the beginning of the proof of 12.1). Therefore,

$$\widetilde{\chi}_{m+1}(a) \leq \sum_{z\in\mathcal{G}} N(z)^{-\beta}\rho_m(s_z^* a s_z) = \sum_{z\in\mathcal{G}} N(z)^{-\beta} A(X, Y, z)\rho(q_z) \leq \rho(q(X, Y)) = \rho(a),$$

where the last inequality is from our hypothesis that $\rho$ is $\beta$-subinvariant. Of course it also follows that $\widetilde{\chi}_{m+1}(a) \leq \rho(a)$ whenever $a$ is a linear combination of the $q(X, Y)$ with positive coefficients. Thus by 11.3 the same holds for any $a \geq 0$ in $\widetilde{\mathcal{Q}} = A_0$.

Assume now that $m \geq 1$. By (iv) applied for $n = m + 1$ (i.e. by (‡)) we have that

$$\widetilde{\chi}_{m+1}(a) \leq \sum_{x\in\mathcal{G}} N(x)^{-\beta}\rho_m(s_x^* a s_x) = \sum_{x\in\mathcal{G}} \rho_m(a s_x s_x^*) \leq \rho_m(a).$$

This concludes the construction of the $\rho_n$ so let us now take up the task of constructing the $\beta$-scaling state $\phi$ mentioned in the statement. Since $C(\widetilde{\Omega}_{\mathcal{T}_A})$ is the closure of $\cup_n A_n$ by 11.7.iv, and $\rho_{n+1}|_{A_n} = \rho_n$, there is a state $\phi$ on $C(\widetilde{\Omega}_{\mathcal{T}_A})$ simultaneously extending all of the $\rho_n$. If $a \in C(\widetilde{\Omega}_{\mathcal{T}_A})$ note that

$$\phi(s_x a s_x^*) = N(x)^{-\beta}\phi(a s_x^* s_x)$$

for all $x \in \mathcal{G}$ by (iii). If moreover $a \in C_0(\Delta^\tau_{x^{-1}})$ then we have

$$\phi(\theta_x(a)) \overset{(6.1.iv)}{=} \phi(s_x a s_x^*) = N(x)^{-\beta}\phi(a s_x^* s_x) \overset{(6.1.iii)}{=} N(x)^{-\beta}\phi(a 1_{\Delta^\tau_{x^{-1}}}) = N(x)^{-\beta}\phi(a),$$

and hence $\phi$ is $\beta$-scaling. Obviously $\phi$ coincides with $\rho$ on $\widetilde{\mathcal{Q}}$.

All of this is meant to work for $\beta < \infty$ but, with the usual interpretation of $N(x)^{-\beta}$, the argument above works also for $\beta = \infty$. Alternatively there is a more straightforward way to prove our statement for $\beta = \infty$ which we would now like to present.

Recall from 11.1.iii that $C(\widetilde{\Omega}_e) \simeq \widetilde{\mathcal{Q}}$ under $\hat{r}$. Identifying these algebras we have that $\rho$ defines a state on $C(\widetilde{\Omega}_e)$ and hence a probability measure on $\widetilde{\Omega}_e$. Extend this to a measure on $\widetilde{\Omega}_{\mathcal{T}_A}$ by declaring that $\widetilde{\Omega}_{\mathcal{T}_A} \setminus \widetilde{\Omega}_e$ has measure zero. This in turn gives the state $\phi$ on $C(\widetilde{\Omega}_{\mathcal{T}_A})$ we are looking for. Precisely, $\phi$ is defined as follows: consider the inclusion $\iota : \widetilde{\Omega}_e \to \widetilde{\Omega}_{\mathcal{T}_A}$ and let $\hat{\iota} : C(\widetilde{\Omega}_{\mathcal{T}_A}) \to C(\widetilde{\Omega}_e)$ be the transposed map. $\phi$ is then the result of the composition

$$C(\widetilde{\Omega}_{\mathcal{T}_A}) \overset{\hat{\iota}}{\longrightarrow} C(\widetilde{\Omega}_e) \overset{\hat{r}}{\longrightarrow} \widetilde{\mathcal{Q}} \overset{\rho}{\longrightarrow} \mathbb{C}.$$

For every $a \in \widetilde{\mathcal{Q}}$ and $\xi \in \widetilde{\Omega}_{\mathcal{T}_A}$ observe that

$$\hat{r}(\hat{\iota}(a))\Big|_\xi = a(r(\xi)) \overset{(11.1.iv)}{=} a(\xi),$$

so that $\hat{r}(\hat{\iota}(a)) = a$ and hence $\phi(a) = \rho(\hat{r}(\hat{\iota}(a))) = \rho(a)$ proving that $\phi$ extends $\rho$.

Observe that for all $x \in \mathcal{G}$ one has that $\hat{\iota}(1_{\Delta^\tau_x})$ is the characteristic function of $\Delta^\tau_x \cap \widetilde{\Omega}_e$, which is the empty set because every $\xi$ in $\Delta^\tau_x$ contains $x$ while the stem of every $\xi \in \widetilde{\Omega}_e$ is trivial. So $\phi(1_{\Delta^\tau_x}) = 0$ for all $x \in \mathcal{G}$. It follows that $\phi(C_0(\Delta^\tau_x)) = \{0\}$ and hence that $\phi$ is an $\infty$-scaling state. □

Putting together 12.1, 12.4, and 12.6 we arrive at one of our main results:

**12.7. Theorem.** *Under 8.1 let $\beta \in (0, \infty]$. Then the correspondence $\phi \mapsto \phi|_{\widetilde{\mathcal{Q}}}$ defines a bijection from the set of $\beta$-scaling states on $C(\widetilde{\Omega}_{\mathcal{T}_A})$ to the set of $\beta$-subinvariant states on $\widetilde{\mathcal{Q}}$.*

## 13. States on $\mathcal{Q}$.

Recall that $\widetilde{\mathcal{Q}}$ is defined to be the *unital* $C^*$-subalgebra of $\widetilde{\mathcal{T}}_A$ generated by the $q_x$. When the emphasis is on $\Omega_{\mathcal{T}_A}$, rather than on $\widetilde{\Omega}_{\mathcal{T}_A}$, it is convenient to work with the algebra $\mathcal{Q} \subseteq C_0(\Omega_{\mathcal{T}_A})$ defined to be the (not necessarily unital) $C^*$-algebra generated by $\{q_x : x \in \mathcal{G}\}$. In the last section we studied states on $C(\widetilde{\Omega}_{\mathcal{T}_A})$ in relation to their restriction to $\widetilde{\mathcal{Q}}$. In order to extend these results to $C_0(\Omega_{\mathcal{T}_A})$ and $\mathcal{Q}$ it would be convenient to know whether states on $C_0(\Omega_{\mathcal{T}_A})$ restrict to states on $\mathcal{Q}$, a fact which is no longer automatic as we are now working with non-necessarily unital $C^*$-algebras. Recall that our matrix $A$ is assumed not to have identically zero *rows*. We will now need to assume 10.1.(COL), i.e. that there are no identically zero *columns*.

**13.1. Proposition.** *Suppose that no column of $A$ is identically zero. Then $\mathcal{Q}$ is an essential subalgebra of $C_0(\Omega_{\mathcal{T}_A})$ in the sense that an approximate identity for $\mathcal{Q}$ is always an approximate identity for $C_0(\Omega_{\mathcal{T}_A})$. Therefore the restriction to $\mathcal{Q}$ of any state on $C_0(\Omega_{\mathcal{T}_A})$ is a state on $\mathcal{Q}$.*

*Proof.* It is clearly enough to show that there is no $\xi \in \Omega_{\mathcal{T}_A}$ such that $a(\xi) = 0$ for all $a \in \mathcal{Q}$. Suppose by contradiction that such a $\xi$ exists. Given $x \in \mathcal{G}$ let $a = q_x$ so that

$$0 = q_x(\xi) = [x^{-1} \in \xi],$$

which implies that $R_\xi(e) = \emptyset$. Suppose first that the stem of $\xi$ is not trivial. In this case there exists $y \in \mathcal{G}$ such that $y \in \xi$. Given that no column of $A$ is zero pick an $x \in \mathcal{G}$ such that $A(x, y) = 1$. Then by Definition 5.1 we have that $x^{-1} \in \xi$ which is a contradiction. The only alternative is then that the stem of $\xi$ is trivial. By [**EL**: 5.12] we conclude that $\xi = \epsilon$ which is again a contradiction since $\epsilon$ was explicitly removed from $\Omega_{\mathcal{T}_A}$. □

We will therefore assume, throughout this section, that no column of $A$ is identically zero keeping, of course, all the other hypothesis in 8.1.

Given a state $\rho$ on $\mathcal{Q}$ it is well known that there exists a unique extension of $\rho$ to a state $\widetilde{\rho}$ on $\widetilde{\mathcal{Q}}$.

**13.2. Definition.** We will say that a state $\rho$ on $\mathcal{Q}$ is $\beta$-invariant (resp. $\beta$-subinvariant) if its canonical extension $\widetilde{\rho}$ is a $\beta$-invariant (resp. $\beta$-subinvariant) state on $\widetilde{\mathcal{Q}}$.

The next result is a generalization of 12.7 to the present context:

**13.3. Theorem.** *Assuming 8.1 and 10.1.(COL) let $\beta \in (0, \infty]$. Then the correspondence $\phi \mapsto \phi|_{\mathcal{Q}}$ defines a bijection from the set of $\beta$-scaling states on $C_0(\Omega_{\mathcal{T}_A})$ to the set of $\beta$-subinvariant states on $\mathcal{Q}$.*

*Proof.* The result follows from 12.7 and 13.1 on noting that states on $C_0(\Omega_{\mathcal{T}_A})$ correspond to states $C(\widetilde{\Omega}_{\mathcal{T}_A})$ whose restriction to $C_0(\Omega_{\mathcal{T}_A})$ is a state (i.e. of norm one), and similarly with respect to $\mathcal{Q}$ and $\widetilde{\mathcal{Q}}$. □

We are now able to give two new characterizations of infinite type states on $C_0(\Omega_{\mathcal{T}_A})$, extending the result obtained in 8.10. We repeat here the conditions of 8.10:

**13.4. Proposition.** *Let $\beta \in (0, \infty)$ and let $\phi$ be a $\beta$-scaling state on $C_0(\Omega_{\mathcal{T}_A})$ corresponding to a measure $\lambda$ on $\Omega_{\mathcal{T}_A}$. Denote by $\rho$ the restriction of $\phi$ to $\mathcal{Q}$ which is a state by 13.1. Then the following are equivalent:*

(i) $\lambda(\Omega_e) = 0$,

(ii) *$\lambda$ is of infinite type,*

(iii) *$\rho$ is a $\beta$-invariant state,*

(iv) $\sum_{x \in \mathcal{G}} N(x)^{-\beta} \rho(q_x) = 1$ *(see 12.2).*

*Proof.* That (i)⇒(ii) was proved in 8.10.

(ii)⇒(iii): Let $X$ and $Y$ be finite subsets of $\mathcal{G}$ and consider the sets

$$S = \{\xi \in \Omega_{\mathcal{T}_A} : x^{-1} \in \xi,\ y^{-1} \notin \xi,\ \ \forall x \in X,\ \forall y \in Y\},$$

and

$$T = \{\xi \in \Omega_{\mathcal{T}_A} : \exists z \in \mathcal{G},\ z \in \xi\ \wedge\ A(X, Y, z) = 1\}.$$

$S$ and $T$ are not necessarily equal but they have exactly the same unbounded elements, as a moment's reflexion based on definition 5.1 will reveal. Given that $\lambda$ is of infinite type, and hence supported in the set of unbounded elements, we must therefore have that $\lambda(S) = \lambda(T)$.

Observe that the characteristic function of $S$ is precisely $q(X, Y)$ while the characteristic function of $T$ is the infinite sum

$$\sum_{z \in \mathcal{G}} A(X, Y, z) p_z.$$

It follows from countable additivity that

$$\phi(q(X, Y)) = \sum_{z \in \mathcal{G}} A(X, Y, z) \phi(p_z).$$

By 8.5 we have

$$\rho(q(X, Y)) = \sum_{z \in \mathcal{G}} A(X, Y, z) N(z)^{-\beta} \rho(q_z),$$

which means that $\rho$ is $\beta$-invariant.

(iii)⇒(iv): Take $X = Y = \emptyset$ above.

(iv)⇒(i): By 8.7.iii we have that

$$\lambda(\Omega^e) = 1 - \sum_{x \in \mathcal{G}} \lambda(\Delta_x^\tau) = 1 - \sum_{x \in \mathcal{G}} \phi(p_x) = 1 - \sum_{x \in \mathcal{G}} N(x)^{-\beta} \phi(q_x) = 0. \qquad \square$$

Given a $\beta$-invariant state $\rho$ on $\mathcal{Q}$ observe that the identity $\sum_{z \in \mathcal{G}} A(X, Y, z) N(z)^{-\beta} \rho(q_z) = \rho(q(X, Y))$ with $X = \{x\}$ and $Y = \emptyset$ becomes

$$\sum_{y \in \mathcal{G}} A(x, y) N(y)^{-\beta} \rho(q_y) = \rho(q_x)$$

for each $x \in \mathcal{G}$. The vector $\big(\rho(q_x)\big)_{x \in \mathcal{G}}$ is therefore a fixed point (i.e. an eigenvector with eigenvalue 1) for the matrix

$$AN^{-\beta} = \{A(x, y) N(y)^{-\beta}\}_{x, y \in \mathcal{G}}.$$

We may take the above expression for $AN^{-\beta}$ as no more than just a definition but note that if $N$ is the diagonal matrix with $N(x, x) := N(x)$ and we interpret $N^{-\beta}$ in the only reasonable way then $AN^{-\beta}$ can be also thought of as the product of $A$ and $N^{-\beta}$.

**13.5. Theorem.** *Under 8.1 let $\beta \in (0, \infty)$. Then the correspondence $\rho \mapsto \big(\rho(q_x)\big)_{x \in \mathcal{G}}$ defines a bijection from the set of $\beta$-invariant states on $\mathcal{Q}$ onto the set of fixed points $v = \big(v_x\big)_{x \in \mathcal{G}}$ for $AN^{-\beta}$ with $v_x \geq 0$ for all $x$ and such that $\sum_{x \in \mathcal{G}} N(x)^{-\beta} v_x = 1$.*

*Proof.* Let $\rho_1$ and $\rho_2$ be $\beta$-invariant states such that $\rho_1(q_x) = \rho_2(q_x)$ for all $x \in \mathcal{G}$. By definition of $\beta$-invariant states it follows that $\rho_1(q(X, Y)) = \rho_2(q(X, Y))$ for all finite sets $X$ and $Y$ of $\mathcal{G}$. So $\rho_1 = \rho_2$, showing our correspondence to be injective.

In order to show that it is also surjective let $v$ be a fixed point as in the statement. Viewing $\mathcal{Q}$ as an ideal in $\widetilde{\mathcal{Q}}$ it is easy to see that the spectrum of $\mathcal{Q}$ is given by

$$\Sigma_{\mathcal{Q}} = \Sigma_A \setminus \{\vec{0}\},$$

where $\vec{0}$ is the zero vector in $2^{\mathcal{G}}$. From our assumption that no column of $A$ is identically zero it then follows that every column of $A$ lies in $\Sigma_{\mathcal{Q}}$.

Define a probability measure on $\Sigma_{\mathcal{Q}}$ by

$$\lambda(S) = \sum_{z\in\mathcal{G}} \big[c_z \in S\big] N(z)^{-\beta} v_z$$

for all Borel subsets $S \subseteq \Sigma_{\mathcal{Q}}$, where $c_z$ refers to the $z^{th}$ column of $A$. By definition $\lambda$ is an atomic measure with atoms the columns of $A$. Each column $c_x$ of $A$ therefore has mass equal to $\sum N(z)^{-\beta} v_z$, where the sum is over the set of $z$'s such that $c_z = c_x$.

Let $\rho$ be the state on $\mathcal{Q}$ given by integration against $\lambda$. Thinking of each $q_x$ as a function on $\widetilde{\Omega}_{\mathcal{T}_A}$, as in 11.2, and observing that $q_x(c_z) = \big[x \in c_z\big] = A(x,z)$, we have

$$\rho(q_x) = \sum_{z\in\mathcal{G}} q_x(c_z) N(z)^{-\beta} v_z = \sum_{z\in\mathcal{G}} A(x,z) N(z)^{-\beta} v_z = v_x.$$

We now wish to show that $\rho$ is $\beta$-invariant. Let therefore $X$ and $Y$ be finite subsets of $\mathcal{G}$ and observe that $q(X,Y)\big|_{c_z} = A(X,Y,z)$. Therefore

$$\rho(q(X,Y)) = \sum_{z\in\mathcal{G}} A(X,Y,z) N(z)^{-\beta} v_z = \sum_{z\in\mathcal{G}} A(X,Y,z) N(z)^{-\beta} \rho(q_z),$$

proving that $\rho$ satisfies the required properties. □

The following summarizes much of what we have discovered so far:

**13.6. Theorem.** *Assuming 8.1 and 10.1.*(COL) *let* $\beta \in (0,\infty]$. *Then the vertical correspondences below are bijective*

| | | |
|---|---|---|
| *$\sigma$–KMS states on $\mathcal{T}_A$ at inverse temperature $\beta$* | | |
| *restriction* $\downarrow$ *to* $C_0(\Omega_{\mathcal{T}_A})$ | | |
| *$\beta$-scaling states on $C_0(\Omega_{\mathcal{T}_A})$* | $\supseteq$ | *infinite type $\beta$-scaling states* |
| *restriction* $\downarrow$ *to* $\mathcal{Q}$ | | *restriction* $\downarrow$ *to* $\mathcal{Q}$ |
| *$\beta$-subinvariant states on $\mathcal{Q}$* | $\supseteq$ | *$\beta$-invariant states on $\mathcal{Q}$* |
| | | *maps $\rho$ to* $\downarrow$ $v = \big(\rho(q_x)\big)_{x\in\mathcal{G}}$ |
| | | *nonnegative fixed points $v$ for $AN^{-\beta}$ with* $\sum_{x\in\mathcal{G}} N(x)^{-\beta} v_x = 1.$ |

*Proof.* We restrict ourselves to pointing out the result relating to each one of the above arrows. The uppermost arrow corresponds to 8.2. The arrow following that is 13.3. On the second column the uppermost arrow is 13.4 and the last one is 13.5. □

## 14. The fixed-source-and-target partition function $Z_{xy}(\beta)$.

In this section we will introduce the third family of Dirichlet series associated to our context.

**14.1. Definition.** Let $x, y \in \mathcal{G}$. The *fixed-source-and-target partition function* relative to the pair of generators $x$ and $y$ for the dynamical system $(\mathcal{T}_A, \sigma, \mathbb{R})$ is the function $Z_{xy}(\beta)$ given by the Dirichlet series

$$Z_{xy}(\beta) = \sum_{\mu \in P_A^{xy}} N(\mu)^{-\beta}, \quad \beta \in (0, \infty),$$

where $P_A^{xy}$ is the set of all admissible words beginning in $x$ and ending in $y$.

As in 10.2 we have:

**14.2. Proposition.** *Let $x_1, x_2, y_1, y_2 \in \mathcal{G}$. Suppose that there are admissible words $\nu, \gamma \in P_A$ such that*

- $\nu_1 = x_1$, $\nu_{|\nu|} = x_2$,
- $\gamma_1 = y_2$, $\gamma_{|\gamma|} = y_1$,

*then for every $\beta \in (0, \infty)$ one has that $Z_{x_2y_2}(\beta) \le K Z_{x_1y_1}(\beta)$, where $K = N(\nu\gamma)^{\beta} N(x_2y_2)^{-\beta}$.*

*Proof.* Considering the injective map $\mu \in P_A^{x_2y_2} \mapsto \nu x_2^{-1} \mu y_2^{-1} \gamma \in P_A^{x_1y_1}$ we have

$$Z_{x_1y_1}(\beta) = \sum_{\mu \in P_A^{x_1y_1}} N(\mu)^{-\beta} \ge \sum_{\mu \in P_A^{x_2y_2}} N(\nu x_2^{-1} \mu y_2^{-1} \gamma)^{-\beta} =$$

$$= N(\nu\gamma)^{-\beta} N(x_2y_2)^{\beta} \sum_{\mu \in P_A^{x_2y_2}} N(\mu)^{-\beta} = N(\nu\gamma)^{-\beta} N(x_2y_2)^{\beta} Z_{x_2y_2}(\beta). \qquad \square$$

Assuming that $A$ is irreducible we have another "solidarity" property (see 10.3) among these Dirichlet series:

**14.3. Proposition.** *Let $A$ be irreducible. Then for every $\beta \in (0, \infty)$ one has that either*

- $Z_{xy}(\beta) < \infty$ *for all* $x, y \in \mathcal{G}$, *or*
- $Z_{xy}(\beta) = \infty$ *for all* $x, y \in \mathcal{G}$.

*Proof.* Follows immediately from 14.2. $\square$

**14.4. Definition.** Under the hypothesis that $A$ is irreducible the abscissa of convergence for each and every one of the Dirichlet series $Z_{xy}(\beta)$ will be called the *fixed-source-and-target critical inverse temperature* and will be denoted $\ddot{\beta}_c$. The set of $\beta$'s where each and every one of these series converge, including $\beta = \infty$, will be called the *interval of fixed-source-and-target super-critical inverse temperatures* and will be denoted $\ddot{I}_c$.

As before, it is obvious that

$$\ddot{\beta}_c \le \dot{\beta}_c \le \beta_c, \quad \text{and} \quad \ddot{I}_c \supseteq \dot{I}_c \supseteq I_c.$$

The relevance of these concepts lies in the following:

**14.5. Theorem.** *Suppose 8.1 and 10.1.(IRR) and let $\beta < \ddot{\beta}_c$. Then there are no $\beta$-scaling states at all on $C_0(\Omega_{\mathcal{T}_A})$ and hence neither are there KMS$_\beta$ states on $\mathcal{T}_A$.*

*Proof.* Assume by contradiction that $\phi$ is a $\beta_0$-scaling state on $C_0(\Omega_{\mathcal{T}_A})$ for some $\beta_0 < \ddot{\beta}_c$. Then by 13.3 the restriction $\rho$ of $\phi$ to $\mathcal{Q}$ is a $\beta_0$-subinvariant state. Let $x \in \mathcal{G}$ and plug $X = \{x\}$ and $Y = \emptyset$ in the definition of subinvariant states (see 12.3) to get

$$\sum_{y\in\mathcal{G}} A(x,y)N(y)^{-\beta_0}\rho(q_y) \le \rho(q_x).$$

This says that the nonnegative vector $v = \big(\rho(q_x)\big)_{x\in\mathcal{G}}$ is a *right 1-subinvariant* vector for the irreducible matrix $AN^{-\beta_0}$ in the sense of [**V**: Section 4]. We will now proceed to show that such a vector cannot exist, therefore arriving at a contradiction.

Unfortunately we cannot just quote the result we need from [**V**: Corollary 1] because of the incompatibility between our point of view which emphasizes Dirichlet series in the variable $\beta$, and Vere-Jones's point of view which emphasizes power series. Nevertheless, proceeding with the necessary care, we may still derive our conclusions from [**V**].

Set the matrix $T$ of [**V**] to be $AN^{-\beta}$ and, according to [**V**: Section 2], let $f_{ij}^{(n)}$ be the "first-entrance probabilities" for $T$, and $F_{ij}(z)$ be the corresponding generating function.

Observe that $F_{ij}(z)$ depends on $\beta$, as $T$ definitely does. Accordingly let us denote by $\mathcal{F}_{ij}(\beta)$ the value of $F_{ij}$ at $z = 1$.

By the right-hand-sided version of [**V**: Lemma 4.1], applied for $\beta = \beta_0$, $r = 1$, and $i = j$ taken to be any fixed element in $\mathcal{G}$, we have that $F_{ii}(1) = \mathcal{F}_{ii}(\beta_0) \le 1$. It is easy to see that $\mathcal{F}_{ii}(\beta)$ is a strictly decreasing function of $\beta$ and hence $\mathcal{F}_{ii}(\beta) < 1$ for all $\beta > \beta_0$.

Observe that the generating function $T_{ij}$, defined near the bottom of page 362 of [**V**], is related to our partition function $Z_{ij}$ by

$$T_{ij}(1) = N(i)^{\beta} Z_{ij}(\beta).$$

Using equation (3) in [**V**], namely $T_{ii}(z) = 1/(1 - F_{ii}(z))$, for $z = 1$ we therefore have that

$$N(i)^{\beta} Z_{ii}(\beta) = \frac{1}{1 - \mathcal{F}_{ii}(\beta)}.$$

Inspired by the idea of the proof of Lemma 2.1 in [**V**], we conclude that $Z_{ii}(\beta)$ has no singularities in the interval $(\beta, \infty)$ because, as seen above, $\mathcal{F}_{ii} < 1$ there. Since $\beta < \ddot{\beta}_c$ we have a contradiction. □

We may now throw some more conclusions into Diagram 10.7 getting the following information about $\beta$-scaling states on $C_0(\Omega_{\mathcal{T}_A})$, and hence also about KMS$_\beta$ states on $\mathcal{T}_A$, again in the case that $A$ is irreducible.

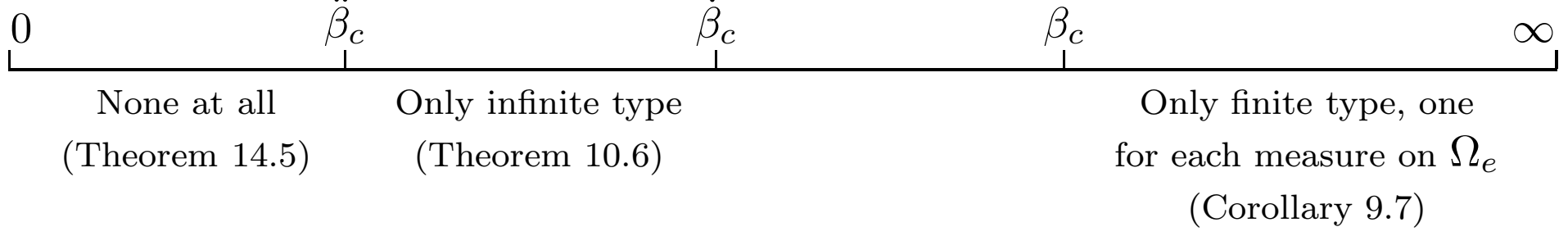

*Diagram 14.6*

## 15. Energy bounded below.

In this section we will prove that there are no $\beta$-scaling states for $\beta < \dot\beta_c$ under the hypothesis that the "energy" parameters $N(x)$ satisfy $\inf_{x\in\mathcal{G}} N(x) > 1$ (see 10.1.(INF)). The main tool to be used is the following Lemma which takes advantage of the fact that a single state $\rho$ on $\mathcal{Q}$ may be used to determine scaling states for different values of $\beta$, as long as $\rho$ remains $\beta$-subinvariant. In particular note that if $\rho$ is $\beta$-subinvariant for some $\beta$ then the same holds for any $\beta' > \beta$, i.e. when the "temperature" $1/\beta$ decreases.

**15.1. Lemma.** (Cooling Lemma) *Assume 10.1.*(COL+INF) *and let $\beta \in (0,\infty)$. Given a $\beta$-scaling state $\phi$ on $C_0(\Omega_{\mathcal{T}_A})$ set $\rho = \phi|_{\mathcal{Q}}$, so that $\rho$ is a $\beta$-subinvariant state on $\mathcal{Q}$. Let $\beta' > \beta$ and observe that $\rho$ is clearly also $\beta'$-subinvariant. Let $\phi'$ be the unique $\beta'$-scaling state on $C_0(\Omega_{\mathcal{T}_A})$ whose restriction to $\mathcal{Q}$ coincides with $\rho$, by 13.3. Then $\phi'$ is of finite type.*

*Proof.* Let $\lambda$ and $\lambda'$ be the measures on $\Omega_{\mathcal{T}_A}$ corresponding to $\rho$ and $\rho'$, respectively. Recall from 8.12 that

$$\lambda'(\Omega_\infty) = \lim_{n\to\infty} \sum_{\mu\in P_A^n} N(\mu)^{-\beta'} \phi'(q_\mu).$$

However $\phi'(q_\mu) = \rho(q_\mu) = \phi(q_\mu)$ so that the above expression for $\lambda'(\Omega_\infty)$ actually depends only on $\beta'$. Let $R = \inf_{x\in\mathcal{G}} N(x)$ and let $\delta = \beta' - \beta$. Observe that for all $x \in \mathcal{G}$ one has

$$N(x)^{-\beta'} = N(x)^{-\delta} N(x)^{-\beta} \leq R^{-\delta} N(x)^{-\beta}.$$

It follows that for all $\mu \in P_A^n$ we have $N(\mu)^{-\beta'} \leq R^{-n\delta} N(\mu)^{-\beta}$ and hence that

$$\lambda'(\Omega_\infty) \leq \lim_{n\to\infty} R^{-n\delta} \sum_{\mu\in P_A^n} N(\mu)^{-\beta}\phi(q_\mu).$$

Observe that for all $n$ one has, using 8.5, that

$$\sum_{\mu\in P_A^n} N(\mu)^{-\beta}\phi(q_\mu) = \sum_{\mu\in P_A^n} \phi(p_\mu) \leq 1$$

because the $p_\mu$ are pairwise orthogonal projections. Since $R^{-n\delta} \to 0$ as $n \to \infty$, we conclude that $\lambda'(\Omega_\infty) = 0$ and hence that $\lambda'$ is of finite type. ☐

As a conclusion we may boost the result obtained in 10.6:

**15.2. Theorem.** *Assume 10.1.*(IRR+INF) *and let $\beta < \dot\beta_c$. Then there are no $\beta$-scaling states at all on $C_0(\Omega_{\mathcal{T}_A})$ and hence neither are there KMS$_\beta$ states on $\mathcal{T}_A$.*

*Proof.* Suppose by contradiction that $\phi$ is a $\beta$-scaling state on $C_0(\Omega_{\mathcal{T}_A})$. Choose $\delta > 0$ such that $\beta' := \beta + \delta < \dot\beta_c$ and, using 13.3, let $\phi'$ be the unique $\beta'$-scaling state such that $\phi'|_{\mathcal{Q}} = \phi|_{\mathcal{Q}}$. Then $\phi'$ is of finite type by 15.1 contradicting 10.6. ☐

The following diagram gives information about $\beta$-scaling states on $C_0(\Omega_{\mathcal{T}_A})$, and hence also about KMS$_\beta$ states on $\mathcal{T}_A$, under the hypothesis of 15.2 improving upon Diagram 10.7:

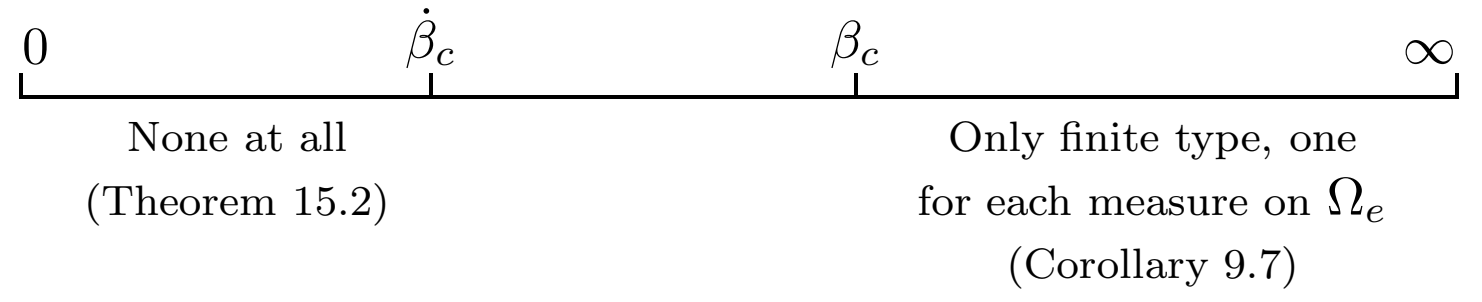

*Diagram 15.3*

Unfortunately we don't have much more to say about the case in which $\beta$ lies in the interval between $\dot\beta_c$ and $\beta_c$. Nevertheless this mysterious interval some times collapses, as in the following situation:

**15.4. Proposition.** *Assume 10.1.*(IRR+FTS)*. Then* $\dot{\beta}_c = \beta_c$ *and* $\dot{I}_c = I_c$.

*Proof.* Let $\{y_1, \ldots, y_n\} \subseteq \mathcal{G}$ be a finite target set as in 10.1.(FTS). Decompose $\mathcal{G}$ in a disjoint union $\mathcal{G} = \bigcup_{i=1}^n \mathcal{G}_i$ such that for every $x \in \mathcal{G}_i$ one has $A(x, y_i) = 1$. Consequently $P_A$ decomposes as the disjoint union $P_A = \{e\} \cup \bigcup_{i=1}^n P_A^{y_i}$ and so for every $\beta$

$$Z(\beta) = \sum_{\mu \in P_A} N(\mu)^{-\beta} = 1 + \sum_{i=1}^n \sum_{\mu \in P_A^{y_i}} N(\mu)^{-\beta} = 1 + \sum_{i=1}^n Z_{y_i}(\beta).$$

Therefore if $\beta \in \dot{I}_c$ we have that $Z_{y_i}(\beta) < \infty$ for all $i$ and hence $Z(\beta) < \infty$ so that $\beta \in I_c$. □

Under the hypotheses of 15.2 and 15.4, i.e. all of the hypotheses listed in 10.1, Diagram 15.3 therefore gives as much information as we could possibly want about $\beta$-scaling states throughout the whole interval $(0, \infty]$, except perhaps at the critical point.

## 16. An example of behavior at the critical point.

In this section we will show that, even if one assumes all of the hypotheses listed in 10.1, there is not much more that can be said in general about the nature of KMS states at the critical inverse temperature $\beta_c$. We will eventually prove that the following antagonistic situations may occur:

(a) The $\text{KMS}_{\beta_c}$ state may be unique and of infinite type.
(b) There may be infinitely many $\text{KMS}_{\beta_c}$ states all of which are of finite type.

In fact in this section we will just give an example of situation (b) since we will later show that situation (a) is the rule for finite irreducible matrices. Let

$$\zeta(\beta) = \sum_{k=1}^{\infty} N_k^{-\beta}$$

be any Dirichlet series which converges at its abscissa of convergence, say $\bar{\beta}$, with $\bar{\beta} \in (0, \infty)$.

Put $\mathcal{G} = \mathbb{N}$. It is relevant to us that $\mathcal{G}$ consist of one element for each term of the above series, plus one more element, namely zero. Accordingly we will write $\mathcal{G}_* = \mathbb{N} \setminus \{0\}$.

Consider the matrix

$$A = \begin{bmatrix} 0 & 1 & 1 & 1 & \ldots \\ 1 & 0 & 0 & 0 & \ldots \\ 1 & 0 & 0 & 0 & \ldots \\ 1 & 0 & 0 & 0 & \ldots \\ \vdots & \vdots & \vdots & \vdots & \ddots \end{bmatrix},$$

whose index set is $\mathcal{G} \times \mathcal{G}$. In other words the $0^{th}$ column and the $0^{th}$ row of $A$ consist of ones, except for $A(0, 0)$ which is zero. All other entries are zero.

Clearly $A$ is irreducible and satisfies 10.1.(COL). Observe that $A$ also satisfies 10.1.(FTS) since for every $x \in \mathcal{G}_*$ one has that $A(x, 0) = 1$, while $A(0, 1) = 1$. That is, the set $\{0, 1\}$ is a finite target set.

Discarding a finite number of terms and relabeling we may suppose that

$$\zeta(\bar{\beta}) = \sum_{k=1}^{\infty} N_k^{-\bar{\beta}} < 2^{\bar{\beta}}. \tag{†}$$

The convergence of the above Dirichlet series implies that $\lim_{k \to \infty} N_k = \infty$ and hence, discarding another finite set of terms, we may suppose that $N_k \geq 2$ for all $k$. Set $N(0) = 2$ and $N(k) = N_k$ for all $k \in \mathcal{G}_*$, so that 10.1.(INF) holds. We are therefore under a situation in which everything in 10.1 holds.

We would now like to compute the partition function $Z_0(\beta)$. In order to do this observe that the admissible words ending in 0 are precisely of the form

$$\mu = \begin{cases} x_1\,0\,x_2\,0\,\ldots\,0\,x_n\,0 & \text{if } |\mu| = 2n, \text{ or} \\ 0\,x_1\,0\,x_2\,0\,\ldots\,0\,x_n\,0 & \text{if } |\mu| = 2n+1, \end{cases}$$

where $\vec{x} = (x_1, x_2, \ldots, x_n)$ is an arbitrary element of $\mathcal{G}_*^n$. Therefore

$$Z_0(\beta) = \sum_{n=0}^{\infty} 2^{-n\beta} \sum_{\vec{x}\in\mathcal{G}_*^n} N(x_1)^{-\beta} \cdots N(x_n)^{-\beta} + \sum_{n=0}^{\infty} 2^{-(n+1)\beta} \sum_{\vec{x}\in\mathcal{G}_*^n} N(x_1)^{-\beta} \cdots N(x_n)^{-\beta}.$$

By considering the summands corresponding to $n = 1$ above we see that $Z_0(\beta)$ diverges when $\zeta(\beta)$ diverges. Therefore the convergence interval for $Z_0(\beta)$ is contained in $[\bar{\beta}, \infty)$. Moreover notice that for all $n \in \mathbb{N}$ we have

$$\sum_{\vec{x}\in\mathcal{G}_*^n} N(x_1)^{-\beta} \cdots N(x_n)^{-\beta} = \left(\sum_{x\in\mathcal{G}_*} N(x)^{-\beta}\right)^n = \zeta(\beta)^n.$$

So,

$$Z_0(\beta) = \sum_{n=0}^{\infty} 2^{-n\beta}\zeta(\beta)^n + \sum_{n=0}^{\infty} 2^{-(n+1)\beta}\zeta(\beta)^n = \left(1 + 2^{-\beta}\right) \sum_{n=0}^{\infty} \left(2^{-\beta}\zeta(\beta)\right)^n.$$

We therefore see that $Z_0(\bar{\beta})$ is a converging geometric series by (†). It follows that $\dot{\beta}_c = \bar{\beta}$ and $\dot{I}_c = [\bar{\beta}, \infty]$. By 15.4 we also have $I_c = [\bar{\beta}, \infty]$.

Combining 15.2 with 9.7 we therefore obtain:

**16.1. Proposition.** *Let $A$, $N$, and $\bar{\beta}$ be given as above. Then*

(i) *For $\beta < \bar{\beta}$ there are no KMS$_\beta$ states on $\mathcal{T}_A$,*

(ii) *For $\beta \geq \bar{\beta}$ the simplex of KMS$_\beta$ states on $\mathcal{T}_A$ is affine-homeomorphic to the simplex of finite measures on $\Omega_e$ such that $Z(\beta, \gamma) = 1$.*

In order to best appreciate this result it is important to observe that for all measures $\gamma$ on $\Omega_e$ one has

$$Z(\beta, \gamma) = \gamma(\Omega_e) + \sum_{x\in\mathcal{G}} Z_x(\beta)\gamma(\Omega_e^x) \leq Z(\beta)\gamma(\Omega_e),$$

and hence, as long as $Z(\beta)$ is finite, that is as long as $\beta \geq \bar{\beta}$, for any finite measure $\gamma$ on $\Omega_e$ the measure $\gamma/Z(\beta, \gamma)$ fits into 16.1.ii. If follows that there are infinitely many KMS states at the critical inverse temperature.

## 17. KMS states on $\mathcal{O}_A$.

Recall that $\widetilde{\mathcal{O}}_A$ is the quotient of $\widetilde{\mathcal{T}}_A$ obtained by imposing relation $\mathrm{CK}_4$ in addition to $\mathrm{CK}_{1-3}$. Clearly the quotient map

$$\Pi : \widetilde{\mathcal{T}}_A \to \widetilde{\mathcal{O}}_A$$

is then covariant for the respective one-parameter automorphism groups. For every KMS state $\psi$ on $\widetilde{\mathcal{O}}_A$ one therefore has that $\psi \circ \Pi$ is a KMS state on $\widetilde{\mathcal{T}}_A$ and hence the simplex of KMS states on $\widetilde{\mathcal{O}}_A$ may be seen as a subset of the KMS states on $\widetilde{\mathcal{T}}_A$. This section is dedicated to giving a characterization of this subset.

Nevertheless it should be observed that occasionally it happens that $\widetilde{\mathcal{T}}_A = \widetilde{\mathcal{O}}_A$ (and hence also $\mathcal{T}_A = \mathcal{O}_A$) and we start by characterizing when exactly this is the case.

**17.1. Proposition.** *For any 0–1 matrix $A = \{A(x, y)\}_{x,y\in\mathcal{G}}$ having no identically zero rows the following are equivalent:*

(i) *Given any neighborhood $V$ of any point $c \in \Sigma_A$ there are infinitely many $j \in \mathcal{G}$ such that the column $c_j$ of $A$ lies in $V$,*

(ii) $\widetilde{\mathcal{T}}_A = \widetilde{\mathcal{O}}_A$.

*Proof.* Suppose that (i) holds. Then the *closure* of set of columns of $A$ within $2^{\mathcal{G}}$, namely $\Sigma_A$, coincides with the set of *accumulation points* of the columns of $A$. Therefore $\widetilde{\Omega}_{\mathcal{T}_A} = \widetilde{\Omega}_{\mathcal{O}_A}$ by [**EL**: 7.7] and hence $\widetilde{\mathcal{T}}_A = \widetilde{\mathcal{O}}_A$. The converse is proven by running this argument backwards. □

It should be remarked that condition 17.1.i comes close to saying that $\Sigma_A$ is a perfect topological space (i.e. that it has no isolated points) except that when one "counts" how many columns there are in a neighborhood one should look at the set of *indices* rather than at the set of columns itself. When all columns are distinct (or repeated at most finitely often) one has that 17.1.i is therefore equivalent to $\Sigma_A$ being perfect. Regardless of the columns being distinct, if $\Sigma_A$ is perfect one clearly has that 17.1.i holds.

Under the above circumstances the study of KMS states on $\mathcal{O}_A$ is therefore identical to the corresponding study for $\mathcal{T}_A$. In the opposite case, however, it is useful to obtain criteria to distinguish, among the KMS states on $\mathcal{T}_A$, which ones factor through $\mathcal{O}_A$.

**17.2. Theorem.** *Assuming 8.1 let $\beta \in (0, \infty]$. Also let*

- *$\psi$ be a KMS$_\beta$ state on $\mathcal{T}_A$,*
- *$\phi$ be the restriction of $\psi$ to $C_0(\Omega_{\mathcal{T}_A})$,*
- *$\lambda$ be the measure on $\Omega_{\mathcal{T}_A}$ representing $\phi$, and*
- *$\rho$ be the restriction of $\phi$ to $\mathcal{Q}$.*

*Then the following are equivalent:*

(i) *there exists a KMS$_\beta$ state $\psi'$ on $\mathcal{O}_A$ such that $\psi = \psi' \circ \Pi$,*

(ii) *$\rho(q(X, Y)) = \sum_{z\in\mathcal{G}} A(X, Y, z)N(z)^{-\beta}\rho(q_z)$ whenever $X, Y \subseteq \mathcal{G}$ are finite and $A(X, Y, z)$ is finitely supported as a function of $z$,*

(iii) *the support of $\lambda$ is contained in the closure of $\Omega_\infty$.*

*Proof.* Assume (i). If $X$ and $Y$ are as in (ii) then $\Pi(q(X, Y)) = \sum_{z\in\mathcal{G}} A(X, Y, z)\Pi(p_z)$ by $\mathrm{CK}_4$ and hence

$$\rho(q(X, Y)) = \psi(q(X, Y)) = \psi'\big(\Pi(q(X, Y))\big) = \psi'\left(\sum_{z\in\mathcal{G}} A(X, Y, z)\Pi(p_z)\right) =$$

$$= \sum_{z\in\mathcal{G}} A(X, Y, z)\psi(p_z) = \sum_{z\in\mathcal{G}} A(X, Y, z)N(z)^{-\beta}\rho(q_z),$$

proving (ii).

Assume (ii) and let $\xi \in \Omega_{\mathcal{T}_A} \setminus \overline{\Omega}_\infty$. Pick a neighborhood $V$ of $\xi$ disjoint from $\Omega_\infty$ which, by [**EL**: 6.2], may be chosen so as to have the form

$$V = \left\{ \begin{array}{lll} \eta \in \Omega_{\mathcal{T}_A}: & \omega \in \eta, & \\ & \omega x^{-1} \in \eta, & \text{for } x \text{ in } X, \\ & \omega y^{-1} \notin \eta, & \text{for } y \text{ in } Y, \\ & \omega z \notin \eta, & \text{for } z \text{ in } Z \end{array} \right\},$$

where $\omega$ is the (finite) stem of $\xi$ and $X$, $Y$, and $Z$ are finite subsets of $\mathcal{G}$, with $X \subseteq R_\xi(\omega)$ and $Y \cap R_\xi(\omega) = \emptyset$. We wish to show that $\lambda(V) = 0$ from which it will follow that $\xi$ is not in the support of $\lambda$ thus proving (iii). Let $U = \alpha_{\omega^{-1}}(V)$ so that

$$U = \left\{ \begin{array}{lll} \eta \in \Omega_{\mathcal{T}_A}: & x^{-1} \in \eta, & \text{for } x \text{ in } X, \\ & y^{-1} \notin \eta, & \text{for } y \text{ in } Y, \\ & z \notin \eta, & \text{for } z \text{ in } Z \end{array} \right\},$$

Because $\lambda$ is $\beta$-scaling we have that

$$\lambda(V) = \lambda(\alpha_\omega(U)) = N(\omega)^{-\beta}\lambda(U)$$

(when $\beta = \infty$ and $\omega \neq e$ this should be interpreted as zero), so it is enough to show that $\lambda(U) = 0$. Note that the characteristic function of $U$ is given precisely by

$$1_U = q(X,Y)\prod_{z\in Z}(1-p_z) = q(X,Y)\left(1 - \sum_{z\in Z} p_z\right) = q(X,Y) - \sum_{z\in Z} q(X,Y)p_z =$$

$$= q(X,Y) - \sum_{z\in Z} A(X,Y,z)p_z,$$

where the last equality follows from $\mathrm{CK}_3$ as shown in the beginning of the proof of 12.1. We claim that $A(X,Y,z) = 0$ for all $z \notin Z$. Arguing by contradiction suppose that $z_0 \notin Z$ and $A(X,Y,z_0) = 1$. Therefore $A(x,z_0) = 1$ for all $x \in X$ and $A(y,z_0) = 0$ for all $y \in Y$.

Pick an infinite admissible word $\nu$ beginning in $z_0$ (which exists because no row of $A$ is zero). By [**EL**: 5.13] there exists $\eta \in \Omega_\infty$ whose stem coincides with $\nu$. Inspecting definition 5.1, observing that $\eta \in \widetilde{\Omega}_{\mathcal{TO}_A}$, and noting that $z_0 \in \eta$ it is easy to show that $\eta \in U$. This contradicts the fact that $U$ and $\Omega_\infty$ are disjoint and hence we see that $A(X,Y,z) = 0$ for all $z \notin Z$ as claimed. Using (ii) we therefore have

$$\lambda(U) = \phi(1_U) = \phi(q(X,Y)) - \sum_{z\in Z} A(X,Y,z)\phi(p_z) =$$

$$= \rho(q(X,Y)) - \sum_{z\in Z} A(X,Y,z)N(z)^{-\beta}\rho(q_z) = 0.$$

This proves that (ii) implies (iii). In order to prove that (iii) implies (i) let $\lambda'$ be the restriction of $\lambda$ to a measure on $\overline{\Omega}_\infty = \Omega_{\mathcal{O}_A}$ (see [**EL**: 7.3]) which is a probability measure by hypothesis. Obviously $\lambda'$ is $\beta$-scaling and hence by 8.2 there exists a $\mathrm{KMS}_\beta$ state $\psi'$ on $\mathcal{O}_A$ whose restriction to $C_0(\Omega_{\mathcal{O}_A})$ is given by integration against $\lambda'$. We claim that $\psi = \psi' \circ \Pi$. Given that $\psi = \phi \circ E$ by 8.2 and similarly for $\psi'$ it is enough to verify that $\psi$ and $\psi' \circ \Pi$ coincide on $C_0(\Omega_{\mathcal{T}_A})$ but this is now obvious. ☐

It should be remarked that every infinite type $\beta$-scaling measure $\lambda$ on $\Omega_{\mathcal{T}_A}$ satisfies 17.2.iii and hence is associated to a $\mathrm{KMS}_\beta$ state on $\mathcal{O}_A$. However, since $\Omega_\infty$ is not necessarily closed, there may exist measures supported in $\overline{\Omega}_\infty$ which are not of infinite type.

## 18. The finite dimensional case.

Throughout this section we assume that $\mathcal{G}$ is a finite set and hence $A$ is a finite matrix. Many simplifications take place under this hypothesis and the results can be stated a bit more conclusively.

Being the closure of the set of columns of $A$ within $2^{\mathcal{G}}$, $\Sigma_A$ is hence a finite space with d(A) points, where:

**18.1. Definition.** We denote by $d(A)$ the number of distinct columns of $A$.

Throughout this section we will assume that no column of $A$ is zero. So the zero vector does not belong to $\Sigma_A$ and it follows from the definition of $\widetilde{\Omega}_{\mathcal{T}_A}$ that $\epsilon \notin \widetilde{\Omega}_{\mathcal{T}_A}$ and hence that $\Omega_{\mathcal{T}_A} = \widetilde{\Omega}_{\mathcal{T}_A}$ and $\Omega_{\mathcal{O}_A} = \widetilde{\Omega}_{\mathcal{O}_A}$. This implies that $\mathcal{T}_A = \widetilde{\mathcal{T}}_A$ as well as that $\mathcal{O}_A = \widetilde{\mathcal{O}}_A$. In other words all algebras and spaces with no tilde coincide with their tilde versions in the last two rows of table 7.1 (the same not holding for the first row because $\epsilon \in \widetilde{\Omega}_{\mathcal{T}\mathcal{O}_A}$ always).

Before we proceed we need the following consequence of the Perron–Frobenius Theorem:

**18.2. Lemma.** *Let $\mathcal{S}$ be the set of all $n \times n$ irreducible (in the sense of* [**Se**: Definition 1.6]*) nonnegative matrices. Let $\mathcal{S}_1$ be the subset of $\mathcal{S}$ formed by the matrices $M$ such that $\sum_{n=0}^{\infty} M^n$ converges. Then $\mathcal{S}_1$ is open in $\mathcal{S}$.*

*Proof.* Let $M \in \mathcal{S}_1$. Then clearly $\sum_{n=0}^{\infty} (rM)^n$ converges for all $r \in (0, 1]$. Therefore $1 - rM$ is invertible for all such $r$ and hence no eigenvalue of $M$ lies in the interval $[1, \infty)$. By the Perron–Frobenius Theorem [**Se**: 1.5] it follows that the spectral radius of $M$ is strictly less than 1. Since the spectrum is lower semicontinuous there exists a neighborhood of $M$ consisting solely of matrices whose spectral radius is less than 1. This neighborhood is therefore contained in $\mathcal{S}_1$. □

Let us first study the three critical inverse temperatures for a finite irreducible matrix. As before we will denote by $N$ the diagonal matrix with $N(x, x) := N(x)$.

**18.3. Proposition.** *Under 8.1 let $A$ be a finite irreducible matrix. Then*

(i) $\ddot{\beta}_c = \dot{\beta}_c = \beta_c < \infty$,
(ii) $\ddot{I}_c = \dot{I}_c = I_c = (\beta_c, \infty]$, *and*
(iii) *the spectral radius of $AN^{-\beta_c}$ is 1.*

*Proof.* Breaking the admissible words according to their final and initial letter we have that

$$Z(\beta) \; = \; 1 + \sum_{y \in \mathcal{G}} Z_y(\beta) \; = \; 1 + \sum_{x,y \in \mathcal{G}} Z_{xy}(\beta).$$

Since $\mathcal{G}$ is finite we therefore have that $Z(\beta) < \infty$ if and only if $Z_y(\beta) < \infty$ for all $y \in \mathcal{G}$ if and only if $Z_{xy}(\beta) < \infty$ for all $x, y \in \mathcal{G}$. Therefore $\ddot{I}_c = \dot{I}_c = I_c$ and hence also $\ddot{\beta}_c = \dot{\beta}_c = \beta_c$.

For a large enough $\beta$ one clearly has that $\sum_{y \in \mathcal{G}} N(y)^{-\beta} < 1$ and hence, using 8.16, we have that $Z(\beta) < \infty$. This shows that $\beta_c < \infty$.

For every $x$ and $y$ in $\mathcal{G}$ note that the $(x, y)$ entry of the formal power series of matrices

$$\sum_{n=0}^{\infty} \left(AN^{-\beta}\right)^n \tag{†}$$

is precisely given by $N(x)^{\beta} Z_{xy}(\beta)$. Therefore that series converges if and only if all $Z_{xy}(\beta) < \infty$, which is the same as saying that $\beta \in \ddot{I}_c$. In other words

$$\ddot{I}_c = \left\{\beta \in (0, \infty) : AN^{-\beta} \in \mathcal{S}_1\right\} \cup \{\infty\},$$

where $\mathcal{S}_1$ is as in 18.2. So $\ddot{I}_c$ is an open set (in the extended real line) and it follows that $\ddot{I}_c = (\ddot{\beta}_c, \infty]$ or, equivalently, that $I_c = (\beta_c, \infty]$.

In order to prove (iii) let, for each $\beta \in \mathbb{R}$, $r(\beta)$ be the spectral radius of $AN^{-\beta}$. If $\beta > \beta_c$ we have seen that (†) converges and hence $r(\beta) \leq 1$. Taking the limit as $\beta \to \beta_c$ we conclude that $r(\beta_c) \leq 1$.

Suppose by contradiction that $r(\beta_c) < 1$. Then we would have that (†) converges for $\beta = \beta_c$ which was ruled out in (ii). □

With this we may give a precise description of the KMS states on $\mathcal{T}_A$:

**18.4. Theorem.** *Under 8.1 let $A$ be a finite irreducible 0–1 matrix. Then:*

(i) *For $\beta > \beta_c$ the KMS$_\beta$ states on $\mathcal{T}_A$ form a simplex of dimension $d(A) - 1$ which is affine homeomorphic to the simplex of all measures $\gamma$ on the finite measure space $\Omega_e$ such that $Z(\beta, \gamma) = 1$.*

(ii) *For $\beta = \beta_c$ there exists precisely one KMS$_\beta$ state $\psi$. Its restriction to $C_0(\Omega_{\mathcal{T}_A})$ is of infinite type and it is determined uniquely by the fact that $(\psi(q_x))_{x\in\mathcal{G}}$ is the unique nonnegative normalized (in the sense that $\sum_{x\in\mathcal{G}} N(x)^{-\beta} v_x = 1$) eigenvector $v$ of the matrix $AN^{-\beta_c}$ with (dominant) eigenvalue 1.*

(iii) *For $\beta < \beta_c$ there are no KMS$_\beta$ states on $\mathcal{T}_A$ at all.*

*Proof.* As observed above $\#\Sigma_A = d(A)$ and hence by 11.1.i one also has that $\#\widetilde{\Omega}_e = d(A)$. By the remark following 8.6 we have $\Omega_e = \widetilde{\Omega}_e \setminus \{\epsilon\}$. Since $A$ is irreducible and hence 10.1.(COL) holds we have that $\epsilon \notin \widetilde{\Omega}_{\mathcal{T}_A} \supseteq \widetilde{\Omega}_e$ so actually $\Omega_e = \widetilde{\Omega}_e$. Therefore $\#\Omega_e = d(A)$.

Given $\beta > \beta_c$ the set of (positive) measures $\gamma$ on $\Omega_e$ with $Z(\beta, \gamma) = 1$ therefore forms a simplex of dimension $d(A) - 1$. Point (i) then follows from 9.7.

As for (iii), this follows from 15.2 given that 10.1.(IRR+INF) are granted.

In order to prove (ii) observe that by 18.3 one has that $\beta_c \notin \dot{I}_c$. From 10.6 we then conclude that all $\beta_c$-scaling states are of infinite type. Using 13.6 we therefore have that the KMS$_{\beta_c}$ states on $\mathcal{T}_A$ correspond bijectively to the normalized (in the above sense) nonnegative fixed points for the matrix $AN^{-\beta_c}$. By 18.3.iii and the Perron–Frobenius Theorem [**Se**: Theorem 1.5] we have that there is exactly one such vector. This concludes the proof. □

Our next result gives a precise description of the KMS states on $\mathcal{O}_A$ in terms of the eigenvalues of $AN^{-\beta}$, even if $A$ is not irreducible. This was first proved in [**EFW**] under the special case that the $N(x)$ are all the same.

**18.5. Theorem.** *Under 8.1 let $A$ be a finite matrix without identically zero columns. Then the KMS$_\beta$ states on $\mathcal{O}_A$ occur exactly at the values of $\beta$ for which there exists a nonnegative vector $v \neq 0$ satisfying $AN^{-\beta}(v) = v$. Given such a $\beta$ the correspondence*

$$\psi \mapsto (\psi(q_x))_{x\in\mathcal{G}}$$

*defines an affine bijection from the simplex of all KMS$_\beta$ states $\psi$ on $\mathcal{O}_A$ to the simplex of all normalized (i.e. $\sum_{x\in\mathcal{G}} N(x)^{-\beta} v_x = 1$) nonnegative solutions of the equation $AN^{-\beta}(v) = v$.*

*Proof.* By 17.2 we know that the KMS states on $\mathcal{O}_A$, equivalently the KMS states on $\mathcal{T}_A$ which factor through $\mathcal{O}_A$, correspond to the $\beta$-scaling measures $\lambda$ on $\Omega_{\mathcal{T}_A}$ supported in the *closure* of $\Omega_\infty$. Under the present hypothesis that $\mathcal{G}$ is finite we claim that $\Omega_\infty$ is closed in $\Omega_{\mathcal{T}_A}$. To see this let $\xi$ be a bounded element of $\Omega_{\mathcal{T}_A}$ with stem $\omega$. Then the set

$$V = \{\eta \in \Omega_{\mathcal{T}_A} : \omega \in \eta \text{ and } \omega z \notin \eta \text{ for all } z \in \mathcal{G}\}$$

is a neighborhood of $\xi$ not intersecting $\Omega_\infty$.

Therefore the measures $\lambda$ mentioned above consist precisely of the infinite type $\beta$-scaling measures. The conclusion then follows from 13.6. □

PARTIAL DYNAMICAL SYSTEMS AND THE KMS CONDITION

BY RUY EXEL AND MARCELO LACA

JOB PAGE

*** PLEASE DISCARD THIS PAGE AND REPLACE IT ***

*** WITH THE LAST FOUR PRINTED PAGES ***